\RequirePackage{fix-cm}
\documentclass[smallextended]{svjour3}       % onecolumn (second format)
\smartqed  % flush right qed marks, e.g. at end of proof
\usepackage{graphicx}
\usepackage{amsmath}
\usepackage[numbers]{natbib}
\usepackage[psamsfonts]{amsfonts} 
\usepackage{pdflscape}
\usepackage{dsfont}    
\usepackage{color}    
\usepackage[usenames,dvipsnames]{xcolor}

\usepackage{hyperref}
\hypersetup{pdfpagemode=FullScreen,  
colorlinks=true,
citecolor=Bittersweet,
linkcolor=Bittersweet,
urlcolor=Bittersweet
}

\newcommand{\mockalph}[1]{}

\newcommand{\cqfd}{\hfill $\square$}

\newcommand{\R}{\mathbb R}

\newcommand{\n}{^{(n)}}

\newcommand{\Xb}{\mathbf{X}}
\newcommand{\Sb}{\mathbf{S}}

\newcommand{\Tb}{\ensuremath{\mathbf{T}}}
\newcommand{\vb}{\ensuremath{\mathbf{v}}}
\newcommand{\xb}{\ensuremath{\mathbf{x}}}

\newcommand{\Ab}{\ensuremath{\mathbf{A}}}
\newcommand{\Bb}{\ensuremath{\mathbf{B}}}

\newcommand{\Ob}{\ensuremath{\mathbf{O}}}

\newcommand{\Cb}{\ensuremath{\mathbf{C}}}

\newcommand{\Ub}{\ensuremath{\mathbf{U}}}

\newcommand{\Yb}{\ensuremath{\mathbf{Y}}}

\newcommand{\thetab}{{\pmb \theta}}

\newcommand{\Thetab}{{\pmb \Theta}}
\newcommand{\Psib}{{\pmb \Psi}}

\newcommand{\Deltab}{{\pmb \Delta}}
\newcommand{\taub}{{\pmb \tau}}

\newcommand{\Gamb}{{\pmb \Gamma}}

\newcommand{\pr}{^{\prime}}

\newcommand{\ny}{n\rightarrow\infty}

\begin{document}

\title{Detecting the direction of a signal on high-dimensional spheres\thanks{Research is supported by the Program of Concerted Research Actions (ARC) of the Universit\'{e} libre de Bruxelles, by a research fellowship from the Francqui Foundation, and by the Cr\'{e}dit de Recherche J.0134.18 of the FNRS (Fonds National pour la Recherche Scientifique), Communaut\'{e} Fran\c{c}aise de Belgique.}
}
\subtitle{Non-null and Le Cam optimality results}


\author{Davy Paindaveine        \and
        Thomas Verdebout
}


\institute{Davy Paindaveine \at Universit\'{e}  libre de Bruxelles,
              ECARES and D\'{e}partement de Math\'{e}matique,
              Avenue \mbox{F.D.} Roosevelt, 50,
              ECARES, CP114/04,
              B-1050, Brussels,
              Belgium\\
              Tel.: +3226503845\\
              \email{dpaindav@ulb.ac.be}           %  \\
           \and
           Thomas Verdebout \at Universit\'{e}  libre de Bruxelles,
              ECARES \mbox{F.D.} D\'{e}partement de Math\'{e}matique,
              Avenue F.D. Roosevelt, 50,
              ECARES, CP114/04,
              B-1050, Brussels,
              Belgium\\
              Tel.: +3226505892\\
              \email{tverdebo@ulb.ac.be}           %  \\
}

\date{Received: date / Accepted: date}

\maketitle

\begin{abstract}
We consider one of the most important problems in directional statistics, namely the problem of testing the null hypothesis that the spike direction~$\thetab$ of a Fisher--von Mises--Langevin distribution on the $p$-dimen\-sional unit hypersphere is equal to a given direction~$\thetab_0$. After a reduction through invariance arguments, we derive local asymptotic normality (LAN) results in a general high-dimensional framework where the dimension~$p_n$ goes to infinity at an arbitrary rate with the sample size~$n$, and where the concentration~$\kappa_n$ behaves in a completely free way with~$n$, which offers a spectrum of problems ranging from arbitrarily easy to arbitrarily challenging ones. We identify various asymptotic regimes, depending on the convergence/divergence properties of~$(\kappa_n)$, that yield different contiguity rates and different limiting experiments. In each regime, we derive Le Cam optimal tests under specified~$\kappa_n$ and we compute, from the Le Cam third lemma, asymptotic powers of the classical Watson test under contiguous alternatives. We further establish LAN results with respect to both spike direction and concentration, which allows us to discuss optimality also under unspecified~$\kappa_n$. To investigate the non-null behavior of the Watson test outside the parametric framework above, we derive its local asymptotic powers through martingale CLTs in the broader, semiparametric, model of rotationally symmetric distributions. A Monte Carlo study shows that the finite-sample behaviors of the various tests remarkably agree with our asymptotic results.  
\end{abstract}

%%%%%%%%%%%%%%%%%%%%%%%%%%%%%%%%%%%%%%%%%%%%%%%%%%%%%%%%%%%%%%%%%%

\section{Introduction}

In directional statistics, the sample space is the unit sphere~$\mathcal{S}^{p-1}=\{\xb\in\R^p:\|\xb\|^2=\xb'\xb=1\}$ in~$\R^p$. By far the most classical distributions on~$\mathcal{S}^{p-1}$ are the Fisher--von Mises--Langevin (FvML) ones; see, e.g.,  \cite{LV17book} or \cite{MarJup2000}. We say that the random vector~$\Xb$ with values in~$\mathcal{S}^{p-1}$ has an ${\rm FvML}_p(\thetab,\kappa)$ distribution, with~$\thetab\in\mathcal{S}^{p-1}$ and~$\kappa\in(0,\infty)$, if it admits the density (throughout, densities on the unit sphere are with respect to the surface area measure)
\begin{equation}
\label{FvMLdensity}	
\xb
\mapsto 
\frac{c_{p,\kappa}}{\omega_{p-1}}
\,
{\rm exp}(\kappa\hspace{.2mm} \xb\pr \thetab) 
,
\end{equation} 
where, denoting as~$\Gamma(\cdot)$ the Euler Gamma function and as~$\mathcal{I}_\nu(\cdot)$ the order-$\nu$ modified Bessel function of the first kind,
$
\omega_{p}
:=
(2\pi^{p/2})/\Gamma\big(\frac{p}{2}\big)
$
is the surface area of~$\mathcal{S}^{p-1}$ 
and 
$$
c_{p,\kappa}
:=
1\, \Big/ \int_{-1}^1 (1-t^2)^{(p-3)/2}  \exp(\kappa t)\, dt
=
\frac{(\kappa/2)^{(p/2)-1}}{\sqrt{\pi}\,\Gamma\big(\frac{p-1}{2}\big) \mathcal{I}_{\frac{p}{2}-1}(\kappa)}
\cdot 
$$
Clearly, $\thetab$ is a location parameter ($\thetab$ is the modal location on the sphere), that identifies the spike direction of the hyperspherical signal. In contrast,~$\kappa$ is a scale or \emph{concentration} parameter: the larger~$\kappa$, the more concentrated the distribution is about the modal location~$\thetab$. As~$\kappa$ converges to zero, $c_{p,\kappa}$ converges to~$c_p:=\Gamma\big(\frac{p}{2}\big)/(\sqrt{\pi}\,\Gamma\big(\frac{p-1}{2}\big))$ and the density in~(\ref{FvMLdensity}) converges to the density~${\bf x} \mapsto 1/\omega_{p}$ of the uniform distribution over~$\mathcal{S}^{p-1}$. The other extreme case, obtained for arbitrarily large values of~$\kappa$, provides distributions that converge to a point mass in~$\thetab$. Of course, it is expected that the larger~$\kappa$, the easier it is to conduct inference on~$\thetab$ --- that is, the more powerful the tests on~$\thetab$ and the smaller the corresponding confidence zones.
 
In this paper, we consider inference on~$\thetab$ and focus on the generic  testing problem for which the null hypothesis~$\mathcal{H}_0:\thetab=\thetab_0$, for a fixed~$\thetab_0\in\mathcal{S}^{p-1}$, is to be tested against~$\mathcal{H}_1:\thetab\neq\thetab_0$ on the basis of a random sample~$\Xb_{n1},\ldots,\Xb_{nn}$ from the ${\rm FvML}_p(\thetab,\kappa)$ distribution --- the triangular array notation antici\-pates non-standard setups where~$p$ (hence, also~$\thetab$) and/or~$\kappa$ will depend on~$n$. Inference problems on $\thetab$ in the low-dimensional case have been considered among others in \cite{ChaHe1993}, \cite{ChaRi01}, \cite{GutLoc1988}, \cite{HeSim1992}, \cite{La02}, \cite{Leyetal2013}, \cite{PaiVer17b} and \cite{Wat1983b}. The related spherical regression problem has been tackled in  \cite{Do03} and \cite{Ri89}, while testing for location on axial frames has been considered in \cite{AJu13}.  

Letting~$\bar{\Xb}_n:=\frac{1}{n}\sum_{i=1}^n \Xb_{ni}$, the most classical test for the testing problem above is the Watson \cite{Wat1983b} test rejecting the null at asymptotic level~$\alpha$ whenever 
\begin{equation}
\label{WatsonLD}	
W_{n}
:=
\frac{n(p-1)\bar{\Xb}_n\pr (\mathbf{I}_{p}-\thetab_0\thetab_0\pr) \bar{\Xb}_n}{1-\frac{1}{n} \sum_{i=1}^n (\Xb_{ni}\pr\thetab_0)^2} 
%=
%\frac{n(p_n-1)\bar{\Xb}_{n}\pr (\mathbf{I}_{p}-\thetab_0\thetab_0\pr) \bar{\Xb}_{n}}{f_{n2}} 
%+
%o_{\rm P}(1)
>
\chi^2_{p-1,1-\alpha}
, 
\end{equation} 
where~$\mathbf{I}_{\ell}$ stands for the $\ell$-dimensional identity matrix and~$\chi^2_{\ell,1-\alpha}$ denotes the $\alpha$-upper quantile of the chi-square distribution with~$\ell$ degrees of freedom. In the classical setup where~$p$ and~$\kappa$ are fixed, the asymptotic properties of the Watson test are well-known, both under the null and under local alternatives; see, e.g.,  \cite{MarJup2000} or \cite{Wat1983b}. Optimality properties in the Le Cam sense have been studied in \cite{PaiVer2015a}. In the non-standard setup where~$\kappa=\kappa_n$ converges to zero, \cite{PaiVer17b} investigated the asymptotic null and non-null behaviors of the Watson test. Interestingly, irrespective of the rate at which~$\kappa_n$ converges to zero (that is, irrespective of how fast the inference problem becomes more challenging as a function of~$n$), the Watson test keeps meeting the asymptotic nominal level constraint and maintains strong optimality properties; see \cite{PaiVer17b} for details. {In the other non-standard, high-concentration, setup where~$\kappa_n$ diverges to infinity, \cite{PaiVer19} showed the Watson test also enjoys strong optimality properties.}  
    
For a fixed dimension~$p$, this essentially settles the investigation of the properties of the Watson test and the study of the corresponding hypothesis testing problem. Nowadays, however, increasingly many applications lead to considering high-dimensional directional data: tests of uniformity on high-dimensional spheres have been studied in \cite{Caietal2013}, \cite{Chi1991}, \cite{Cueetal2009} and \cite{PaiVer17a}, while high-dimensional FvML distributions (or mixtures of high-dimensional FvML distributions) have been considered in magnetic resonance, gene-expression, and text mining; see, among others, \cite{banerjee2003generative}, \cite{Banetal2005} and \cite{Dry2005}. This motivates considering the high-dimensional spherical location problem, based on a random sample~$\Xb_{n1},\ldots,\Xb_{nn}$ from the ${\rm FvML}_{p_n}(\thetab_n,\kappa_n)$ distribution, with~($p_n$) diverging to infinity (the dimension of~$\thetab_n$ then depends on~$n$, which justifies the notation). In this context, it was proved in \cite{Leyetal2015} that the Watson test is robust to high-dimensionality in the sense that, as~$p_n$ goes to infinity
\vspace{-.4mm}
 with~$n$, this test still has asymptotic size~$\alpha$ under~$\mathcal{H}_0\n:\thetab_n=\thetab_{n0}$. This does not require any condition on the concentration sequence~$(\kappa_n)$ nor on the rate at which~$p_n$ goes to infinity, hence covers arbitrarily easy problems ($\kappa_n$ large) and arbitrarily challenging ones ($\kappa_n$ small), as well as moderately high dimensions and ultra-high dimensions. On its own, however, this null robustness result is obviously far from sufficient to motivate using the Watson test in high dimensions, as it might very well be that robustness under the null is obtained at the expense of power (in the extreme case, the Watson test, in high dimensions, might actually asymptotically behave like the trivial $\alpha$-level test that randomly rejects the null with probability~$\alpha$). 

These considerations raise many interesting questions, among which:
are there alternatives under which the Watson test is consistent in high dimensions? What are the less severe alternatives (if any) under which the Watson test exhibits non-trivial asymptotic powers? Is the Watson test rate-optimal or, on the contrary, are there tests that show asymptotic powers under less severe alternatives than those detected by the Watson test? Does the Watson test enjoy optimality properties in high dimensions? As we will show, answering these questions will require considering several regimes fixing how the concentration~$\kappa_n$ behaves as a function of the dimension~$p_n$ and sample size~$n$. Our results, that will crucially depend on the regime considered, are extensive in the sense that they answer the questions above in all possible regimes.

Our results will rely on two different approaches. (a) The first approach is based on Le Cam's asymptotic theory of statistical experiments. While this theory is very general, it does not directly apply in the present context since the high-dimensional spherical location problem involves a parametric space, namely~$\{(\thetab_n,\kappa):\thetab_n\in\mathcal{S}^{p_n-1}, \kappa\in(0,\infty)\}$, that depends on~$n$ (through~$p_n$). We solve this by exploiting the invariance properties of the testing problem considered. In the image of the model by the corresponding maximal invariant, indeed, the parametric space does not depend on~$n$ anymore, which opens the door to studying the problem through the Le Cam approach. We derive stochastic second-order expansions of the resulting log-likelihood ratios, which is the main technical ingredient to establish the local asymptotic normality (LAN) of the invariant model. The LAN property takes different forms and involves different contiguity rates depending on the regime that is considered. In each regime, we determine the Le Cam optimal test for the problem considered and apply the Le Cam third lemma to obtain the asymptotic powers of this test and of the Watson test. This allows us to determine the regime(s) in which the Watson test is Le Cam optimal, or only rate-optimal, or not even rate-optimal. While this is first done under specified concentration~$\kappa_n$, we further provide LAN results with respect to both location and concentration to be able to discuss optimality under unspecified~$\kappa_n$.  
(b) %In regimes where the Watson test is not rate-optimal, that is, in regimes where it is blind to contiguous alternatives, this first, Le Cam, approach leaves open the question of the existence of alternatives that can be detected by the Watson test. 
While our investigation in~(a) will fully characterize the asymptotic optimality properties of the Watson test in the FvML case, it will not provide any insight on the non-null behavior of this test outside this stringent parametric framework. This motivates complementing our investigation by a second approach, based on martingale CLTs. We will consider a broad, semiparametric, model, namely the class of rotationally symmetric distributions, and will identify the alternatives (if any) under which the Watson test will show non-trivial asymptotic powers in high dimensions. Again, this requires considering various regimes according to the concentration pattern.   

The outline of the paper is as follows. In Section~\ref{secLAN}, we consider the high-dimensional version of the FvML spherical location problem. In Section~\ref{subsecinvariance}, we describe the invariance approach that allows us to later rely on Le Cam's asymptotic theory of statistical experiments. In Section~\ref{subseckappaspec}, we provide a stochastic second-order expansion of the resulting invariant log-likelihood ratios and prove, in various regimes that we identify, that these invariant models are locally asymptotically normal. This allows us to derive the corresponding Le Cam optimal tests for the specified concentration problem and to study the non-null asymptotic behavior of the Watson test in the light of these results. In Section~\ref{subseckappaunspec}, we tackle the unspecified concentration problem through the derivation of LAN results that are with respect to both location and concentration. In Section~\ref{secBilling}, we conduct a systematic investigation of the non-null asymptotic properties of the Watson test in the broader context of rotationally symmetric distributions. This in particular confirms the FvML results obtained in Section~\ref{secLAN}. In Section~\ref{secSimu}, we conduct a Monte Carlo study to investigate how well the finite-sample behaviors of the various tests reflect our theoretical asymptotic results. In Section~\ref{secFinaComments}, we summarize the results obtained in the paper and shortly discuss research perspectives. Finally, an appendix contains all proofs.

%%%%%%%%%%%%%%%%%%%%%%%%%%%%%%%%%%%%%%%%%%%%%%%%%%%%%

\section{Invariance and Le Cam optimality} 
\label{secLAN} 
 
As already mentioned in the introduction, the high-dimensional spherical location problem requires considering triangular arrays of observations of the form~$\Xb_{ni}$, $i=1,\ldots,n$, $n=1,2,\ldots$ For any sequence~$(\thetab_n)$ 
\vspace{-.6mm}
 such that~$\thetab_n$ belongs to~$\mathcal{S}^{p_n-1}$ for any~$n$ and any sequence~$(\kappa_n)$ in~$(0,\infty)$, we denote as~${\rm P}\n_{\thetab_n,\kappa_n}$ the hypothesis under which~$\Xb_{ni}$, $i=1,\ldots,n$, form a random sample from the ${\rm FvML}_{p_n}(\thetab_n,\kappa_n)$ distribution. The resulting sequence of statistical models is then associated with
\begin{equation}
	\label{FvMLstatmodels}
	\mathcal{P}\n
=
\Big\{
{\rm P}\n_{\thetab_n,\kappa}: (\thetab_n,\kappa)\in \Thetab_n:=\mathcal{S}^{p_n-1}\times (0,\infty) 
\Big\}
\end{equation}
(the index in the parameter~$\thetab_n$ in principle is superfluous but is used here to stress the dependence of this  parameter on~$p_n$, hence on~$n$).
\vspace{-.5mm}
 The spherical 
 location problem consists in testing
 \vspace{-.7mm}
 the null hypothesis~$\mathcal{H}_0\n: \thetab_n=\thetab_{n0}$ against the alternative~$\mathcal{H}_0\n: \thetab_n\neq\thetab_{n0}$, where~$(\thetab_{0n})$ is a fixed sequence such that~$\thetab_{0n}$ belongs to~$\mathcal{S}^{p_n-1}$ for any~$n$. Clearly,~$\thetab_n$ is the parameter of interest, whereas~$\kappa_n$ plays the role of a nuisance. 
Our main objective in this section is to derive Le Cam optimality \vspace{-.2mm}
 results for this problem, referring to sequences of local alternatives of the form~${\rm P}\n_{\thetab_n,\kappa_n}$, with~$\thetab_{n}=\thetab_{n0}+\nu_n \taub_n$, where the sequence~$(\nu_n)$ and the bounded sequence~$(\taub_{n})$, respectively in~$(0,\infty)$ and~$\R^{p_n}$, are such that~$\thetab_{n}\in\mathcal{S}^{p_n-1}$ for any~$n$, which imposes that
\begin{equation}
\label{constraintnun}
	\thetab_{n0}\pr\taub_n
	=
	-\frac{1}{2} \nu_n \|\taub_n\|^2 
\end{equation}
for any~$n$; {throughout, ``the sequence~$(\taub_n)$ in~$\R^{p_n}$ is bounded" means that~$\taub_n\in\R^{p_n}$ for any~$n$ and that~$\|\taub_n\|=O(1)$ as~$n\to\infty$}. The obvious lack of identifiability of~$\nu_n$ and $\taub_n$ will be no problem in the sequel  (only the locally perturbed parameter values~$\thetab_n=\thetab_0 +\nu_n \taub_n$ are of interest, hence not the individual quantities~$\nu_n$ and~$\taub_n$ themselves) and this form of local alternatives is actually the standard one in the Le Cam theory; see, e.g., Chapter~6 in \cite{LeCam2000} or Definition~7.14 in~\cite{van2000}. Whenever local asymptotic powers will be considered below, we will assume that~$\|\taub_n\|$ is~$O(1)$ without being~$o(1)$, so that~$\nu_n$ will characterize (the rate of) the severity of the local alternatives~$\thetab_n=\thetab_0 +\nu_n \taub_n$ (the slower~$\nu_n$ goes to zero, the more severe the corresponding local alternatives). 

Since the sequence of ``statistical experiments" 
%(in Le Cam's terminology) 
associated with~(\ref{FvMLstatmodels}) involves parametric spaces~$\Thetab_n$ that depend on~$n$, applying Le Cam's theory will require the following reduction of the problem through invariance arguments.

%%%%%%

\subsection{Reduction through invariance}
\label{subsecinvariance}

Denoting as~$SO_p(\thetab)$ 
\vspace{-.5mm}
 the collection of $p \times p$ orthogonal matrices satisfying~$\Ob\thetab=\thetab$, the null hypothesis~$\mathcal{H}_0\n$ is invariant under the group~$\mathcal{G}_{n},_{\,\circ}$ collecting the transformations 
%$=\{g\n_{\Ob}:\Ob\in SO_{p_n}(\thetab_{n0}) \}$ of transformations defined through
$$
(\Xb_{n1},\ldots,\Xb_{nn})
\mapsto 
g_{n\Ob}(\Xb_{n1},\ldots,\Xb_{nn}) 
=
(\Ob\Xb_{n1},\ldots,\Ob\Xb_{nn}) 
,
%\hspace{2mm}
%\Ob\in SO_{p_n}(\thetab_{n0})
%.
$$
with~$\Ob\in SO_{p_n}(\thetab_{n0})$. The transformation~$g_{n\Ob}$ induces a transformation of the parametric space~$\Thetab_n$ defined through~$%\bar{g}\n_{\Ob}
(\thetab_n,\kappa)\mapsto (\Ob\thetab_n,\kappa)$. The orbits of the resulting induced group are~$\mathcal{C}_{u,\kappa}(\thetab_{n0}):=\{ \thetab_n \in \mathcal{S}^{p_n-1} : \thetab_n'\thetab_{n0}=u\}\times\{\kappa\}$, with~$u\in[-1,1]$ and~$\kappa\in(0,\infty)$. 
%, where~$\mathcal{C}_{\thetab_{n0},u}$ is the hyperspherical parallel~$$.  
In such a context, the invariance principle (see, e.g., \cite{Lehetal2005}, Chapter~6) leads to restricting to tests~$\phi_n$ that are invariant with respect to the group~$\mathcal{G}_{n},_{\circ}$. 
%, in the sense that~$\phi\n\circ g\n_{\Ob}=\phi\n$ for any~$\Ob\in SO_{p_n}(\thetab_{n0})$. %This ensures that~$\phi\n$ rejects the null on the basis of~$(\Xb_{n1},\ldots,\Xb_{nn})$ if and only if~$\phi\n$ on any $g\n_{\Ob}$-rotated version of~$(\Xb_{n1},\ldots,\Xb_{nn})$. 
Denoting as~$\Tb_n=\Tb_n(\Xb_{n1},\ldots,\Xb_{nn})$ a maximal invariant statistic for~$\mathcal{G}_{n},_{\circ}$, the class of invariant tests coincides with the class of~$\Tb_n$-measurable tests. Invariant tests thus are to be defined in the image
\begin{equation}
	\label{InvariantFvMLstatmodels}
	\mathcal{P}^{(n)\Tb_n}
=
\Big\{
{\rm P}^{(n)\Tb_n}_{u,\kappa}: (u,\kappa)\in \Psib:=[-1,1]\times (0,\infty)
\Big\}
\end{equation}
of the model~$\mathcal{P}\n$ by~$\Tb_n$,
\vspace{-.8mm}
  where~${\rm P}^{(n)\Tb_n}_{u,\kappa}$ denotes the common distribution of~$\Tb_n$ under any~${\rm P}^{(n)}_{\thetab_n,\kappa}$ with~$(\thetab_n,\kappa)\in \mathcal{C}_{u,\kappa}(\thetab_{n0})$. Unlike the original sequence of statistical experiments in~(\ref{FvMLstatmodels}), the invariant one in~(\ref{InvariantFvMLstatmodels}) involves a fixed parametric space~$\Psib$, which makes it in principle possible to rely on Le Cam's asymptotic theory. 

Now, the original local log-likelihood ratios
$\log(d{\rm P}^{(n)}_{\thetab_n,\kappa_n}/d{\rm P}^{(n)}_{\thetab_{n0},\kappa_n})$ associated with the generic local alternatives~$\thetab_{n}=\thetab_{n0}+\nu_n \taub_n$ above correspond, in view of~(\ref{constraintnun}), to the invariant local log-likelihood ratios 
\begin{equation}
\label{ilr}	
\Lambda^{(n){\rm inv}}_{\thetab_n/\thetab_{n0};\kappa_n}
:=
\log \frac{d{\rm P}^{(n)\Tb_n}_{1-\nu_n^2 \|\taub_n\|^2/2,\kappa_n}}{d{\rm P}^{(n)\Tb_n}_{1,\kappa_n}}
\cdot
\end{equation}
Deriving local asymptotic normality (LAN) results requires investigating the asymptotic behavior of such invariant log-likelihood ratios, which in turn requires evaluating the corresponding likelihoods. While obtaining a closed-form expression for~$\Tb_n$ and its distribution is a very challenging task, these likelihoods can be obtained from Lemma~2.5.1 in~\cite{Gir1996}, which, denoting as~$m_n$ the surface area measure on~$\mathcal{S}^{p_n-1}\times\ldots\times \mathcal{S}^{p_n-1}$ ($n$ times), yields 
\begin{eqnarray}
 \frac{d{\rm P}^{(n)\Tb_n}_{1-\nu_n^2 \|\taub_n\|^2/2,\kappa_n}}{d{m_n}} 
&=&
\int_{SO_{p_n}(\thetab_{n0})}
 \frac{d{\rm P}^{(n)}_{\thetab_{n},\kappa_n}}{d{m_n}}
 (\Ob\Xb_{n1},\ldots,\Ob\Xb_{nn}) 
\,
d\Ob
\nonumber
\\[2mm]
&= & 
\int_{SO_{p_n}(\thetab_{n0})}
\prod_{i=1}^n 
\bigg(
\frac{c_{p_n,\kappa_n}}{\omega_{p_n-1}}
\exp\big( \kappa_n (\Ob\Xb_{ni})\pr \thetab_{n} \big) 
\bigg)
\,d\Ob
\nonumber
\\[2mm]
&= & 
\frac{c_{p_n,\kappa_n}^n}{\omega_{p_n-1}^n}
\,
\int_{SO_{p_n}(\thetab_{n0})}
\exp\big( n\kappa_n \bar{\Xb}_{n} \Ob\pr \thetab_{n} \big) 
\,d\Ob
,
\label{tsj}
\end{eqnarray}
where integration is with respect to the Haar measure on $SO_{p_n}(\thetab_{n0})$. Note that~(\ref{tsj}) shows that
\vspace{-1.1mm}
  the invariant null probability measure~${\rm P}^{(n)\Tb_n}_{1,\kappa_n}$ coincides with the original  null probability measure~${\rm P}^{(n)}_{\thetab_{n0},\kappa_n}$. In other words, it is only for non-null probability measures that the invariance reduction above is non-trivial.

%%%%%%%%

\subsection{Optimal testing under specified~$\kappa_n$}
\label{subseckappaspec}

The main ingredient needed to obtain LAN results is Theorem~\ref{TheoLAQ} below, that provides a stochastic second-order expansion of the invariant log-likelihood ratios in~(\ref{ilr}). 
To state this theorem, we need to introduce the following notation. We will refer to the decomposition 
%\begin{equation}
%	\label{pdfUFvML}
$\Xb_{ni}
=
U_{ni} \thetab_{n0}
+
V_{ni} \Sb_n
$,
%\end{equation}
with
$$
U_{ni}
%=U_{ni}(\thetab_n)
=\Xb_{ni}'\thetab_{n0}
,
\ \
V_{ni}
%=V_{ni}(\thetab_n)
=(1-U_{ni}^2)^{1/2}
\ \
\textrm{ and }
\ \
\Sb_{ni}
=
\frac{(\mathbf{I}_{p_n}-\thetab_{n0}\thetab_{n0}\pr)\Xb_{ni}}{\|(\mathbf{I}_{p_n}-\thetab_{n0}\thetab_{n0}\pr)\Xb_{ni}\|}
,
$$
as the tangent-normal decomposition of~$\Xb_{ni}$ with respect to~$\thetab_{n0}$. Under the hypothesis~${\rm P}\n_{\thetab_{n0},\kappa_n}$,
% (equivalently, under~${\rm P}^{(n)\Tb_n}_{1,\kappa_n}$), 
$U_{ni}$ has probability density function
\begin{equation}
\label{pdfUFvML}	
u\mapsto c_{p_n,\kappa_n} (1-u^2)^{(p_n-3)/2}  \exp(\kappa_n u) \, \mathbb{I}[u\in[-1,1]]
,
\end{equation}
where~$\mathbb{I}[A]$ denotes the indicator function of~$A$, $\Sb_{ni}$ is uniformly distributed over the ``equator"~$\{\xb\in\mathcal{S}^{p_n-1}:\xb'\thetab_{n0}=0\}$, and~$U_{ni}$ and~$\Sb_{ni}$ are mutually independent. Throughout, 
\vspace{-.6mm}
we will denote as~$e_{n\ell}={\rm E}[U_{ni}^\ell]$ and~$\tilde{e}_{n\ell}={\rm E}[(U_{ni}-e_{n1})^\ell]$, $\ell=1,2,\ldots$ the non-central and 
\vspace{-.6mm}
central moments of~$U_{ni}$ under~${\rm P}\n_{\thetab_{n0},\kappa_n}$, and as~$f_{n\ell}={\rm E}[V_{ni}^{\ell}]$ the corresponding non-central moments of~$V_{ni}$. Although this is not stressed in the notation, these moments clearly depend on~$p_n$ and~$\kappa_n$; for instance,
\begin{equation}
	\label{linken1en2tildekappan}
e_{n1}
=
\frac{\mathcal{I}_{\frac{p_n}{2}}(\kappa_n)}{\mathcal{I}_{\frac{p_n}{2}-1}(\kappa_n)}
,
\ \ 
\tilde{e}_{n2}
= 
1-\frac{p_n-1}{\kappa_n}\, e_{n1}-e_{n1}^2
\ \textrm{ and } \
f_{n2}
=
\frac{p_n-1}{\kappa_n}\, e_{n1}
\end{equation}
(this readily follows from~(2)--(3) in \cite{Sch1978} by using the standard properties of exponential families; see also Lemma~S.2.1 in \cite{PaiVer17asupp}).
 We can now state the stochastic second-order expansion result of the invariant log-likelihood ratios in~(\ref{ilr}).

\begin{theorem}
\label{TheoLAQ}
Let~$(p_n)$ be a sequence of integers that diverges to infinity and~$(\kappa_n)$ be an arbitrary sequence in~$(0,\infty)$. Let~$(\thetab_{n0})$, $(\nu_n)$ and~$(\taub_n)$ be sequences such that $\thetab_{n0}$ and~$\thetab_n=\thetab_{n0}+\nu_n\taub_n$ belong to~$\mathcal{S}^{p_n-1}$ for any~$n$, with~$(\taub_n)$ bounded and~$(\nu_n)$ such that 
\begin{equation}
\label{taillenu}
\nu_n^2 
=
O\Big(\frac{\sqrt{p_n}}{n \kappa_n e_{n1}}\Big)
.
\end{equation}
Then, letting
$$
Z_n
:=
\frac{\sqrt{n}(\bar{\Xb}_{n}\pr \thetab_0-e_{n1})}{\sqrt{\tilde{e}_{n2}}}
\quad\textrm{and}\quad%
\widetilde{W}_n
:=
\frac{W_{n}-(p_n-1)}{\sqrt{2(p_n-1)}}
,
$$
we have that
% $($equivalently, under~${\rm P}^{(n)\Tb_n}_{1,\kappa_n})$,
	\begin{eqnarray*}
\lefteqn{	
\hspace{-3mm} 
\Lambda^{(n){\rm inv}}_{\thetab_n/\thetab_{n0};\kappa_n}
%=
%\log 
%\Bigg( 
%\frac{d{\rm P}^{(n)\Tb_n}_{1-\frac{1}{2} \nu_n^2 \|\taub_n\|^2,\kappa_n}}{d{\rm P}^{(n)\Tb_n}_{1,\kappa_n}}
%\Bigg)
\!\!=
-
\frac{1}{2} \sqrt{n} \kappa_n \nu_n^2 \sqrt{\tilde{e}_{n2}}\, \|\taub_n\|^2  Z_n 
+ 
\frac{n \kappa_n\nu_n^2 e_{n1}}{\sqrt{2}p_n^{1/2}}
\|\taub_n\|^2 \Big(1- \frac{1}{4} \nu_n^2 \|\taub_n\|^2\Big)
  \widetilde{W}_n
}
\\[2mm]
& &  
\hspace{10mm} 
-
\frac{1}{8}
n \kappa_n\nu_n^4 e_{n1} \|\taub_n\|^4
-
\frac{n^2 \kappa_n^2 \nu_n^4 e_{n1}^2}{4p_n}  \|\taub_n\|^4 \Big(1- \frac{1}{4} \nu_n^2 \|\taub_n\| ^2\Big)^2 
  +o_{\rm P}(1)
,
\end{eqnarray*}
as~$n\to\infty$ under~${\rm P}^{(n)}_{\thetab_{n0},\kappa_n}$.
\end{theorem}
\vspace{1mm}

Recalling that
\vspace{-.4mm}
 the log-likelihood ratio~$\Lambda^{(n){\rm inv}}_{\thetab_n/\thetab_{n0};\kappa_n}$ refers to the local perturbation~$\thetab_n'\thetab_{n0}=1-\nu_n^2 \|\taub_n\|^2/2$ of the null reference value~$\thetab_{n0}'\thetab_{n0}=1$, the result in Theorem~\ref{TheoLAQ} essentially shows that the invariant model considered enjoys a \emph{local asymptotic quadraticity} (LAQ) structure
\vspace{-.5mm}
  in the vicinity of the null hypothesis~$\mathcal{H}_0\n:\thetab_n=\thetab_{n0}$; see, e.g., \cite{LeCam2000}, page~120. Actually, quadraticity, which is supposed to be in the increment~$-\nu_n^2 \|\taub_n\|^2/2$, only holds for arbitrarily small values of this increment, hence only in regimes where~$\nu_n$ will converge to zero (in regimes below where, in contrast,~$\nu_n$ will be constant, the non-flat manifold structure of the hypersphere actually prevents a standard quadraticity property). This LAQ result
  %, which is the main ingredient required to establish the \emph{local asymptotic normality} (LAN) of the model, 
  hints that optimal testing for the specified-$\kappa_n$ problem at hand is obtained by rejecting the null for small values of~$Z_n$ (that is, when~$\bar{\Xb}_n$ and~$\thetab_{n0}$ project far from each other onto the axis~$\pm\thetab_{n0}$), for large values of~$\widetilde{W}_n$ (that is, when~$\bar{\Xb}_n$ and~$\thetab_{n0}$ project far from each other onto the orthogonal complement to~$\thetab_{n0}$ in~$\R^{p_n}$), or, more generally, for large values of a hybrid test statistic of the form
$$
Q_n^{\mu,\lambda}
=
\mu\widetilde{W}_n + \lambda (-Z_n),
$$
with non-negative weights~$\mu$ and~$\lambda$. While any~$Q_n^{\mu,\lambda}$ provides a reasonable test statistic for the problem at hand, only one set of weights will yield a Le Cam optimal test and, interestingly, this set of weights depends on the way~$\kappa_n$ behaves with~$p_n$ and~$n$. This will be one of the many consequences of the following LAN result.

\begin{theorem}
\label{TheoLAN}
Let~$(p_n)$ be a sequence of integers that diverges to infinity,~$(\kappa_n)$ be a sequence in~$(0,\infty)$, and~$(\thetab_{n0})$ be a sequence such that $\thetab_{n0}$ belongs to~$\mathcal{S}^{p_n-1}$ for any~$n$. 
 Then, there exist a sequence~$(\nu_n)$ in~$(0,\infty)$ and a sequence of random variables~$(\Delta_n)$ that is asymptotically normal with zero mean and variance~$\Gamma$ under~${\rm P}^{(n)}_{\thetab_{n0},\kappa_n}$ such that, for any bounded sequence~$(\taub_n)$ such that~$\thetab_n=\thetab_{n0}+\nu_n\taub_n$ belongs to~$\mathcal{S}^{p_n-1}$ for any~$n$,
$$
\hspace{-5mm} 
\Lambda^{(n){\rm inv}}_{\thetab_n/\thetab_{n0};\kappa_n}
%\log 
%\Bigg( 
%\frac{d{\rm P}^{(n)\Tb_n}_{\thetab_n,\kappa_n}}{d{\rm P}^{(n)\Tb_n}_{\thetab_{n0},\kappa_n}}
%\Bigg)
=
\|\taub_n\|^2  \Delta_n 
-
\frac{1}{2} \|\taub_n\|^4 \Gamma
+o_{\rm P}(1)
$$
as~$n\to\infty$ under~${\rm P}^{(n)}_{\thetab_{n0},\kappa_n}\!$. 
If (i)~$\kappa_n/p_n\to \infty$, then
$$
\nu_n
=
\frac{p_n^{1/4}}{\sqrt{n\kappa_n}}
%\frac{p_n^{3/4}}{\sqrt{n} \kappa_n \sqrt{f_{n2}}}
,
\quad
\Delta_n
=
\frac{\widetilde{W}_n}{\sqrt{2}}
,
\quad
\textrm{ and } 
\quad
\Gamma=\frac{1}{2}
;
$$
if 
(ii)~$\kappa_n/p_n\to \xi>0$, then, letting~$c_{\xi}:=\frac{1}{2}+\sqrt{\frac{1}{4} + \xi^2}$,
$$
\nu_n
=
\frac{\sqrt{c_\xi}\, p_n^{3/4}}{\sqrt{n} \kappa_n}
%\frac{p_n^{3/4}}{\sqrt{n} \kappa_n \sqrt{f_{n2}}}
,
\quad
\Delta_n
=
\frac{\widetilde{W}_n}{\sqrt{2}}
,
\quad
\textrm{ and } 
\quad
\Gamma=\frac{1}{2}
;
$$
if (iii)~$\kappa_n/p_n\to 0$ with~$\sqrt{n} \kappa_n/p_n\to \infty$, then
$$
\nu_n
=
\frac{p_n^{3/4}}{\sqrt{n} \kappa_n}
%\frac{p_n^{3/4}}{\sqrt{n} \kappa_n \sqrt{f_{n2}}}
,
\quad
\Delta_n
=
\frac{\widetilde{W}_n}{\sqrt{2}}
,
\quad
\textrm{ and } 
\quad
\Gamma=\frac{1}{2}
;
$$
if (iv)~$\sqrt{n} \kappa_n/p_n\to \xi >0$, then
$$
\nu_n
=
\frac{p_n^{3/4}}{\sqrt{n} \kappa_n}
%\nu_n
%=
%\frac{p_n^{3/4}}{\sqrt{n} \kappa_n \sqrt{f_{n2}}}
,
\quad
\Delta_n
=
\frac{\widetilde{W}_n}{\sqrt{2}}
-
\frac{Z_n}{2\xi}
,
\quad
\textrm{ and } 
\quad
\Gamma
=
\frac{1}{2}
+
\frac{1}{4\xi^2}
;
$$
if (v)~$\sqrt{n} \kappa_n/p_n\to 0$ with~$\sqrt{n} \kappa_n/\sqrt{p_n}\to \infty$, then
$$
\nu_n
=
\frac{p_n^{1/4}}{n^{1/4} \sqrt{\kappa_n}}
,
\quad
\Delta_n
=
-\frac{Z_n}{2}
,
\quad
\textrm{ and } 
\quad
\Gamma
=
\frac{1}{4}
;
$$
if (vi)~$\sqrt{n} \kappa_n/\sqrt{p_n}\to \xi>0$, then
$$
\nu_n
=
1
,
\quad
\Delta_n
=
-\frac{\xi Z_n}{2}
,
\quad
\textrm{ and } 
\quad
\Gamma
=
\frac{\xi^2}{4}
;
$$
finally, if (vii)~$\sqrt{n} \kappa_n/\sqrt{p_n}\to 0$, then, even with~$\nu_n=1$, 
the invariant log-likelihood ratio~$\Lambda^{(n){\rm inv}}_{\thetab_n/\thetab_{n0};\kappa_n}
$ is~$o_{\rm P}(1)$ as~$n\to\infty$ under~${\rm P}^{(n)}_{\thetab_{n0},\kappa_n}\!$. 
\end{theorem}
\vspace{3mm}

%% $$
%% \frac{\widetilde{W}_n}{\sqrt{2}}
%%-
%%\frac{Z_n}{2\xi}
%%\approx
%% \frac{W^*_n}{\sqrt{2}}
%%-
%%\frac{Z_n}{2\hat{\xi}}
%%=
%%\frac{n(p_n-1)f_{n2}^{-1}\bar{\Xb}_n\pr (\mathbf{I}_{p}-\thetab_0\thetab_0\pr) \bar{\Xb}_n-(p_n-1)}{2\sqrt{p_n-1}}
%%-
%%\frac{p_n Z_n}{2\sqrt{n}\kappa_n}
%% $$
%% $$
%% =
%%\frac{n \kappa_n e_{n1}^{-1}\bar{\Xb}_n\pr (\mathbf{I}_{p}-\thetab_0\thetab_0\pr) \bar{\Xb}_n-(p_n-1)}{2\sqrt{p_n-1}}
%%-
%%\frac{p_n Z_n}{2\sqrt{n}\kappa_n}
%% $$
%% 
% 
% 
% 
% \newpage
 
In the image model~(\ref{InvariantFvMLstatmodels}), 
\vspace{-.6mm}
 the spherical location problem consists in testing~$\mathcal{H}_0\n:u=1$ against~$\mathcal{H}_0\n:u<1$. In the localized at~$u=u_0=1$ experiments, parametrized by~$u=u_0-\frac{1}{2}\nu_n^2\|\taub_n\|^2$ 
 \vspace{-.6mm}
as in~(\ref{ilr}), this reduces to testing~$\mathcal{H}_0\n:\|\taub_n\|=0$ against~$\mathcal{H}_1\n:\|\taub_n\|>0$. In any given regime~(i)--(vii) from Theorem~\ref{TheoLAN}, it directly follows from this theorem that a locally asymptotically most powerful test for this problem --- hence, \emph{locally asymptotically most powerful invariant} test for the original spherical location problem --- rejects the null at asymptotic level~$\alpha$ whenever
\begin{equation}
	\label{opttest}
\Delta_n/\sqrt{\Gamma} > \Phi^{-1}(1-\alpha)
,
\end{equation}
where~$\Phi$ denotes the cumulative distribution function of the standard normal distribution (in the rest of the paper, the term ``\emph{optimal}" will refer to this particular Le Cam optimality concept). A routine application of the Le Cam third lemma then shows that, in each regime,
\vspace{-.3mm}
 the asymptotic distribution of~$\Delta_n$, under the corresponding contiguous
 \vspace{-.5mm}
 alternatives~${\rm P}^{(n)}_{\thetab_{n0}+\nu_n\taub_n,\kappa_n}$ with~$\|\taub_n\|\to t$, is normal with mean~$\Gamma t^2$ and variance~$\Gamma$, so that the resulting asymptotic power of the optimal test in~(\ref{opttest}) is
\begin{eqnarray}
\lefteqn{
\lim_{n\to\infty}
{\rm P}^{(n)}_{\thetab_{n0}+\nu_n\taub_n,\kappa_n}
\big[
\Delta_n/\sqrt{\Gamma} > \Phi^{-1}(1-\alpha)
\big]
}
\nonumber
\\[2mm]
& & 
\hspace{28mm} 
=
1-\Phi\Big(
\Phi^{-1}(1-\alpha)
-
\sqrt{\Gamma} t^2
\Big)
.
\label{powgeneraloptimal}	
\end{eqnarray}

In each regime~(i)--(vii),~$\nu_n$ is the contiguity rate, which implies that the least severe alternatives 
\vspace{-.4mm}
 under which a test may have non-trivial asymptotic powers are of the form~${\rm P}^{(n)}_{\thetab_{n0}+\nu_n\taub_n,\kappa_n}$, with a sequence~$(\|\taub_n\|)$ that is~$O(1)$ but not~$o(1)$. Theorem~\ref{TheoLAN} shows that this contiguity rate depends on the regime considered and does so in a monotonic fashion, which is intuitively reasonable: the larger~$\kappa_n$ (that is, the easier the inference problem), the faster~$\nu_n$ goes to zero, that is, the less severe the alternatives that can be detected by rate-consistent tests. Because the unit sphere~$\mathcal{S}^{p_n-1}$ has a fixed diameter, $\nu_n= 1$ characterizes the most severe alternatives that can be considered. In regime~(vi), no tests will therefore be consistent under such most severe alternatives, while, in regime~(vii), the distribution is so close to the uniform distribution on~$\mathcal{S}^{p_n-1}$ that no tests can show non-trivial asymptotic powers under such alternatives, so that even the trivial $\alpha$-test is optimal.  

One of the most striking consequences of Theorem~\ref{TheoLAN} is that the optimal test depends on the regime considered. In regimes~(v)--(vii), the optimal test in~(\ref{opttest}) rejects the null when~$Z_n<\Phi^{-1}(\alpha)$; of course, this optimality is degenerate in regime~(vii), where any invariant test with asymptotic level~$\alpha$ would also be optimal. In contrast, the optimal $\alpha$-level test in regimes~(i)--(iii) rejects the null when
$$ 
\widetilde{W}_n
=
\frac{W_n-(p_n-1)}{\sqrt{2(p_n-1)}}
>
\Phi^{-1}(1-\alpha)
.
$$
Since the chi-square distribution with~$p-1$ degrees of freedom converges, after standardization via its mean~$p-1$ and standard deviation~$\sqrt{2(p-1)}$, to the standard normal distribution as~$p$ diverges to infinity, this test is asymptotically equivalent to the Watson test in~(\ref{WatsonLD}), based obviously on the dimension~$p=p_n$ at hand. This shows that, in regimes~(i)--(iii), the traditional, low-dimensional, Watson test is optimal in high dimensions. In regime~(iv), which is at the frontier between these regimes where the optimal test is the Watson test and those where the optimal test is based on~$Z_n$, the optimal test is quite naturally based on a linear combination of~$\widetilde{W}_n$ and~$Z_n$.  

Finally, the Le Cam third lemma allows us to derive the asymptotic non-null behavior of the Watson test under the contiguous alternatives considered in any regime~(i)--(vii). In regimes~(i)--(iv),
\vspace{-.7mm}
  the limiting powers under contiguous alternatives of the form~${\rm P}^{(n)}_{\thetab_{n0}+\nu_n\taub_n,\kappa_n}$, with~$\|\taub_n\|\to t$, are given by 
\begin{equation}
\label{asymppoweroptimalitoiii}
	1-\Phi\bigg(
\Phi^{-1}(1-\alpha)
-
\frac{t^2}{\sqrt{2}}
\bigg)
.
\end{equation}
In regimes~(i)--(iii), the Watson test is the optimal test and these asymptotic powers are equal to those in~(\ref{powgeneraloptimal}), whereas in regime~(iv), the Watson test is only rate-consistent, as the corresponding asymptotic powers of the optimal test are
\begin{equation}
\label{asymppoweroptimaliv}
1-\Phi\bigg(
\Phi^{-1}(1-\alpha)
-
t^2
\sqrt{\frac{1}{2}
+
\frac{1}{4\xi^2}}
\
\bigg)
.
\end{equation} 
In regimes~(v)--(vi), the Le Cam third lemma shows that the limiting powers of the Watson test, still under the corresponding contiguous alternatives, are equal to the nominal level~$\alpha$, so that the Watson test is not even rate-consistent in those regimes. Finally, as already discussed, the Watson test is optimal in regime~(vii), but trivially so since the trivial $\alpha$-test there also is. 

%%%%%%%

\subsection{Optimal testing under unspecified~$\kappa_n$}
\label{subseckappaunspec}

The optimal test in regimes (i)--(iii), namely the Watson test, is a genuine test in the sense that it can be applied on the basis of the observations only. In contrast, the optimal tests in regimes~(iv)--(vi) are ``oracle" tests since they require knowing the values of~$
e_{n1}$ and~$\tilde{e}_{n2}$, or equivalently (see~(\ref{linken1en2tildekappan})), the value of the concentration~$\kappa_n$. This concentration, however, can hardly be assumed to be specified in practice, so that it is natural to wonder what is the optimal test, in regimes~(iv)--(vi), when~$\kappa_n$ is treated as a nuisance parameter.

We first focus on regime~(iv). There, the concentration~$\kappa_n$ is asymptotically of the form~$\kappa_n=p_n \xi/\sqrt{n}$ for some $\xi >0$. Within regime~(iv), $\xi$, obviously, is a perfectly valid alternative concentration parameter. Inspired by the classical treatment of asymptotically optimal inference in the presence of nuisance parameters (see, e.g., \cite{Bic1998}), this suggests studying the asymptotic behavior of invariant log-likelihood ratios of the form
\begin{equation*}
\label{ilrbis}	
\Lambda^{(n){\rm inv}}_{\thetab_n, \kappa_{n,s}/\thetab_{n0},\kappa_{n}}
:=
\log \frac{d{\rm P}^{(n)\Tb_n}_{1-\nu_n^2 \|\taub_n\|^2/2,\kappa_{n,s}}}{d{\rm P}^{(n)\Tb_n}_{1,\kappa_{n}}},
\end{equation*}
where~$\kappa_{n,s}:=p_n(\xi+ \vartheta_n s)/\sqrt{n}$
%=\frac{p_n (\xi+\frac{1}{\sqrt{p_n}}s)}{\sqrt{n}}
%=\frac{p_n \xi}{\sqrt{n}}+ \frac{\sqrt{p_n} s}{\sqrt{n}}$
is a suitable sequence of perturbed concentrations. We have the following result.

\begin{theorem}\label{LANunspec}
Let~$(p_n)$ be a sequence of integers that diverges to infinity with $p_n=o(n^2)$ as $\ny$. 
Let~$\kappa_{n}:=p_n \xi/\sqrt{n}$, with~$\xi>0$, and~$\kappa_{n,s}:=p_n(\xi+ s/\sqrt{p_n})/\sqrt{n}$, where~$s$ is such that $\xi+ s/\sqrt{p_n} >0$ for any~$n$.
Let the sequence~$(\thetab_{n0})$ in~$\mathcal{S}^{p_n-1}$ and the bounded sequence~$(\taub_n)$ in~$\R^{p_n}$ be such that $\thetab_{n0}$ and~$\thetab_n=\thetab_{n0}+\nu_n\taub_n$, with the~$\nu_n$ below, belong to~$\mathcal{S}^{p_n-1}$ for any~$n$. Then, putting ${\bf t}_n:=(\| \taub_n\|^2, s)\pr$, 
$$
\nu_n
:=
\frac{p_n^{3/4}}{\sqrt{n} \kappa_{n}}
,
\ \,
\Deltab_n
:=
\left(
\begin{array}{c} 
\frac{\widetilde{W}_n}{\sqrt{2}} -\frac{Z_n}{2\xi} \\[1mm]
Z_{n}
\end{array}
\right)
,
\
\textrm{ and } 
\ 
\Gamb
:= 
\left(\begin{array}{cc} 
\frac{1}{2}+\frac{1}{4\xi^2} & -\frac{1}{2\xi} \\[2mm] 
-\frac{1}{2\xi} & 1 \end{array}\right)
,
$$
we have
\begin{equation}
\label{LAQdoubleiv}	
\Lambda^{(n){\rm inv}}_{\thetab_n, \kappa_{n,s}/\thetab_{n0},\kappa_{n}}={\bf t}_n\pr \Deltab_n- \frac{1}{2} {\bf t}_n\pr{\Gamb}{\bf t}_n+o_{\rm P}(1)
\end{equation}
as $\ny$ under~${\rm P}^{(n)}_{\thetab_{n0},\kappa_n}$, where~$\Deltab_n$, under the same sequence of hypotheses, is asymptotically normal with mean zero and covariance matrix~$\Gamb$. 
\end{theorem}

Theorem \ref{LANunspec} shows that, in regime (iv), the sequence of high-dimensional FvML experiments is jointly LAN in the location and concentration parameters. The corresponding Fisher information matrix~$\Gamb=(\Gamma_{ij})$ is not diagonal, which entails that the unspecification of the concentration parameter has asymptotically a positive cost when performing inference on the location parameter. In the present joint LAN framework, Le Cam optimal inference for location under unspecified concentration is to be based (see again \cite{Bic1998}) on the residual of the regression (in the limiting Gaussian shift experiment) of the location part~$\Delta_{n1}$ of the central sequence~$\Deltab=(\Delta_{n1},\Delta_{n2})'$ with respect to the concentration part~$\Delta_{n2}$, that is, is to be based on the \emph{efficient central sequence}
\begin{equation}
	\label{effcen}
\Delta_{n1}^*
:=
\Delta_{n1}- \frac{\Gamma_{12}}{\Gamma_{22}}\, \Delta_{n2}
.
\end{equation}
Under the null, $\Delta_{n1}^*$ is asymptotically normal with mean zero and variance~$\Gamma_{11}^*
\linebreak
=\Gamma_{11}-\Gamma_{12}^2/\Gamma_{22}$, and the Le Cam optimal location test under unspecified~$\kappa_n$ rejects the null at asymptotical level~$\alpha$ when
$$
\Delta_{n1}^* / \sqrt{\Gamma_{11}^*} 
=
\widetilde{W}_n
> \Phi^{-1}(1-\alpha)
.
$$
As a corollary, provided that~$p_n=o(n^2)$, the unspecified-$\kappa_n$ optimal test in regime~(iv) is the Watson test. Consequently, the difference between the local asymptotic powers in~(\ref{asymppoweroptimalitoiii}) and~(\ref{asymppoweroptimaliv}), associated with the Watson test and the specified-$\kappa_n$ optimal test in regime~(iv), respectively, can be interpreted as the asymptotic cost of the unspecification of the concentration when performing inference on location in the regime considered. Note that the optimal specified-$\kappa_n$ test and optimal unspecified-$\kappa_n$ test exhibit the same consistency rates, so that the cost of not knowing~$\kappa_n$ lies in the difference of powers these tests show under contiguous alternatives.

We now turn to regime~(vi), where the concentration~$\kappa_n$ is asymptotically of the form~$\kappa_n=\sqrt{p_n}\xi/\sqrt{n}$. In this regime, taking~$\nu_n=1$ (as in Theorem~\ref{TheoLAN}) and perturbed concentrations of the form~$\kappa_{n,s}:=\sqrt{p_n}(\xi+s)/\sqrt{n}$, it is easy to show, by working along the same lines as in the proof of Theorem \ref{LANunspec}, that the sequence of experiments is also jointly LAN in location and concentration, this time without any condition on~$p_n$. The corresponding central sequence and 
Fisher information matrix are
\begin{equation} 
\label{infofive}
\Deltab_n
:=
\bigg(
\!
\begin{array}{c} 
\Delta_{n1} \\[1mm]
\Delta_{n2}
\end{array} 
\!
\bigg)
=
\bigg(
\!
\begin{array}{c} 
%-\frac{Z_n}{2} \\[1mm]
-Z_n/2 \\[1mm]
Z_n
\end{array} 
\!
\bigg)
\quad
\textrm{and}
\quad
\Gamb
:= 
%\bigg(
%\begin{array}{cc} 
%\frac{1}{4} & -\frac{1}{2} \\[1mm]
%-\frac{1}{2} & \ 1 
%\end{array} 
%\bigg)
\bigg(
\begin{array}{cc} 
1/4 & -1/2 \\[1mm]
-1/2 & \, 1 
\end{array} 
\bigg)
.
\end{equation}
The collinearity between the location part~$\Delta_{n1}$ and concentration part~$\Delta_{n2}$ of the central sequence implies that the efficient central sequence~$\Delta_{n1}^*$  is zero in regime~(vi). As a result, for the unspecified concentration problem, no test can detect alternatives in~$\nu_n=1$ in regime~(vi), which is in line with the corresponding trivial asymptotic powers of the Watson test in Section~\ref{subseckappaspec}. Since~$\nu_n = 1$ provides the most severe location alternatives than can be considered, we conclude that, for the unspecified concentration problem, no test in regime~(vi) can do asymptotically better than the trivial $\alpha$-level test that randomly rejects the null with probability~$\alpha$. Under unspecified~$\kappa_n$, thus, the Watson test is optimal in regime~(vi), too, even if it is in a degenerate way.

Finally, we consider regime~(v), where the situation is more complicated.
% and where the cost for not knowing~$\kappa_n$ is expected to be intermediate between those observed in regimes~(iv) and~(vi). Regime~(v)  
This regime is associated with~$\kappa_n=p_n r_n\xi/\sqrt{n}$, where~$\xi>0$ and~$(r_n)$ is a 
\vspace{-.7mm}
positive sequence satisfying~$r_n=o(1)$ and~$r_n \sqrt{p_n} \to \infty$. If one takes $\nu_n=p_n^{1/4}/(n^{1/4} \sqrt{\kappa_n})$ (still as in Theorem~\ref{TheoLAN}) and considers perturbed concentrations of the form~$\kappa_{n,s}=p_n r_n(\xi+s/(\sqrt{p_n}r_n))/\sqrt{n}$, then it can be shown 
%(by working along the same lines as in the proof of Theorem \ref{LANunspec}) 
that, provided that~$p_n=o(n^2r_n^{-6})$, the resulting sequence of experiments is still jointly LAN in location and concentration, with the same central sequence and Fisher information matrix as in~(\ref{infofive}). Consequently, the corresponding efficient central sequence~$\Delta_{n1}^*$ is zero again, so that no unspecified-$\kappa_n$ test can detect deviations from the null hypothesis at the $\nu_n$-rate in regime~(v). %, which is in line with the corresponding trivial asymptotic powers of the Watson in Section~\ref{subseckappaspec}. 
Unlike in regime~(vi), however, alternatives that are more severe than the contiguous ones can be considered in regime~(v). As a consequence, several important questions  are left wide open in regime~(v) for the unspecified-$\kappa_n$ problem: (1)~are there alternatives that can be detected by an unspecified-$\kappa_n$ test? (2) If so, what are the least severe ones that can be detected by such a test and (3) what is the Le Cam optimal test (if any)? (4) Are there alternatives that can be detected by the Watson test? (5) Does this test enjoy any Le Cam optimality property in this regime?  

%As we showed in the previous paragraph, a joint LAN result in the original location and concentration parameter does not provide answers to these questions. 
To answer these questions, one needs to orthogonalize the parameter of interest~$u=\thetab_n'\thetab_{n0}$ and concentration parameter~$\kappa$. In regime~(iv), this orthogonalization was achieved, within the LAN framework of Theorem~\ref{LANunspec}, by the efficient central sequence in~(\ref{effcen}). 
\vspace{-.5mm}
 In regime~(v), where the consistency rates of the $Z_n$ and~$\widetilde{W}_n$ tests do not match, this approach does not work and it is needed to perform orthogonalization by introducing explicitly a new parametrization ({such an orthogonalization through reparametrization is suitable when Fisher information matrices are singular; see, e.g., \cite{HL14}}). The following LAN result relates to this new parametrization of the statistical experiments at hand, that involves the same parameter of interest~$u=\thetab_n'\thetab_{n0}$ and the alternative concentration parameter~$\bar{\kappa}_n=\kappa_n/u$ (of course, this reparametrization requires restricting to the hemisphere associated with~$u>0$, which still allows us to consider ``local" alternatives).

\begin{theorem}\label{LANunspecregime5}
Let~$(p_n)$ be a sequence of integers that diverges to infinity with $p_n=o(n^2 r_n^{-4})$ as $\ny$, where~$(r_n)$ is a positive real sequence such that~$r_n=o(1)$ and~$\sqrt{p_n}r_n\to\infty$. 
Let~$(\thetab_{n0})$ be a sequence in~$\mathcal{S}^{p_n-1}$ and~$(\taub_n)$ be a bounded sequence in~$\R^{p_n}$ such that~$\thetab_n=\thetab_{n0}+\nu_n\taub_n$, with the~$\nu_n$ below, belongs to~$\mathcal{S}^{p_n-1}$ for any~$n$. 
Let~$\kappa_{n}:=p_n r_n \xi/\sqrt{n}$, with~$\xi>0$ and
$$
\kappa_{n,s,\taub_n}
=
\frac{p_n r_n(\xi+s/(\sqrt{p_n}r_n))}{\sqrt{n}(1-{\textstyle{\frac{1}{2}}}\nu_n^2\|\taub_n\|^2)}
=:
\frac{\rho_n p_n r_n}{\sqrt{n}}(\xi+s/(\sqrt{p_n}r_n)) 
,
$$
%$$
%\textcolor{red}{ 
%\kappa_{n,s,\taub_n}-\kappa_{n}
%=
%\frac{(\rho_n-1) p_n r_n}{\sqrt{n}} \xi
%+
%\frac{\rho_n p_n r_n}{\sqrt{n}} s/(\sqrt{p_n}r_n)
%=
%\frac{\rho_n p_n r_n\nu_n^2\|\taub_n\|^2}{2\sqrt{n}} \xi
%+
%\frac{\rho_n \sqrt{p_n}}{\sqrt{n}} s
%}
%$$
%$$
%\textcolor{red}{
%=
%\frac{\rho_n \sqrt{p_n}}{\sqrt{n}}
%\bigg[
%\frac{\sqrt{p_n} r_n\nu_n^2\|\taub_n\|^2}{2} \xi
%+
%s
%\bigg]
%}
%$$
where~$s$ is such that~$\xi+s/(\sqrt{p_n}r_n)>0$ for any~$n$. Assume that, still with the~$\nu_n$ below, ${\textstyle{\frac{1}{2}}}\nu_n^2\|\taub_n\|^2$ is upper-bounded by~$1-\delta$ for some~$\delta>0$. 
Then, putting  
$$
(a) 
\qquad
\nu_n
=
\frac{p_n^{3/4}}{\sqrt{n} \kappa_{n}}
,
\ \,
\ \
\Cb_n
= 
\bigg(\begin{array}{cc} 
1\ & 0 \\[1mm] 
0\ & 1 \end{array}\bigg)
,
\hspace{3mm} 
\Deltab_n
=
\bigg(
\begin{array}{c} 
\frac{\widetilde{W}_n}{\sqrt{2}} \\[2mm]
Z_{n}
\end{array}
\bigg)
,
\
\textrm{ and } 
\ 
\Gamb
= 
\bigg(\begin{array}{cc} 
\frac{1}{2}\ & 0 \\[1.5mm] 
0\ & 1 \end{array}\bigg)
,
$$
$$
\hspace{1.5mm} 
(b) 
\quad
\nu_n
=
1
,
\ 
\, 
\Cb_n
= 
\bigg(\begin{array}{cc} 
\xi^2\big(1-\frac{\|\taub_n\|^2}{4}\big)\ & 0 \\[2mm] 
0\ & 1 \end{array}\bigg)
,
\hspace{1mm} 
\Deltab_n
=
\bigg(
\begin{array}{c} 
\frac{\widetilde{W}_n}{\sqrt{2}} \\[2mm]
Z_{n}
\end{array}
\bigg)
,
\,
\textrm{ and } 
\
\Gamb
= 
\bigg(\begin{array}{cc} 
\frac{1}{2}\ & 0 \\[1.5mm] 
0\ & 1 \end{array}\bigg)
,
$$
or
$$
\hspace{-10mm} 
(c) 
\qquad
\nu_n
=
1
,
\ \,
\ \
\Cb_n
= 
\bigg(\begin{array}{cc} 
1\ & 0 \\[1mm] 
0\ & 1 \end{array}\bigg)
,
\hspace{3mm} 
\Deltab_n
=
\bigg(
\begin{array}{c} 
0 \\[1mm]
Z_{n}
\end{array}
\bigg)
,
\
\textrm{ and } 
\ 
\Gamb
= 
\bigg(\begin{array}{cc} 
0\ & 0 \\[1mm] 
0\ & 1 \end{array}\bigg)
,
$$
depending on whether~(a) $\rho_n p_n^{1/4}r_n\to\infty$, (b) ${\rho_n p_n^{1/4}r_n\to 1}$, or (c) $\rho_n p_n^{1/4}r_n=o(1)$, respectively, we have, with ${\bf t}_n:=(\| \taub_n\|^2, s)\pr$,
\begin{equation}
\label{LAQdoubleiv}	
\Lambda^{(n){\rm inv}}_{\thetab_n, \kappa_{n,s,\taub_n}/\thetab_{n0},\kappa_{n}}={\bf t}_n\pr \Cb_n\Deltab_n- \frac{1}{2} {\bf t}_n\pr\Cb_n^2{\Gamb}{\bf t}_n+o_{\rm P}(1)
\end{equation}
as $\ny$ under~${\rm P}^{(n)}_{\thetab_{n0},\kappa_n}$, where~$\Deltab_n$, under the same sequence of hypotheses, is asymptotically normal with mean zero and covariance matrix~$\Gamb$. 
\end{theorem}
\vspace{3mm}

The block-diagonality of the three Fisher information matrices~$\pmb\Gamma$ in this result confirms that the new parametrization achieves orthogonalization in regime~(v). More importantly, Theorem~\ref{LANunspecregime5} allows us to answer the open questions above. In this purpose, the key observation is that   the problem of testing the null hypothesis~$\mathcal{H}_0:u=1$ against the alternative~$\mathcal{H}_1:u<1$ under unspecified~$\kappa_n$ in the original parametrization is strictly equivalent to the problem of testing the null hypothesis~$\mathcal{H}_0:u=1$ against the alternative~$\mathcal{H}_1:u<1$ under unspecified~$\bar{\kappa}_n$ in the new parametrization. Therefore, the $s\equiv 0$ version of Theorem~\ref{LANunspecregime5} establishes the following: in regime~(v$_a$), which refers to case~(a) in this result, the Watson test is Le Cam optimal for the \mbox{unspecified-$\kappa_n$} problem and 
\vspace{-.4mm}
 will show non-trivial asymptotic powers under alternatives associated with~$p_n^{3/4}/(\sqrt{n} \kappa_{n})$ (the Le Cam third lemma readily implies that these asymptotic powers are equal to those in~(\ref{asymppoweroptimalitoiii})). In regime~(v$_c$), no unspecified-$\kappa_n$ test can detect even the most severe alternatives associated with~$\nu_n=1$. In the boundary case of regime~(v$_b$), the situation is more complex, as the sequence of statistical experiments there is not LAN. Yet, the result shows that the least severe alternatives that can be detected by an unspecified-$\kappa_n$ test are those associated with~$\nu_n=1$ and that the Watson test is rate-consistent. Theorem~\ref{LANunspecregime5}(b) also shows that the Watson test is Le Cam optimal for small departures~$\taub_n$ of the null hypothesis (this follows from the fact that the usual LAN property is obtained for small~$\|\taub_n\|$); we refer to Theorem~4.1(iii) in~\cite{Remy2019} for a similar phenomenon in low dimensions. This thoroughly answers the questions~(1)--(5) raised above. 

Wrapping up, we proved that the Watson test is optimal in regimes~(i)--(iii) only for the specified concentration problem and that it is optimal in \emph{all} regimes in the more important unspecified concentration one (in regimes~(iv)--(v$_a$), optimality requires a constraint on~$p_n$ that is at most~$p_n=o(n^2)$, and optimality is only local in~$\taub_n$ in regime~(v$_b$)). The asymptotic cost due to the unspecification of the concentration is nil in regimes~(i)--(iii) (and~(vii)), affects limiting powers but not consistency rates in regime~(iv), and is in terms of consistency rates in regimes~(v)--(vi). {Table~\ref{thethetable} provides a summary of the optimality results we obtained both for the specified-$\kappa_n$ and unspecified-$\kappa_n$ problems.

\begin{center}
\begin{table}[h]
 \begin{center}
\begin{tabular}{|c|c|c|c|}
\hline
&&&   \\
$\#$ & Regime & {\emph{$\kappa_n$ specified}} &  \emph{$\kappa_n$ unspecified} \\
&&&  \\
 \hline
 && &  \\
(i)& $\kappa_n/p_n\to \infty$ &$ {\widetilde{W}_n}$  & $\widetilde{W}_n$\\
&&&\\
(ii)& $\kappa_n/p_n\to \xi>0$ & $ {\widetilde{W}_n}$& $ {\widetilde{W}_n}$ \\
&&& \\
(iii)& $\kappa_n/p_n\to 0$ with~$\sqrt{n} \kappa_n/p_n\to \infty$ & $ {\widetilde{W}_n}$ & $ {\widetilde{W}_n}$  \\
&&& \\
(iv)& $\sqrt{n} \kappa_n/p_n\to \xi >0$ &  $\frac{\widetilde{W}_n}{\sqrt{2}} 
-
\frac{Z_n}{2\xi}$ & $ {\widetilde{W}_n}$ ($\star$) \\
&&& \\
(v$_a$)& $\sqrt{n} \kappa_n/p_n\to 0$ with~$\sqrt{n} \kappa_n/p_n^{3/4}\to \infty$ & $Z_n$ & $ {\widetilde{W}_n}$ ($\star$) \\
&&& \\
(v$_b$)& $\sqrt{n} \kappa_n/p_n^{3/4}\to \xi >0$ & $Z_n$ & $ {\widetilde{W}_n}$ ($\dagger$) \\
&&& \\
(v$_c$)& $\sqrt{n} \kappa_n/p_n^{3/4}\to 0$ with~$\sqrt{n} \kappa_n/\sqrt{p_n}\to \infty$ & $Z_n$ & $\emptyset $ \\
&& & \\
(vi) & $\sqrt{n}\kappa_n/\sqrt{p_n} \to \xi >0$  & $Z_n$ & $\emptyset $ \\
 && & \\
(vii) & $\sqrt{n}\kappa_n/\sqrt{p_n} \to 0$ & $\emptyset$ & $\emptyset $\\
 &&& \\
\hline
\end{tabular}
  \end{center}
 \caption{The test statistics on which locally asymptotically optimal tests are based in the various asymptotic regimes for both the specified-$\kappa_n$ and unspecified-$\kappa_n$ problems. The symbol~$\emptyset$ means that no test can detect even the most severe alternatives associated with~\mbox{$\nu_n=1$}. The symbol~$\star$ indicates that the result is obtained provided that~$p_n=o(n^2)$ (for the $\star$ in regime~(v$_a$), the constraint is actually milder than~$p_n=o(n^2)$; see Theorem~\ref{LANunspecregime5} for details). The symbol~$\dagger$ stresses that, in the non-standard limiting experiment obtained in regime~(v$_b$) for unspecified~$\kappa_n$, Le Cam optimality is achieved only locally in~$\taub_n$.} 
 \label{thethetable} 
\end{table}
\end{center}

%%%%%%%%%%%%%%%%%%%%%%%%%%%%%%%%%%%%%%%%%%%%%%
 
\section{Non-null investigation via martingale CLTs}
\label{secBilling}

The results above thoroughly describe the asymptotic non-null and optimality properties of the Watson test in the FvML case and provide a strong motivation to use this test in this specific parametric framework. While the Watson test remains valid (in the sense that it still meets the asymptotic nominal level constraint) under much broader distributional assumptions, it is unclear how well this test behaves under high-dimensional non-FvML alternatives (we refer to \cite{PaiVer2015a}, \cite{PaiVer17b} and~\cite{PaiVer19} for an extensive study of the low-dimensional case). In this section, we therefore investigate, through a different approach relying on martingale CLTs, the non-null behavior of the Watson test under general rotationally symmetric distributions.  

Recall that the distribution of a random vector~$\Xb$ with values in~$\mathcal{S}^{p-1}$ is \emph{rotationally symmetric about~$\thetab(\in\mathcal{S}^{p-1})$} if~$\Ob\Xb$ and~$\Xb$ share the same distribution for any~$\Ob\in SO_{\thetab}(p)$, and that it is \emph{rotationally symmetric} if it is rotationally symmetric about some~$\thetab$ in~$\mathcal{S}^{p-1}$. Clearly, if~$\Xb$ has an ${\rm FvML}_p(\thetab,\kappa)$ distribution, then it is rotationally symmetric about~$\thetab$, so that the distributional context considered in this section will encompass the one in Section~\ref{secLAN}. Parallel to what was done there, we will % Irrespective of the fact that~$\Xb$ is rotationally symmetric about~$\thetab$ or not, we will often
refer to the decomposition 
$ 
\Xb
=
U \thetab
+
V \Sb
,
$
with~$U=\Xb'\thetab$, $V=\sqrt{1-U^2}$ and~$\Sb=(\mathbf{I}_p-\thetab\thetab\pr)\Xb/\|(\mathbf{I}_p-\thetab\thetab\pr)\Xb\|$, as the tangent-normal decomposition of~$\Xb$ with respect to~$\thetab$. 
If~$\Xb$ is rotationally symmetric about~$\thetab$, then $\Sb$ is uniformly distributed over~$\{\xb\in\mathcal{S}^{p-1}:\xb'\thetab=0\}$ and is independent of~$U$. The distribution of~$\Xb$ is then fully determined by~$\thetab$ and by the cumulative distribution function~$F$ of~$U$, which justifies denoting the corresponding distribution as~${\rm Rot}_p(\thetab,F)$. In the sequel, we tacitly restrict to classes of rotationally symmetric distributions making~$\thetab$ identifiable, which typically only excludes distributions satisfying~${\rm Rot}_p(-\thetab,F)={\rm Rot}_p(\thetab,F)$. 
%Tacitly, we will assume that~$F$ is such that We will throughout tacitly assume that~$\{u\in [-1,1]: F(u)<1-\lim_{t\stackrel{<}{\to} -u} F(t)\}$ is non-empty and has a negative lower bound. It is easy to check that this is the minimal assumption ensuring identifiability of~$\thetab$ under rotational symmetry, in the sense that if~${\rm Rot}_p(\thetab,F)={\rm Rot}_p(\etab,G)$, then we must have~$\thetab=\etab$ (hence also~$F=G$). 

We consider then a triangular array of observations of the form~$\Xb_{ni}$, $i=1,\ldots,n$, $n=1,2,\ldots$, where~$\Xb_{n1},\ldots,\Xb_{nn}$ form a random sample from the rotationally symmetric
\vspace{-.6mm}
  distribution~${\rm Rot}_{p_n}(\thetab_n,F_n)$. The corresponding hypothesis, that will be
  \vspace{-.5mm}
  denoted as~${\rm P}\n_{\thetab_n,F_n}$ involves a sequence of integers~$(p_n)$ diverging to infinity, a sequence~$(\thetab_n)$ such that~$\thetab_n\in\mathcal{S}^{p_n-1}$ for any~$n$, and a sequence~$(F_n)$ of cumulative
\vspace{-.6mm}
  distribution functions over~$[-1,1]$. In this framework, the spherical location problem
  \vspace{-.5mm}
 consists in testing
$
\mathcal{H}\n_{0}:\thetab_n=\thetab_{n0}
$
against
%\qquad
%\textrm{ against }
%\qquad
$\mathcal{H}\n_{1}:\thetab_n\neq \thetab_{n0}
$,
where~$(\thetab_{n0})$ is a fixed null parameter sequence. Parallel to the notation that was used in the FvML case, we will write~$e_{n\ell}$ and~$\tilde{e}_{n\ell}$, $\ell=1,2,\ldots$ for the non-central 
\vspace{-.6mm}
 and central moments of~$F_n$, respectively. These are the moments, under~${\rm P}\n_{\thetab_n,F_n}$, of the quantity~$U_{n1}=\Xb_{n1}\pr\thetab_n$ in the tangent-normal decomposition of~$\Xb_{n1}$ with respect to~$\thetab_n$. The corresponding non-central moments of $V_{n1}=\sqrt{1-U_{n1}^2}$ will still be denoted as~$f_{n\ell}$. 

%In this section, we investigate the non-null asymptotic properties of the high-dimensional Watson test under suitable ``local" alternatives. These alternatives will be of the form~${\rm P}\n_{\thetab_{n},F_n}$, where~$\thetab_{n}=\thetab_{n0}+\nu_n \taub_n$, where the sequence~$(\nu_n)$ in~$\R^+_0$ and the bounded sequence~$(\taub_n)$ in~$\R^{p_n}$, are such that~$\thetab_{n}\in\mathcal{S}^{p_n-1}$ for any~$n$, and where~$(F_n)$ is still a sequence of cumulative distribution functions over~$[-1,1]$. The constraint~$\thetab_{n}\in\mathcal{S}^{p_n-1}$ imposes that
%\begin{equation}
%\label{constraintnun}
%	\thetab_{n0}\pr\taub_n
%	=
%	-\frac{1}{2} \nu_n \|\taub_n\|^2 
%\end{equation}
%for any~$n$. 

%%%%%%%%%%%%

%\newpage

Using the notation~$V_{ni}$ and~$\Sb_{ni}$ from the tangent-normal decomposition of~$\Xb_{ni}$ with respect to the null location~$\thetab_{n0}$, the Watson test statistic rewrites
$$
	\widetilde{W}_n  
	=
	\frac{W_{n}-(p_n-1)}{\sqrt{2(p_n-1)}}
	=
	\frac{\sqrt{2(p_n-1)}}{\sum_{i=1}^n V_{ni}^2}
	\sum_{1\leq i<j\leq n} V_{ni} V_{nj} \Sb_{ni}\pr\Sb_{nj}
,
$$
where~$W_{n}$ denotes the Watson test statistic in~(\ref{WatsonLD}) based on the null location~$%\thetab_0=
\thetab_{n0}$. Under the null and under appropriate local alternatives, it is expected that~$\widetilde{W}_n$ is asymptotically equivalent in probability to 
$$
	W^*_n
	:=
	\frac{\sqrt{2(p_n-1)}}{n f_{n2}}
	\sum_{1\leq i<j\leq n} V_{ni} V_{nj} \Sb_{ni}\pr\Sb_{nj}
,
$$
so that an important step in the investigation of the non-null properties of~$\widetilde{W}_n$ is the study of the non-null behavior of~$W^*_n$. A classical martingale central limit theorem (see, e.g., Theorem~35.12 in \cite{Bil1995}) provides the following result.

\begin{theorem}
	\label{TheoBillingH1}
Let~$(p_n)$ be a sequence of integers that diverges to infinity and~$(\thetab_{n0})$ be a sequence such that $\thetab_{n0}$ belongs to~$\mathcal{S}^{p_n-1}$ for any~$n$. 
Let~$(F_n)$ be a sequence of cumulative distribution functions on~$[-1,1]$ such that~(a) $f_{n2}>0$ for any~$n$,~(b) $f_{n4}/f_{n2}^2=o(n)$ and~(c)~$\sqrt{p_n}e_{n2}=o(1)$.  
Then, we have the following, where, in each case, $(\taub_n)$ refers to an arbitrary sequence such that~$\thetab_n=\thetab_{n0}+\nu_n\taub_n$ belongs to~$\mathcal{S}^{p_n-1}$ for any~$n$ and such that~$(\|\taub_n\|)$ converges to~$t(\in[0,\infty)):$
	\begin{enumerate}
		\item[] (i)--(iii) if (i)~$\sqrt{n}e_{n1}\to \infty$, if (ii)~$\sqrt{n}e_{n1}\to\xi>0$, or if (iii)~$\sqrt{n}e_{n1}
		\linebreak
		\vspace{-1mm}
		\to 0$ with~$\sqrt{n}p_n^{1/4}e_{n1}\to \infty$,
		then
$$
	W^*_n 
	\stackrel{\mathcal{D}}{\longrightarrow}
	\mathcal{N}\bigg(
	\frac{t^2}{\sqrt{2}}
	,
	1
	\bigg)
$$
	under~${\rm P}\n_{\thetab_{n0}+\nu_n \taub_n,F_n}$, with~$\nu_n=\sqrt{f_{n2}}/(\sqrt{n}p_n^{1/4}e_{n1})$; in cases~(i)--(ii), the constraint~(c) above is superfluous;
	\vspace{3mm}
		\item[] (iv) if~$\sqrt{n}p_n^{1/4}e_{n1}\to \xi>0$, 
		then
$$
	W^*_n 
	\stackrel{\mathcal{D}}{\longrightarrow}
	\mathcal{N}\bigg(
	\frac{\xi^2t^2}{\sqrt{2}}
	\big(1-\frac{t^2}{4} \big)
	,1
	\bigg)
$$
under~${\rm P}\n_{\thetab_{n0}+\nu_n \taub_n,F_n}$, with~$\nu_n=1$;
	\vspace{3mm}
		\item[] (v) if~$\sqrt{n}p_n^{1/4}e_{n1}=o(1)$, then
			$$
	W^*_n 
	\stackrel{\mathcal{D}}{\longrightarrow}
	\mathcal{N}(0,1)
	$$
	under~${\rm P}\n_{\thetab_{n0}+\nu_n \taub_n,F_n}$, with~$\nu_n=1$. 
	\end{enumerate}
\end{theorem}
\vspace{3mm}

To obtain the corresponding non-null results for the Watson test statistic~$\widetilde{W}_n$, we need to prove that~$\widetilde{W}_n$ and~$W^*_n$ are indeed asymptotically equivalent in probability. The following result does so in the, possibly non-null, general rotationally symmetric context considered (in the \emph{FvML} case, the \emph{null} version of this result was established when proving the results of Section~\ref{secLAN}; see the proof of Lemma~\ref{LemJointNormality}).

\begin{theorem}
	\label{TheoGeneralSlutzky}
Let~$(p_n)$ be a sequence of integers that diverges to infinity and~$(\thetab_{n0})$ be a sequence such that $\thetab_{n0}$ belongs to~$\mathcal{S}^{p_n-1}$ for any~$n$. 
 Let~$(F_n)$ be a sequence of cumulative distribution functions on~$[-1,1]$ such that~(a) $f_{n2}>0$ for any~$n$ and (b) $f_{n4}/f_{n2}^2=o(n)$. Then, with~$(\nu_n)$ and~$(\taub_n)$ as in Theorem~\ref{TheoBillingH1}, we have that, in each regime~(i)--(v) considered there, 
	$$
	\widetilde{W}_n 
	=
	W^*_n + o_{\rm P}(1)
$$
	as~$n\to\infty$ under~${\rm P}\n_{\thetab_{n0}+\nu_n \taub_n,F_n}\!$.
	\end{theorem}
\vspace{1mm}
 
Of course, Theorem~\ref{TheoGeneralSlutzky} readily implies that Theorem~\ref{TheoBillingH1} still holds if one substitutes~$\widetilde{W}_n$ for~$W_n^*$. Rather than restating the result explicitly, we present the following corollary, which focuses on the FvML case.

\begin{corollary}
	\label{CorolBillingH1FvML}
Let~$(p_n)$ be a sequence of integers that diverges to infinity,~$(\kappa_n)$ be a sequence in~$(0,\infty)$, and~$(\thetab_{n0})$ be a sequence such that $\thetab_{n0}$ belongs to~$\mathcal{S}^{p_n-1}$ for any~$n$. 
Then, we have the following, where in each case~$(\taub_n)$ refers to an arbitrary sequence such that~$\thetab_n=\thetab_{n0}+\nu_n\taub_n$ belongs to~$\mathcal{S}^{p_n-1}$ for any~$n$ and such that~$(\|\taub_n\|)$ converges to~$t(\in[0,\infty))\!:$
	\begin{enumerate}
		\item[] (i) if~$\sqrt{n}\kappa_n/p_n^{3/4}\to \infty$, then 
$$
	W^*_n 
	\stackrel{\mathcal{D}}{\longrightarrow}
	\mathcal{N}\bigg(
	\frac{t^2}{\sqrt{2}}
	,1
	\bigg)
$$
	under~${\rm P}\n_{\thetab_{n0}+\nu_n \taub_n,\kappa_n}$, with~$\nu_n=p_n^{3/4}/(\sqrt{n}\kappa_n \sqrt{f_{n2}}\,)$; 
	\vspace{3mm}
		\item[] (ii) if~$\sqrt{n}\kappa_n/p_n^{3/4} \to \xi>0$, 
		then
$$
	W^*_n 
	\stackrel{\mathcal{D}}{\longrightarrow}
	\mathcal{N}\bigg(
	\frac{\xi^2t^2}{\sqrt{2}} \Big(1-\frac{t^2}{4}\Big)
	,1
	\bigg)
$$
under~${\rm P}\n_{\thetab_{n0}+\nu_n \taub_n,\kappa_n}$, with~$\nu_n=1$;
	\vspace{3mm}
		\item[] (iii) if~$\sqrt{n}\kappa_n/p_n^{3/4} = o(1)$, then
			$$
	W^*_n  
	\stackrel{\mathcal{D}}{\longrightarrow}
	\mathcal{N}(0,1)
	$$
	under~${\rm P}\n_{\thetab_{n0}+\nu_n \taub_n,\kappa_n}$, with~$\nu_n=1$. 
	\end{enumerate}
\end{corollary} 
\vspace{3mm}

It is interesting to comment on how this relates to the results of the previous section: Corollary~\ref{CorolBillingH1FvML}(i) covers the regimes (i)--(iv) %(Theorems~\ref{TheoLAN}--\ref{LANunspec}) 
and~(v$_a$).
% (Theorem~\ref{LANunspecregime5}). 
In view of the asymptotic behavior of~$f_{n2}$ in these regimes (see Lemma~\ref{LemEnFvML}), Corollary~\ref{CorolBillingH1FvML}(i) confirms the consistency rates of the Watson test in Theorem~\ref{TheoLAN}--\ref{LANunspecregime5}, as well as the corresponding asymptotic powers obtained in~(\ref{asymppoweroptimalitoiii}) through the Le Cam third lemma. 
Corollary~\ref{CorolBillingH1FvML}(ii) relates to regime~(v$_b$),
% (Theorem~\ref{LANunspecregime5}), 
where the Watson test can only see the ``fixed" alternatives associated with~$\nu_n= 1$, with limiting power
\begin{equation}
\label{nonmonotonicpattern}
1-\Phi\bigg(
\Phi^{-1}(1-\alpha) 
-
	\frac{\xi^2t^2}{\sqrt{2}} \Big(1-\frac{t^2}{4}\Big)
\bigg) 
\end{equation}
(note that this limiting power can be obtained both by using Corollary~\ref{CorolBillingH1FvML}(ii) or by applying the Le Cam third lemma in Theorem~\ref{LANunspecregime5}, {even if the second approach will provide the result only for alternatives associated with~$t<\sqrt{2}$, that is, for alternatives in the open hemisphere centered at the null location}). 

The limiting power in~(\ref{nonmonotonicpattern}) increases monotonically from the nominal level~$\alpha$ (for~$t=0$, where the underlying location is the null one) to its maximal value (achieved at~$t=\sqrt{2}$, that is, when the true location is orthogonal to the null one), then decreases monotonically to~$\alpha$ (this limiting value being obtained when the true location is antipodal to the null location). This non-monotonic pattern of the asymptotic power in this regime is a direct consequence of the nature of the Watson test that, as already mentioned, rejects the null when~$\bar{\Xb}_{n}$ and~$\thetab_{n0}$ project
\vspace{-.3mm}
  far from each other onto the orthogonal complement to~$\thetab_{n0}$ in~$\R^{p_n}$. Finally, Corollary~\ref{CorolBillingH1FvML}(iii) indicates that, for~$\sqrt{n}\kappa_n/p_n^{3/4} = o(1)$, there are no alternatives under which the Watson test can show asymptotic powers larger than the nominal level~$\alpha$, which is perfectly in line with results obtained in the previous section for the corresponding regimes, namely for regimes~(v$_c$), (vi) and~(vii).

%%%%%%%%%%%%%%%%%%%

%
%
%
%
%\begin{equation}
%\label{hsl}
%\nu_n^2
%=
%\frac{1}{n^{1/2}p_n^{1/2} e_{n1}}
%(=o(1))
%,
%\end{equation} 
%we have ${\rm E}[R_{n3}]=o(1)$ and ${\rm Var}[R_{n1}]=o(1)={\rm Var}[R_{n3}]$, so that~$R_n=R_{n0}+o_{\rm P}(1)$. Full agreement with regime~(5) for the FvML. Here, if one takes
%$$
%\nu_n^2 
%=
%O\Big(\frac{p_n^{1/2}}{n \kappa_n e_{n1}}\Big)
%=
%O\Big(\frac{p_n^{1/2}}{n p_n e_{n1}^2}\Big)
%=
%O\Big(\frac{1}{n p_n^{1/2} e_{n1}^2}\Big)
%,
%$$
%which are more severe than those in~(\ref{hsl}), then Watson shows non-trivial asymptotic powers. 
%
%(d) $n^{1/2} e_{n1}\to 0$ with $n^{1/2}p_n^{1/2} e_{n1}=O(1)$. 
%$$
%\nu_n
%=
%1
%,
%$$ 
%we have ${\rm E}[R_{n3}]=o(1)$ and ${\rm Var}[R_{n1}]=o(1)={\rm Var}[R_{n3}]$, so that~$R_n=R_{n0}+o_{\rm P}(1)$. Full agreement with regimes~(6)--(7) for the FvML. In this case, of course, no more severe alternatives that could be seen by Watson can be provided. 

%%%%%%%%%%%%%%%%%%%%%%%%%%%%%%%%%%%%%%%%%%%%%%

\section{Simulations}
\label{secSimu} 

This section reports the results of a Monte Carlo study we conducted to see how well the finite-sample behavior of the various tests reflect the  asymptotic findings in Theorems~\ref{TheoLAN}--\ref{LANunspecregime5} and Corollary~\ref{CorolBillingH1FvML}. To compare the results for different values of~$p/n$ (note that most aforementioned asymptotic findings allow~$p_n$ to go to infinity at an arbitrary rate), we conducted three simulations, for~$(n,p)=(800,200)$, $(n,p)=(400,400)$, and~$(n,p)=(200,800)$, respectively. In each simulation, we generated, for every combination of~$r=(i),\ldots,(iv),(v_a),(v_b),(vi),(vii)$ and~$\ell=0,1,\ldots,L=5$, a collection of~$M=1,000$ independent random samples of size~$n$ from the~$p$-variate FvML distribution with location
\begin{eqnarray*}
\thetab_{n,r,\ell}
&:=\,&
%\Big(
%1- \frac{2\nu^2_{n,r}  \ell^2}{L^2},\frac{2\nu_{n,r} \ell}{L}
%\sqrt{1 -\frac{\nu_{n,r}^2  \ell^2}{L^2}},0,\ldots,0
%\Big)' 
%\\[2mm]
%&=& 
(1,0,\ldots,0)'
+
\nu_{n,r}
\bigg(
- \frac{2\nu_{n,r} \ell^2}{L^2},\frac{2 \ell}{L}
\Big(1 -\frac{\nu_{n,r}^2 \ell^2}{L^2}\Big)^{1/2},0,\ldots,0
\bigg)' 
\\[2mm]
&\,=:&
\thetab_{n0}
+
\nu_{n,r} \taub_{n,r,\ell}
\in\mathcal{S}^{p_n-1}
\end{eqnarray*}
and concentration~$\kappa_{n,r}$. 
The index~$r$ allows to consider the various regimes from Theorem~\ref{TheoLAN} (associated with the~$\kappa_{n,r}$ used). In each case, we considered the corresponding local alternatives (associated with~$\nu_{n,r}$) from the same theorem. More precisely, we used
\begin{itemize}
	\item[$\bullet$] $\kappa_{n,(i)}:=p_n^2$, $\nu_{n,(i)}=p_n^{1/4}/\sqrt{n\kappa_n}$,
	\item[$\bullet$] $\kappa_{n,(ii)}:=p_n$, $\nu_{n,(ii)}=\sqrt{c_1}p_n^{3/4}/(\sqrt{n}\kappa_n)$,  
	\item[$\bullet$] $\kappa_{n,(iii)}:=p_n/n^{1/4}$, $\nu_{n,(iii)}=p_n^{3/4}/(\sqrt{n}\kappa_n)$, 
	\item[$\bullet$] $\kappa_{n,(iv)}:=p_n/\sqrt{n}$, $\nu_{n,(iv)}=p_n^{3/4}/(\sqrt{n}\kappa_n)$,  
	\item[$\bullet$] $\kappa_{n,(v_a)}:=p_n^{7/8}/\sqrt{n}$, $\nu_{n,(v_a)}=p_n^{1/4}/(n^{1/4}\sqrt{\kappa_n})$, 
	\item[$\bullet$] $\kappa_{n,(v_b)}:=p_n^{3/4}/\sqrt{n}$, $\nu_{n,(v_b)}=p_n^{1/4}/(n^{1/4}\sqrt{\kappa_n})$,  
	\item[$\bullet$] $\kappa_{n,(vi)}:=\sqrt{p_n}/\sqrt{n}$, $\nu_{n,(vi)}=1$, and 
	\item[$\bullet$] $\kappa_{n,(vii)}:=p_n^{1/4}/\sqrt{n}$, $\nu_{n,(vii)}=1$.    
\end{itemize}
The value~$\ell=0$ corresponds to the null hypothesis~$\mathcal{H}_{0}\n:\thetab_n=\thetab_{n0}$, whereas the values~$\ell=1,\ldots,5$ provide increasingly severe alternatives. For each sample, we performed three tests, all at asymptotic level~$\alpha=5\%$, namely 
(a) the Watson test rejecting the null when
$$
W_{n}
=
\frac{n(p-1)\bar{\Xb}_n\pr (\mathbf{I}_{p}-\thetab_{n0}\thetab_{n0}\pr) \bar{\Xb}_n}{1-\frac{1}{n} \sum_{i=1}^n (\Xb_{ni}\pr\thetab_{n0})^2} 
>
\chi^2_{p-1,1-\alpha}
,
$$
(b) the $Z_n$-test rejecting the null when 
$$
Z_n=\frac{\sqrt{n}(\bar{\Xb}_{n}\pr \thetab_{n0}-e_{n1})}{\sqrt{\tilde{e}_{n2}}}
<\Phi^{-1}(\alpha)
,
$$
and (c) the hybrid test rejecting the null when 
$$
H_n:=
\bigg(
\frac{\widetilde{W}_n}{\sqrt{2}}
-
\frac{Z_n}{2\xi_n}
\bigg)
\Big/
\sqrt{\frac{1}{2}
+
\frac{1}{4\xi_n^2}}
>\Phi^{-1}(1-\alpha)
,
$$
where~$\xi_n:=\sqrt{n}\kappa_{n}/p_n$ is based on the (unknown) concentration~$\kappa_{n}$ depending on the regime~$r$ at hand. In each regime from Theorem~\ref{TheoLAN}, this hybrid test is clearly expected to behave as the corresponding optimal specified-$\kappa_n$ test. {We stress that the tests~(b)--(c) address the specified-$\kappa_n$ problem only, whereas the Watson test~(a) addresses both the specified-$\kappa_n$ and unspecified-$\kappa_n$ problems.} 

Plots of the resulting rejection frequencies are provided in Figures~\ref{Fig1} to~\ref{Fig3}, for~$(n,p)=(800,200)$, $(n,p)=(400,400)$ and $(n,p)=(200,800)$, respectively. In each case, the asymptotic powers, obtained from~(\ref{powgeneraloptimal})--(\ref{asymppoweroptimaliv}), are also plotted. Clearly, irrespective of the three values of~$p/n$ considered, the rejection frequencies of the tests are in an excellent agreement with the corresponding asymptotic powers. Also, the results confirm the adaptive nature of the hybrid test, that throughout is the most powerful test.

\begin{figure}[htbp!]
\hspace*{-3pt}\includegraphics[width=\textwidth]{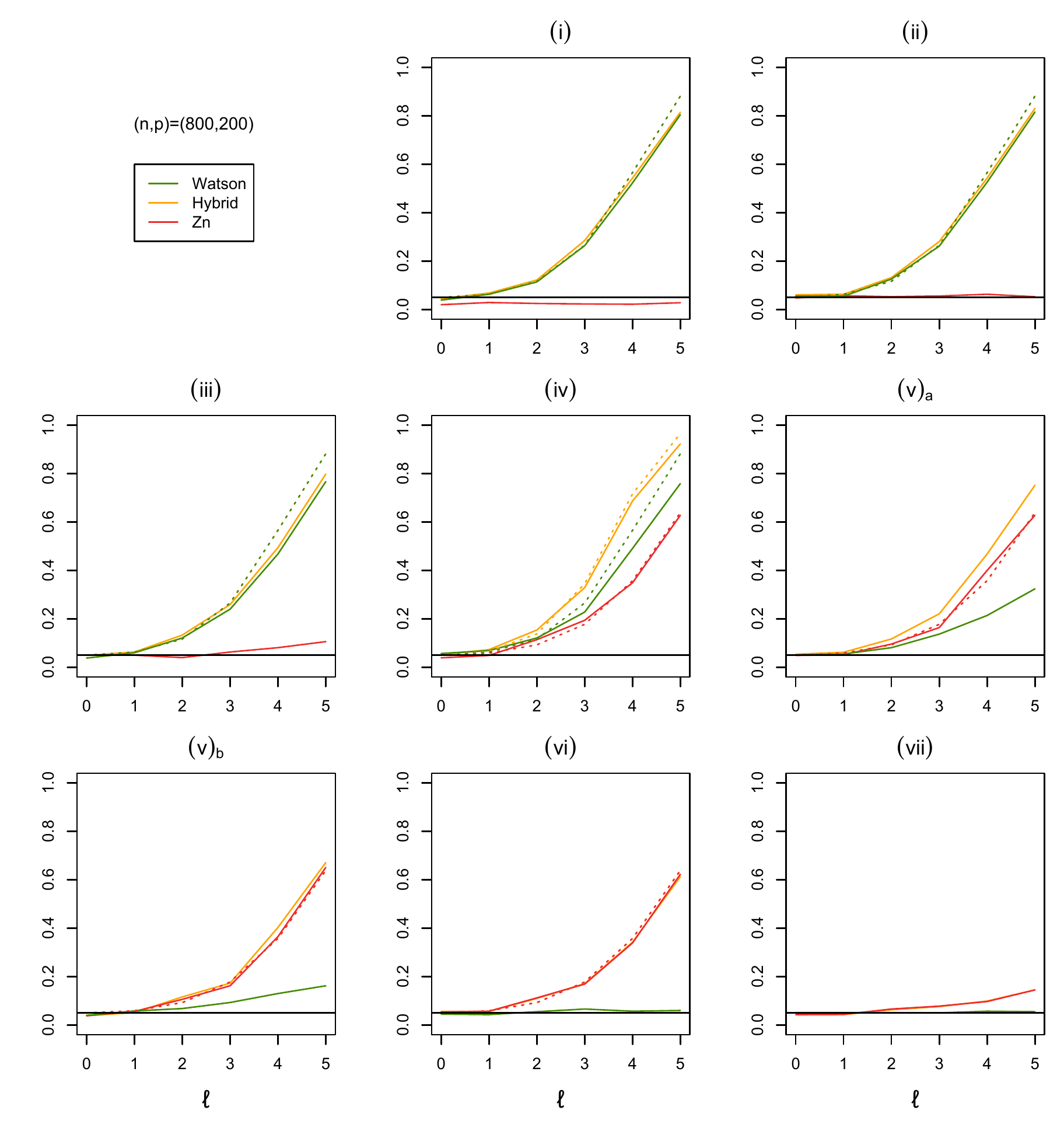}    
\vspace{-1mm}
 \caption{Rejection frequencies
 \vspace{-.5mm}
  (solid lines), 
out of $M=1,\!000$ independent replications, of the Watson test (green), the hybrid test (orange)
 and the $Z_n$-based test (red) for $\mathcal{H}_{0}\n:\thetab_n=\thetab_{n0}=(1,0,\ldots,0)'\in\R^p$, under the null $(\ell=0)$ and under increasingly severe $p$-dimensional FvML alternatives $(\ell=1,\ldots,5)$; here, the sample size is~$n=800$ and the dimension is~$p=200$. The regimes~$(i),\ldots,(vii)$ fix the way the underlying concentration~$\kappa_n$ is chosen as a function of~$n$ and~$p$. In each regime, the corresponding contiguous alternatives from Theorem~\ref{TheoLAN} are used; see Section~\ref{secSimu} for details. The corresponding asymptotic powers are plotted in each case (dashed lines).  
}
\label{Fig1}   
\end{figure}

%%%%%%%%%%%%%%%%%%%%%%%%%%%%%%%%%%%%%%%%%%%%%%

\begin{figure}[htbp!]
\hspace*{-3pt}\includegraphics[width=\textwidth]{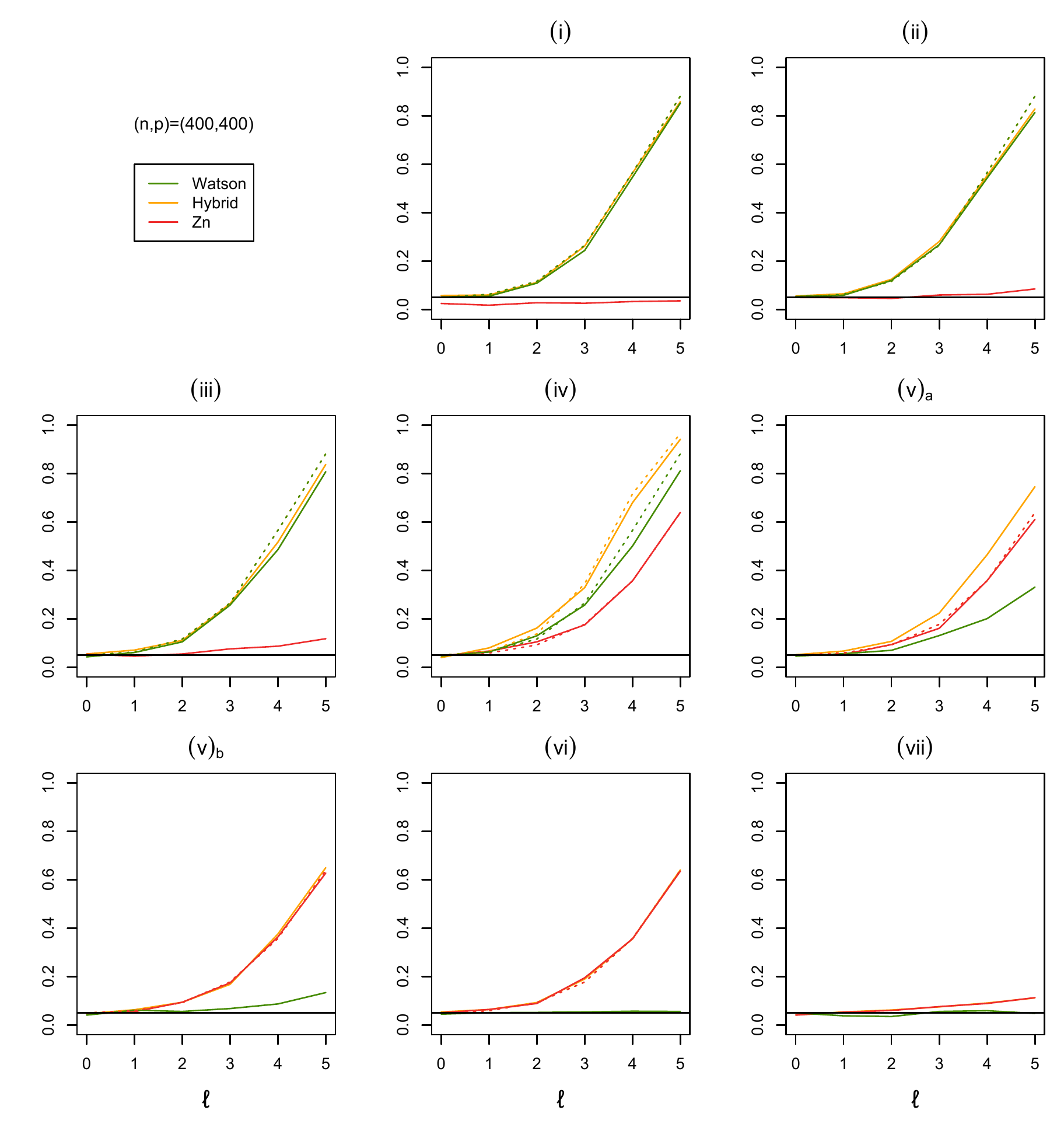}    
\vspace{-1mm}
 \caption{Same results as in Figure~\ref{Fig1}, but for sample size~$n=400$ and dimension~$p=400$. 
}
\label{Fig2}  
\end{figure}
 
%%%%%%%%%%%%%%%%%%%%%%%%%%%%%%%%%%%%%%%%%%%%%% 

\begin{figure}[htbp!]
\hspace*{-3pt}\includegraphics[width=\textwidth]{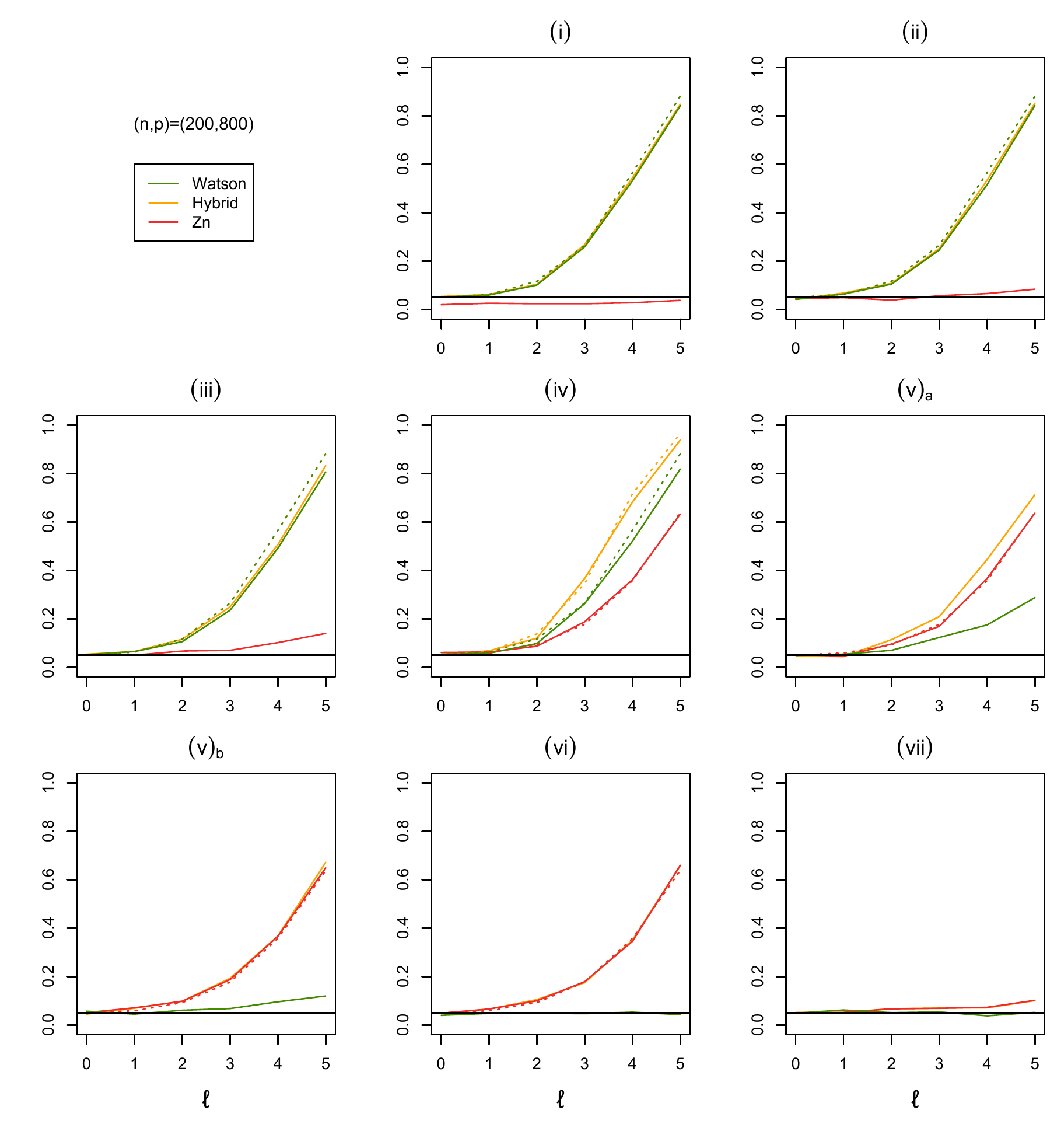}    
\vspace{-1mm}
 \caption{Same results as in Figures~\ref{Fig1}--\ref{Fig2}, but for sample size~$n=200$ and dimension~$p=800$. 
}
\label{Fig3}  
\end{figure}

%%%%%%%%%%%%%%%%%%%%%%%%%%%%%%%%%%%%%%%%%%%%%%

To illustrate similarly the results of Theorem~\ref{LANunspecregime5} and Corollary~\ref{CorolBillingH1FvML}, we focused on the regimes (v$_a$)--(v$_b$) above, but considered the corresponding more severe alternatives. More precisely, we here took  
\begin{itemize}
	\item[$\bullet$] $\kappa_{n,(v_a)}:=p_n^{7/8}/\sqrt{n}$, $\nu_{n,(v_a)}=p_n^{3/4}/(\sqrt{n}\kappa_n)$, and
	\item[$\bullet$] $\kappa_{n,(v_b)}:=p_n^{3/4}/\sqrt{n}$, $\nu_{n,(v_b)}=1$.   
\end{itemize}
The rejection frequencies of the same three tests as above, still based on $M=1,\!000$ independent replications, are provided in Figure~\ref{Fig4}. For the Watson test, the agreement between rejection frequencies and asymptotic powers is perfect in regime~(v$_b$) (where the non-monotonic asymptotic power pattern is confirmed), but is less so in regime~(v$_a$); at the finite dimensions / sample sizes considered, this may be explained by the fact that the regimes~(v$_a$)--(v$_b$) are close to each other, so that the empirical powers of the Watson test in regime~(v$_a$) tends to be pulled to the ones in regime~$(vi)$.

\begin{figure}[htbp!]
%\hspace*{-3pt}
\begin{center}
	\includegraphics[width=.85\textwidth]{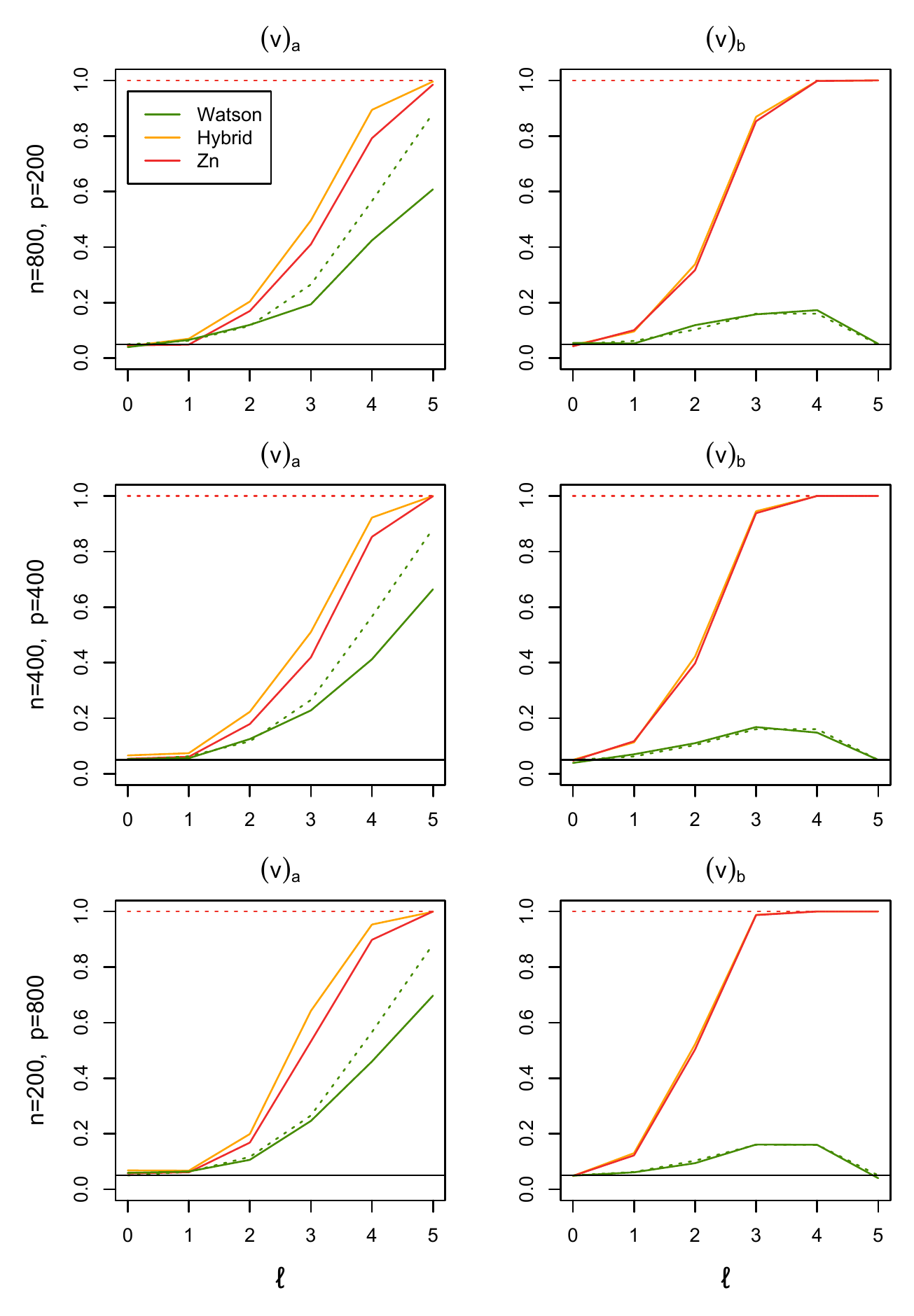}   
\end{center}
\vspace{-1mm}
 \caption{Rejection frequencies
 \vspace{-.5mm}
  (solid lines), out of $M=1,\!000$ independent replications, of the Watson test (green), the hybrid test (orange) and the $Z_n$-based test (red) for $\mathcal{H}_{0}\n:\thetab_n=\thetab_{n0}=(1,0,\ldots,0)'\in\R^p$, under the null $(\ell=0)$ and under increasingly severe $p$-dimensional FvML alternatives $(\ell=1,\ldots,5)$; the couples~$(n,p)$ used are those from Figures~\ref{Fig1}--\ref{Fig3}. Here, we focus on the regimes~(v$_a$)--(v$_b$) and consider the more severe alternatives associated with Theorem~\ref{LANunspecregime5} and Corollary~\ref{CorolBillingH1FvML}; see Section~\ref{secSimu} for details.  The corresponding asymptotic powers are plotted in each case (dashed lines).  
}
\label{Fig4}  
\end{figure}

%%%%%%%%%%%%%%%%%%%%%%%%%%%%%%%%%%%%%%%%%%%%%%
 
%\clearpage

\section{Summary and research perspectives}
\label{secFinaComments} 

In the present paper, we tackled the problem of testing, in high dimensions, the null hypothesis that the spike direction~$\thetab$ of a rotationally symmetric distribution is equal to a given direction~$\thetab_0$. Under FvML distributional assumptions, we showed that, after resorting to the invariance principle, the sequence of statistical experiments at hand is LAN. More precisely, we identified seven regimes, according to the way the underlying concentration parameter~$\kappa_n$ depends on~$n$ and~$p_n$, each leading to a specific limiting experiment, with its own central sequence, Fisher information and contiguity rate (interestingly, these heterogeneous contiguity rates precisely quantify how difficult the problem gets for low concentration situations). As a result, the Le Cam optimal test (more precisely, the locally asymptotically most powerful invariant test) depends on the regime considered. In regimes where~$\sqrt{n}\kappa_n/p_n\to\infty$, the classical Watson test is optimal, whereas in regimes where~$\sqrt{n}\kappa_n/p_n=O(1)$, the optimal test is an oracle test that explicitly involves the unknown value of the underlying concentration~$\kappa_n$. If~$\sqrt{n}\kappa_n/p_n\to \xi>0$, then the Watson test fails to be optimal but is still rate-consistent, whereas if~$\sqrt{n}\kappa_n/p_n=o(1)$, then it is not even rate-consistent. In all cases, we obtained from the Le Cam third Lemma the asymptotic powers of the corresponding optimal tests and of the Watson test under contiguous alternatives. All results above allow the dimension~$p_n$ to go to infinity arbitrarily slowly or arbitrarily fast as a function of~$n$, hence cover moderately high dimensions as well as ultra-high dimensions.  

Optimality above refers to the specified-$\kappa_n$ version of the testing problem considered. Since the concentration~$\kappa_n$ can hardly be assumed to be known in practice, however, optimality results for the corresponding unspecified-$\kappa_n$ problem are more relevant. For this problem, the Watson test of course remains optimal in regimes where~$\sqrt{n}\kappa_n/p_n\to\infty$. But remarkably, for \mbox{unspecified}~$\kappa_n$, the Watson test is also optimal in regimes where~$\sqrt{n}\kappa_n/p_n=O(1)$, sometimes under the condition that~$p_n=o(n^2)$ (on an even weaker condition on~$p_n$); we refer to Table~\ref{thethetable} and to Theorems~\ref{LANunspec}--\ref{LANunspecregime5} for details.  

Our work opens several perspectives for future research. (a) First, while we derived non-null results for the Watson test also outside the FvML distributional setup, all our optimality results are limited to the FvML case. A natural question is therefore whether or not the strong optimality properties of the Watson test extend away from the FvML case. The low-dimensional investigation conducted in \cite{PaiVer17b} leads us to conjecture that optimality would also hold away from the FvML case, at least in low concentration patterns. Establishing this would require expanding invariant log-likelihood ratios taking a much more complicated form than in the FvML case. This calls for entirely different techniques, hence is beyond the scope of the present paper. (b) Second, we would like to mention that our results are also relevant in a Euclidean (i.e., non-directional) context. They indeed characterize the asymptotic efficiency of sign tests for the direction~$\thetab$ of a skewed single-spiked distribution in~$\R^{p}$, that is, a distribution whose projection along~$\thetab$ is skewed and whose projection onto the orthogonal complement to~$\thetab$ is spherically symmetric. This skewed version of the corresponding classical, elliptical, problem is natural in a signal detection framework, where the signal at hand is quite naturally maximal in direction~$\thetab$ and minimal in the opposite direction~$-\thetab$. While our results exhaustively address the question of efficiency of sign tests for this problem (that is, of tests that involve the observations only through their direction form the center of the distribution), it would be of interest to also consider the efficiency of more general testing procedures.

%%%%%%%%%%%%%%%%%%%%%%%%%%

\appendix

\section{Technical proofs for Section~2}
\label{SecSupProofLAN}

%%%%%%%%%%%%%%%%%%%%%%%%%%
%%%%%%%%%%%%%%%%%%%%%%%%%%
%%%%%%%%%%%%%%%%%%%%%%%%%%

The proof of  Theorem~\ref{TheoLAQ} requires the following preliminary results.

\begin{lemma}
\label{LemSlutzkyFvML}
Let~$(p_n)$ be a sequence of integers diverging to infinity and~$(\kappa_n)$ be an arbitrary sequence in~$(0,\infty)$. Let~$L_n:=\sum_{i=1}^n V_{ni}^2/(n f_{n2})$, where we used 
\vspace{-.3mm}
 the notation~$V_{ni}=(1-(\Xb_{ni}'\thetab_{n0})^2)^{1/2}$ and~$f_{n2}={\rm E}[V_{n1}^2]$.
Then, 
${\rm E}
\big[
(
L_n
-
1
)^2
\big]
=
o(p_n^{-1})
$
as~$n\to \infty$ under~${\rm P}\n_{\thetab_{n0},\kappa_n}\!$. 
\end{lemma}

{\sc Proof of Lemma~\ref{LemSlutzkyFvML}.}
Since
\begin{eqnarray*}
\lefteqn{
\hspace{-13mm}
{\rm E}
\Bigg[
\bigg(
\frac{\sum_{i=1}^n V_{ni}^2}{n f_{n2}}
-
1 
\bigg)^2
\Bigg]
=
\frac{1}{f_{n2}^2}
\,
{\rm E}
\Bigg[
\Bigg(
\frac{1}{n} \sum_{i=1}^{n} V_{ni}^2 - {\rm E}[V_{n1}^2]
\Bigg)^2
\Bigg]
}
\nonumber
\\[2mm]
& &
\hspace{5mm}
=
\frac{1}{f_{n2}^2}
\,
{\rm Var}\Bigg[\frac{1}{n} \sum_{i=1}^{n} V_{ni}^2\Bigg]
=
\frac{{\rm Var}[V_{n1}^2]}{nf_{n2}^2}
=
\frac{f_{n4}-f_{n2}^2}{nf_{n2}^2}
\end{eqnarray*}
(recall that~$f_{n4}:={\rm E}[V_{n1}^4]$), it is sufficient to prove that 
\begin{equation}
\frac{f_{n4}-f_{n2}^2}{f_{n2}^2}
=
O(p_n^{-1})
.
\end{equation}
Now, the expression for~$f_{n4}/f_{n2}^2
%={\rm E}[V_{n1}^{4}]/({\rm E}[V_{n1}^{2}])^2
$ in page~82 of \cite{Leyetal2015} yields
\begin{eqnarray*}
\bigg|
\frac{f_{n4}-f_{n2}^2}{f_{n2}^2}
	\bigg|
&=&
\bigg|
\frac{(p_n+1)\mathcal{I}_{\frac{p_n}{2}+1}(\kappa_n)\mathcal{I}_{\frac{p_n}{2}-1}(\kappa_n)}{(p_n-1)(\mathcal{I}_{\frac{p_n}{2}}(\kappa_n))^2}
-
1
\bigg|
\\[2mm]
&=&
\bigg|\frac{(p_n+1)(\mathcal{I}_{\frac{p_n}{2}+1}(\kappa_n)\mathcal{I}_{\frac{p_n}{2}-1}(\kappa_n)-(\mathcal{I}_{\frac{p_n}{2}}(\kappa_n))^2)}{(p_n-1)(\mathcal{I}_{\frac{p_n}{2}}(\kappa_n))^2}
+
\frac{2}{p_n-1}
\bigg|
\\[2mm]
&\leq&
\frac{3|\mathcal{I}_{\frac{p_n}{2}+1}(\kappa_n)\mathcal{I}_{\frac{p_n}{2}-1}(\kappa_n)-(\mathcal{I}_{\frac{p_n}{2}}(\kappa_n))^2|}{(\mathcal{I}_{\frac{p_n}{2}}(\kappa_n))^2}
+
\frac{2}{p_n-1}
\cdot
\end{eqnarray*}
Since~$|\mathcal{I}_{\frac{p_n}{2}+1}(\kappa_n)\mathcal{I}_{\frac{p_n}{2}-1}(\kappa_n)-(\mathcal{I}_{\frac{p_n}{2}}(\kappa_n))^2|\leq (\mathcal{I}_{\frac{p_n}{2}}(\kappa_n))^2/(\frac{p_n}{2}+1)$ (see~(3.1)--(3.2) in \cite{JosBis1991}), the result follows.
\cqfd
\vspace{1mm}

%%%%%%%%%%%%%%%%%%%%%%%%%%%%%%%%%%%%%%%%%%%%%%%%%%%%%%

\begin{lemma}
\label{LemJointNormality} 
Let~$(p_n)$ be a sequence of integers that diverges to infinity and~$(\kappa_n)$ be an arbitrary sequence in~$(0,\infty)$. Let~$(\thetab_{n0})$ be a sequence such that~$\thetab_{n0}$ belongs to~$\mathcal{S}^{p_n-1}$ for any~$n$. Consider the random variables~$\widetilde{W}_n$ and~$Z_n$ introduced in Theorem~\ref{TheoLAQ}. Then, $(\widetilde{W}_n,Z_n)'$ is asymptotically standard bivariate normal under~${\rm P}^{(n)}_{\thetab_{n0},\kappa_n}\!$.  
\end{lemma}

{\sc Proof of Lemma~\ref{LemJointNormality}.}
Throughout the proof, expectations and variances are under~${\rm P}^{(n)}_{\thetab_{n0},\kappa_n}$ and stochastic convergences are as~$n\to\infty$ under the same sequence of hypotheses, whereas $U_{ni}$, $V_{ni}$ and~$\Sb_{ni}$ refer to the tangent-normal decomposition of~$\Xb_{ni}$ with respect to~$\thetab_{n0}$. Letting then 
$$
W^*_n
=
\frac{\sqrt{2(p_n-1)}}{nf_{n2}}
\sum_{1\leq i < j\leq n} V_{ni} V_{nj} \Sb_{ni}\pr\Sb_{nj} 
,
$$
assume that $(W^*_n,Z_n)'$ is asymptotically standard bivariate normal. Then,
%$$
%\widetilde{W}_n
%=
%\frac{W_{n}-(p_n-1)}{\sqrt{2(p_n-1)}}
%=
%\frac{\sqrt{2(p_n-1)}}{\sum_{i=1}^n V_{ni}^2}
%\sum_{1\leq i < j\leq n} V_{ni} V_{nj} \Sb_{ni}\pr\Sb_{nj} 
%$$
\begin{eqnarray}
\widetilde{W}_n
-
W^*_n
&=&
\Bigg[
\frac{\sqrt{2(p_n-1)}}{\sum_{i=1}^n V_{ni}^2}
-
\frac{\sqrt{2(p_n-1)}}{nf_{n2}}
\Bigg]
\sum_{1\leq i < j\leq n} V_{ni} V_{nj} \Sb_{ni}\pr\Sb_{nj} 
\label{differencedeW}
\\[2mm]
&=&
\Bigg[
1
-
\frac{\sum_{i=1}^n V_{ni}^2}{nf_{n2}}
\Bigg]
\times 
\frac{nf_{n2}}{\sum_{i=1}^n V_{ni}^2}
\times
\Bigg(
\frac{\sqrt{2(p_n-1)}}{nf_{n2}}
\sum_{1\leq i < j\leq n} V_{ni} V_{nj} \Sb_{ni}\pr\Sb_{nj} 
\Bigg)
\nonumber
\\[2mm]
&=&
\frac{1-L_n}{L_n} 
\,
W^*_n
,
\nonumber
\end{eqnarray}
where~$L_n$ was introduced in Lemma~\ref{LemSlutzkyFvML}. This lemma implies that~\mbox{$L_n-1$}, hence also~$(1-L_n)/L_n$, is~$o_{\rm P}(1)$. If~$(W^*_n,Z_n)'$ is indeed asymptotically standard bivariate normal, then we conclude that~$\widetilde{W}_n-W^*_n$ is~$o_{\rm P}(1)$, so that~$(\widetilde{W}_n,Z_n)'$ itself is asymptotically standard bivariate normal. 

It is therefore sufficient to show that $(W^*_n,Z_n)'$ is asymptotically standard bivariate normal. We will do this by fixing~$\gamma$ and~$\eta$ such that~$\gamma^2+\eta^2=1$ and by using a classical martingale Central Limit Theorem to show that~$D_n:=\gamma W^*_n + \eta Z_n$ is asymptotically standard normal. 
To do so, let ${\cal F}_{n\ell}$ be the $\sigma$-algebra generated by $\Xb_{n1}, \ldots, \Xb_{n \ell}$ and denote by ${\rm E}_{n\ell}[.]$ the conditional expectation with respect to ${\cal F}_{n\ell}$. Define
$
D_{n\ell}
:=
{\rm E}_{n\ell} [D_n]-{\rm E}_{n,\ell-1} [D_n] 
$
for $\ell=1,\ldots,n$ and $D_{n\ell}=0$ for~$\ell>n$. It is then easy to check that~$D_{n\ell}=\gamma W^*_{n\ell} + \eta Z_{n\ell}$, with  
$$
W^*_{n\ell}
:=
\frac{\sqrt{2(p_n-1)}}{n f_{n2}}
\,
 \sum_{i=1}^{\ell-1} V_{ni} V_{n\ell} \Sb_{ni}\pr\Sb_{n\ell}
\quad \textrm{ and } \quad
Z_{n\ell}
:= 
\frac{U_{n\ell}-e_{n1}}{\sqrt{n\tilde{e}_{n2}}} 
%=
%\frac{\Xb_{n\ell}\pr \thetab_{n0}-e_{n1}}{\sqrt{n\tilde{e}_{n2}}} 
%(\Xb_{n\ell}\pr \thetab_{n0}-e_{n1})
$$
for~$\ell=1,\ldots,n$ and $W^*_{n\ell}=0=Z_{n\ell}$ for~$\ell>n$ ($W^*_{n1}$ is also to be understood as zero). 
%\textcolor{red}{Trivially, both 
%% the $D_{n\ell}$'s have finite second-order moments. 
%$\sum_{\ell=1}^\infty D_{n\ell}$ and $\sum_{\ell=1}^\infty \sigma^2_{n\ell}$, where we let $\sigma^2_{n\ell}:={\rm E}_{n,\ell-1} [D_{n\ell}^2]$, converge with probability~1}. 
To conclude from the martingale Central Limit Theorem in Theorem~35.12 from \cite{Bil1995} that~$D_n=\sum_{\ell=1}^\infty D_{n\ell}$ is indeed asymptotically standard normal, we need to show that
(a) $\sum_{\ell=1}^n \sigma^2_{n\ell}\to 1$ in probability, with~$\sigma^2_{n\ell}:={\rm E}_{n,\ell-1} [D_{n\ell}^2]$, and that 
(b) $
\sum_{\ell=1}^n {\rm E}[D_{n\ell}^2 \, {\mathbb I}[| D_{n\ell}| > \varepsilon]]\to 0
$  
for any~$\varepsilon>0$. 
Clearly, for~$\ell=1,\ldots,n$,    
\begin{eqnarray}
\sigma^2_{n\ell}
&=&
\gamma^2 {\rm E}_{n,\ell-1} [(W^*_{n\ell})^2] 
+ \eta^2 {\rm E}_{n,\ell-1} [Z_{n\ell}^2] 
+ 2\gamma \eta {\rm E}_{n,\ell-1} [W^*_{n\ell} Z_{n\ell}]
\nonumber
\\[2mm]
&=&
\gamma^2 {\rm E}_{n,\ell-1} [(W^*_{n\ell})^2] 
+ \frac{\eta^2}{n}
, 
\label{decompovar}
\end{eqnarray}
so that~(a) follows from Lemma~A.1 in \cite{Leyetal2015}. We may thus focus on~(b). 
%Below, $C$ is a generic constant that may change from line to line. 
Since
$$
{\rm E}_{n,\ell-1}[(W^*_{n\ell})^2] 
%=
%\frac{2(p_n-1)}{n^2 ({\rm E}[v_{n1}^2])^2}  \sum_{i,j=1}^{\ell-1} V_{ni} V_{nj} {\rm E}[v_{n\ell}^2] \Sb_{ni}\pr {\rm E}[\Sb_{n\ell}\Sb_{n\ell}\pr]\Sb_{nj} 
=
2 (n^2 f_{n2})^{-1} \sum_{i,j=1}^{\ell-1} V_{ni} V_{nj} \Sb_{ni}\pr \Sb_{nj} 
,
$$
we obtain~${\rm Var}[D_{n\ell}]
%={\rm E}[D^2_{n\ell}]
={\rm E}[\sigma^2_{n\ell}]
=\gamma^2 {\rm E}[(W^*_{n\ell})^2] + (\eta^2/n)
=2\gamma^2 (\ell-1)/n^2 + (\eta^2/n)
\leq 2/n$, which yields that there exists a constant~$C$ such that, for any~$\varepsilon>0$,   
\begin{eqnarray*}
	\lefteqn{
\hspace{-5mm} 
	\sum_{\ell=1}^n {\rm E}[D_{n\ell}^2 \; {\mathbb I}[| D_{n\ell}| > \varepsilon]]
\leq 
 \sum_{\ell=1}^n 
\sqrt{{\rm E}[D_{n\ell}^4]} 
\,
\sqrt{{\rm P}[|D_{n\ell}| > \varepsilon]} 
}
\\[2mm]
&\leq &
\frac{1}{\varepsilon}
 \sum_{\ell=1}^n 
\sqrt{{\rm E}[D_{n\ell}^4]} 
\,
\sqrt{{\rm Var}[D_{n\ell}]} 
\leq 
\frac{\sqrt{2}}{\sqrt{n}\varepsilon}
 \sum_{\ell=1}^n 
\sqrt{{\rm E}[D_{n\ell}^4]} 
\\[2mm]
&\leq &
\frac{C}{\sqrt{n}\varepsilon}
 \sum_{\ell=1}^n 
\sqrt{{\rm E}[(W^*_{n\ell})^4]}
+
\frac{C}{\sqrt{n}\varepsilon}
 \sum_{\ell=1}^n 
 \sqrt{{\rm E}[Z_{n\ell}^4]} 
.
\end{eqnarray*}
From~(A.9) in \cite{Leyetal2015}, we then obtain
\begin{eqnarray*}
	\sum_{\ell=1}^n {\rm E}[D_{n\ell}^2 \; {\mathbb I}[| D_{n\ell}| > \varepsilon]]
&\leq&
\frac{C}{\sqrt{n}\varepsilon}
 \sum_{\ell=1}^n 
\sqrt{\frac{12}{n^4}
\bigg(
\ell 
\,
\frac{f_{n4}^2}{f_{n2}^4}
+
\ell^2 
\,
\frac{f_{n4}}{f_{n2}^2}
\bigg)
}
+
\frac{C\sqrt{n}}{\varepsilon}
\sqrt{\frac{\tilde{e}_{n4}}{n^2\tilde{e}_{n2}^2}}
\\[2mm]
&\leq&
\frac{\sqrt{12}C}{\varepsilon}
\sqrt{
\Big(\frac{f_{n4}}{nf_{n2}^2}\Big)^2
+
\frac{f_{n4}}{nf_{n2}^2} 
}
+
\frac{C}{\varepsilon}
\sqrt{\frac{\tilde{e}_{n4}}{n\tilde{e}_{n2}^2}}
\cdot
\end{eqnarray*}
The result therefore follows from the fact that
%, irrespective of~$p_n$ and~$\kappa_n$, 
both~$f_{n4}/f_{n2}^2$ and~$\tilde{e}_{n4}/\tilde{e}_{n2}^2$ are upper-bounded by a universal constant; see Theorem~S.2.1 in \cite{PaiVer17asupp}. 
\cqfd
\vspace{1mm}

%%%%%%%%%%%%%%%%%%%%%%%%%

\begin{lemma}
\label{LemGeneralDLFvML}
Let~$(\nu_n)$ be a sequence in~$(0,\infty)$ that diverges to~$\infty$, $(a_n)$, $(b_n)$ be sequences in~$(0,\infty)$ such that~$\liminf a_n>0$, $b_n/\nu_n\to \xi\in[0,\infty)$ and~$b_n^6=o(a_n^4\nu_n^5)$. Let~$T_n$ be a sequence of random variables that is~$O_{\rm P}(1)$. Then, writing,
$$
H_{\nu}(x)
:=
\frac{\int_{-1}^1 (1-t^2)^{\nu-\frac{1}{2}}  \exp(x t)\, dt}{\int_{-1}^1 (1-t^2)^{\nu-\frac{1}{2}}\, dt}
=
\frac{\Gamma(\nu+1)\mathcal{I}_{\nu}(x)}{(x/2)^{\nu}}
%=
%\frac{c_{2(\nu+1)}}{c_{2(\nu+1),x}}
,
$$ 
we have that
$$
a_n^2 \log H_{\nu_n}\Big(\frac{b_n T_n}{a_n}\Big)
=
\frac{b_n^2 T_n^2}{4\nu_n}
-
\frac{b_n^4 T_n^4}{32\nu_n^3a_n^2}
+
\frac{\xi^2 T_n^2}{4}
+
o_{\rm P}(1)
$$
as~$n\to \infty$.
\end{lemma}

{\sc Proof of Lemma~\ref{LemGeneralDLFvML}.}
The proof is based on the bounds
$$
S_{\nu+\frac{1}{2},\nu+\frac{3}{2}}(x) \leq \log H_\nu(x) \leq S_{\nu,\nu+2}(x)
$$
for any~$x>0$, with
$S_{\alpha,\beta}(x) 
:=
%\sqrt{x^2+\beta^2}-\beta-\alpha \log(\alpha+\sqrt{x^2+\beta^2})+\alpha \log(\alpha+\beta)
\sqrt{x^2+\beta^2}-\beta-\alpha \log((\alpha+\sqrt{x^2+\beta^2})/(\alpha+\beta))
$; see~(5) in \cite{HotGru2014}. Consider 
$$
G_\nu(x)
:=
\log H_\nu(x) - \frac{x^2}{4\nu} + \frac{x^4}{32\nu^3}
+\frac{x^2}{4\nu^2}
,
$$
along with its resulting lower and upper bounds
$$
G^{\rm low}_\nu(x)
:=
S_{\nu+\frac{1}{2},\nu+\frac{3}{2}}(x) - \frac{x^2}{4\nu} + \frac{x^4}{32\nu^3} 
+\frac{x^2}{4\nu^2}
$$
and
$$
G^{\rm up}_\nu(x)
:=
S_{\nu,\nu+2}(x) - \frac{x^2}{4\nu} + \frac{x^4}{32\nu^3}
+\frac{x^2}{4\nu^2}
\cdot
$$
We prove the lemma by establishing that
\begin{equation}
	\label{toshowforchristmas}
a_n^2 G^{\rm low/up}_{\nu_n}\Big(\frac{b_n T_n}{a_n}\Big)
=
o_{\rm P}(1)
.
\end{equation}
To do so, we expand the log term in~$G^{\rm low/up}_\nu(x)$ as~$\log x=(x-1)-\frac{1}{2}(x-1)^2+\frac{1}{3c^3}(x-1)^3$ with~$c\in(1,x)$ (note that the argument of these log terms is larger than or equal to one), and we write~$G^{\rm low/up}_\nu(x)=G^{\rm low/up,1}_\nu(x)+	G^{\rm low/up,2}_\nu(x)$, with 
\begin{eqnarray*}
\lefteqn{
G^{\rm low,1}_\nu(x)
:=
	\sqrt{x^2+({\textstyle{\nu+\frac{3}{2}}})^2} - ({\textstyle{\nu+\frac{3}{2}}}) 
}
	\\[2mm]
	& & 
	\hspace{0mm} 
		-({\textstyle{\nu+\frac{1}{2}}}) 
\Bigg[
\Bigg( \frac{({\textstyle{\nu+\frac{1}{2}}})+\sqrt{x^2+({\textstyle{\nu+\frac{3}{2}}})^2}}{2(\nu+1)} -1  \Bigg)
\\[1mm]
	& & 
	\hspace{13mm} 
-\frac{1}{2}
\Bigg( \frac{({\textstyle{\nu+\frac{1}{2}}})+\sqrt{x^2+({\textstyle{\nu+\frac{3}{2}}})^2}}{2(\nu+1)} -1  \Bigg)^2
\Bigg]
%\\[1mm]
%	& & 
%	\hspace{3mm} 
- \frac{x^2}{4\nu} + \frac{x^4}{32\nu^3}
+\frac{x^2}{4\nu^2}
,
\end{eqnarray*}
\begin{eqnarray*}
\lefteqn{
\hspace{-5mm} 
G^{\rm up,1}_\nu(x)
:=
	\sqrt{x^2+(\nu+2)^2} - (\nu+2) 
}
	\\[2mm]
	& & 
	\hspace{0mm} 
		-\nu 
\Bigg[
\Bigg( \frac{\nu+\sqrt{x^2+(\nu+2)^2}}{2(\nu+1)} -1  \Bigg)
\\[1mm]
	& & 
	\hspace{13mm} 
-\frac{1}{2}
\Bigg( \frac{\nu+\sqrt{x^2+(\nu+2)^2}}{2(\nu+1)} -1  \Bigg)^2
\Bigg]
%\\[1mm]
%	& & 
%	\hspace{23mm} 
- \frac{x^2}{4\nu} + \frac{x^4}{32\nu^3}
+\frac{x^2}{4\nu^2}
,
\end{eqnarray*}
$$
%	\hspace{-23mm} 
G^{\rm low,2}_\nu(x)
:=
	- \frac{{\textstyle{\nu+\frac{1}{2}}}}{3(c^{\rm low})^3} 
\Bigg( \frac{({\textstyle{\nu+\frac{1}{2}}})+\sqrt{x^2+({\textstyle{\nu+\frac{3}{2}}})^2}}{2(\nu+1)} -1  \Bigg)^3
,
$$
and
$$
%	\hspace{-38mm} 
G^{\rm up,2}_\nu(x)
:=
	- \frac{\nu}{3(c^{\rm up})^3} 
\Bigg( \frac{\nu+\sqrt{x^2+(\nu+2)^2}}{2(\nu+1)} -1  \Bigg)^3
.
$$

Routine yet tedious computations allow to show that
\begin{eqnarray}
\lefteqn{
G^{\rm low,1}_\nu(x)
=
\frac{x^2}{4\nu^2(\nu+1)}
+
\frac{(4\nu^2+5\nu+2)x^4}{32\nu^3(\nu+1)^2(\nu+2)}
}
\nonumber
\\[2mm]
& & 
\hspace{8mm} 
+\Bigg( 1 - \frac{4(1+\frac{2}{\nu})(1+\frac{3}{2\nu})}{\Big( (1+\frac{3}{2\nu})+\sqrt{(\frac{x}{\nu})^2+(1+\frac{3}{2\nu})^2}\,\Big)^2} \Bigg) \frac{x^4}{32(\nu+1)^2(\nu+2)}
\label{Glow1claim}
\end{eqnarray}
and
\begin{eqnarray}
\hspace{-14mm} 
\lefteqn{
\hspace{-10mm} 
G^{\rm up,1}_\nu(x)
=
\frac{x^2}{4\nu^2(\nu+1)}
+
\frac{(4\nu^2+5\nu+2)x^4}{32\nu^3(\nu+1)^2(\nu+2)}
}
\nonumber
\\[2mm]
& & 
\hspace{13mm} 
+
\Bigg( 1 - \frac{4}{\Big(1+\sqrt{(\frac{x}{\nu+2})^2+1}\, \Big)^2} \Bigg) \frac{x^4}{32(\nu+1)^2(\nu+2)}
\cdot
\label{Gup1claim}
\end{eqnarray}
Since both~$c^{\rm low}$ and~$c^{\rm up}$ are larger than one, we easily obtain
\begin{equation}
	\label{Glow2claim}
|G^{\rm low,2}_\nu(x)|
\leq
{\textstyle{ \Big((1+\frac{3}{2\nu})+\sqrt{(\frac{x}{\nu})^2+(1+\frac{3}{2\nu})^2}\Big)^{-3} }}
 \frac{(\nu+\frac{1}{2})x^6}{24\nu^3(\nu+1)^3}
\end{equation}
and
\begin{equation}
	\label{Gup2claim}
\hspace{-3mm} 
|G^{\rm up,2}_\nu(x)|
\leq
{\textstyle{ \Big((1+\frac{2}{\nu})+\sqrt{(\frac{x}{\nu})^2+(1+\frac{2}{\nu})^2}\Big)^{-3} }}
 \frac{x^6}{24\nu^2(\nu+1)^3}
\cdot
\end{equation}
Using the mean value theorem to control the last term in the righthand sides of~(\ref{Glow1claim})--(\ref{Gup1claim}), it directly follows from~(\ref{Glow1claim})--(\ref{Gup2claim}) that, under the assumptions of the lemma, 
$$
a_n^2 G^{{\rm low/up},1}_{\nu_n}\Big(\frac{b_n T_n}{a_n}\Big)
=
o_{\rm P}(1)
\quad \textrm{and} \quad 
a_n^2 G^{{\rm low/up},2}_{\nu_n}\Big(\frac{b_n T_n}{a_n}\Big)
=
o_{\rm P}(1)
,
$$ 
which proves~(\ref{toshowforchristmas}), hence establishes the result.
\cqfd
\vspace{3mm}

%%%%%%%%%%%%%%%%%%%%

{\sc Proof of Theorem~\ref{TheoLAQ}.}
Throughout the proof, distributions and expectations are  under~${\rm P}^{(n)}_{\thetab_{n0},\kappa_n}$ and stochastic convergences are as~$n\to\infty$ under the same sequence of hypotheses. By using the fact that~$\Ob\thetab_{n0}=\thetab_{n0}=\Ob\pr\thetab_{n0}$ for any~$\Ob\in SO_{p_n}(\thetab_{n0})$ and by decomposing~$\taub_n$ into~$(\taub_n\pr \thetab_{n0}) \thetab_{n0}+ \Pi_{\thetab_{n0}} \taub_n$, with~$\Pi_{\thetab_{n0}}:=\mathbf{I}_{p_n}-\thetab_{n0}\thetab_{n0}\pr$, (\ref{tsj}) yields
\begin{eqnarray*}
\lefteqn{
 \frac{d{\rm P}^{(n)\Tb_n}_{1-\frac{1}{2} \nu_n^2 \|\taub_n\|^2,\kappa_n}}{d{m_n}} 
% \frac{d{\rm P}^{(n)\Tb_n}_{\thetab_{n},\kappa_n}}{d{m_n}} 
=
\frac{c_{p_n,\kappa_n}^n}{\omega_{p_n-1}^n}
\,
\int_{SO_{p_n}(\thetab_{n0})}
\exp\big( n \kappa_n \bar{\Xb}_{n}\pr \Ob\pr(\thetab_{n0}+\nu_n\taub_n) \big) 
\,d\Ob
}
\\[2mm]
& & \hspace{-2mm} 
=\,
\frac{c_{p_n,\kappa_n}^n}{\omega_{p_n-1}^n}
\exp\big( n \kappa_n \bar{\Xb}_{n}\pr \thetab_{n0} \big) 
\int_{SO_{p_n}(\thetab_{n0})}
\exp\big( n \kappa_n\nu_n \bar{\Xb}_{n}\pr \Ob\pr\taub_n \big) 
\,d\Ob
\\[2mm]
& & \hspace{-2mm} 
=\,
\frac{c_{p_n,\kappa_n}^n}{\omega_{p_n-1}^n}
\exp\big( n \kappa_n \bar{\Xb}_{n}\pr \thetab_{n0} \big) 
\int_{SO_{p_n}(\thetab_{n0})}
\exp\big( n \kappa_n\nu_n \bar{\Xb}_{n}\pr 
[(\taub_n\pr \thetab_{n0}) \thetab_{n0}+ \Ob\pr \Pi_{\thetab_{n0}} \taub_n] \big) 
\,d\Ob	
\\[2mm]
& & \hspace{-2mm} 
=\,
\frac{c_{p_n,\kappa_n}^n}{\omega_{p_n-1}^n}
\exp\big( n \kappa_n (1+\nu_n (\taub_n\pr \thetab_{n0})) \bar{\Xb}_{n}\pr \thetab_{n0} \big) 
\int_{SO_{p_n}(\thetab_{n0})}
\exp\big( n \kappa_n\nu_n \bar{\Xb}_{n}\pr \Ob\pr\Pi_{\thetab_{n0}} \taub_n \big) 
\,d\Ob	
.
\end{eqnarray*}
Now, since ${\bf O}\pr \Pi_{\thetab_{n0}}={\bf O}\pr \Pi_{\thetab_{n0}}^2=\Pi_{\thetab_{n0}} {\bf O}\pr \Pi_{\thetab_{n0}}$,
\begin{eqnarray*}
\lefteqn{
\int_{SO_{p_n}(\thetab_{n0})}
\exp\big( n \kappa_n\nu_n \bar{\Xb}_{n}\pr \Ob\pr\Pi_{\thetab_{n0}} \taub_n \big) 
\,d\Ob	
%=
%\int_{SO_{p_n}(\thetab_{n0})}
%\exp\big( n \kappa_n\nu_n \bar{\Xb}_{n}\pr \Pi_{\thetab_{n0}} \Ob\pr \Pi_{\thetab_{n0}} \taub_n \big) 
%\,d\Ob	
}
\\[2mm]
& &
\hspace{-2mm} 
=
\int_{SO_{p_n}(\thetab_{n0})}
\exp\big( n \kappa_n\nu_n \bar{\Xb}_{n}\pr \Pi_{\thetab_{n0}} \Ob\pr \Pi_{\thetab_{n0}} \taub_n \big) 
\,d\Ob	
\\[2mm]
& &
\hspace{-2mm} 
=\,
\int_{SO_{p_n}(\thetab_{n0})}
\exp\bigg( n \kappa_n\nu_n 
\|\Pi_{\thetab_{n0}}\taub_n\|
\|\Pi_{\thetab_{n0}}\bar{\Xb}_{n}\|
\bigg(\frac{\Pi_{\thetab_{n0}}\bar{\Xb}_{n}}{\|\Pi_{\thetab_{n0}}\bar{\Xb}_{n}\|}\bigg)\pr  
\bigg(
\Ob\pr \frac{\Pi_{\thetab_{n0}}\taub_n}{\|\Pi_{\thetab_{n0}}\taub_n\|}\bigg)
\bigg) 
\,d\Ob	
%
%\\[2mm]
%& &
%\hspace{1mm} 
%=\,
%{\rm E}
%\Big[
%\exp\bigg( n \kappa_n\nu_n 
%\|\Pi_{\thetab_{n0}}\taub_n\|
%\|\Pi_{\thetab_{n0}}\bar{\Xb}_{n}\|
%\bigg(\frac{\Pi_{\thetab_{n0}}\bar{\Xb}_{n}}{\|\Pi_{\thetab_{n0}}\bar{\Xb}_{n}\|}\bigg)\pr  
%\Ub
%\bigg)
%\,|\, 
%\Xb_{n1},\ldots,\Xb_{nn}
%\Big]
%%
\\[3mm]
& &
\hspace{-2mm} 
=\, 
{\rm E}
\big[
\exp\big( n \kappa_n\nu_n 
\|\Pi_{\thetab_{n0}}\taub_n\|
\|\Pi_{\thetab_{n0}}\bar{\Xb}_{n}\|
\vb_n\pr
\Sb
\big)
\,|\, 
\Xb_{n1},\ldots,\Xb_{nn}
\big]
,
\end{eqnarray*}
where~$\Sb$ is uniformly distributed over~$\mathcal{S}^{p_n-1}_{\thetab_{n0}}:=\{\xb\in\mathcal{S}^{p_n-1}: \xb'\thetab_{n0}=0\}$ and where~$\vb_n\in\mathcal{S}^{p_n-1}_{\thetab_{n0}}$ is {arbitrary}. Since~$\vb_n'\Sb$ has density~$t\mapsto c_{p_n-1}(1-t^2)^{\frac{p_n-4}{2}} \, \mathbb{I}[t\in[-1,1]]$, with~$c_{p_n-1}=1/\int_{-1}^1 (1-t^2)^{\frac{p_n-4}{2}}\,dt$, this yields 
\begin{eqnarray*}
\lefteqn{
\int_{SO_{p_n}(\thetab_{n0})}
\exp\big( n \kappa_n\nu_n \bar{\Xb}_{n}\pr \Ob\pr\Pi_{\thetab_{n0}} \taub_n \big) 
\,d\Ob	
}
\\[1mm]
& &
\hspace{10mm} 
= \,
c_{p_n-1}\int_{-1}^1 {\rm exp}\big(n \kappa_n \nu_n  \|\Pi_{\thetab_{n0}}\taub_n\|
\|\Pi_{\thetab_{n0}}\bar{\Xb}_{n}\|  t \big)(1-t^2)^{\frac{p_n-4}{2}}\, dt
\\[3mm] 
& & 
\hspace{10mm} 
= \,
H_{\frac{p_n-3}{2}}(n \kappa_n \nu_n  \|\Pi_{\thetab_{n0}}\taub_n\|
\|\Pi_{\thetab_{n0}}\bar{\Xb}_{n}\|)
.
\end{eqnarray*}
Summing up,
\begin{eqnarray}
\lefteqn{	
 \frac{d{\rm P}^{(n)\Tb_n}_{1-\frac{1}{2} \nu_n^2 \|\taub_n\|^2,\kappa_n}}{d{m_n}} 
% \frac{d{\rm P}^{(n)\Tb_n}_{\thetab_{n},\kappa_n}}{d{m_n}} 
= 
\frac{c_{p_n,\kappa_n}^n}{\omega_{p_n-1}^n}
\exp\big( n \kappa_n \bar{\Xb}_{n}\pr \thetab_{n0} \big)
}
\nonumber
\\[2mm]
& & 
\hspace{5mm} 
\times \,
\exp\big( n \kappa_n \nu_n (\taub_n\pr \thetab_{n0}) \bar{\Xb}_{n}\pr \thetab_{n0} \big) 
H_{\frac{p_n-3}{2}}(n \kappa_n \nu_n  \|\Pi_{\thetab_{n0}}\taub_n\|
\|\Pi_{\thetab_{n0}}\bar{\Xb}_{n}\|)
.
\label{invlikeproof}
\end{eqnarray}
Now, with the quantity~$L_n$ introduced in Lemma~\ref{LemSlutzkyFvML}, we have
\begin{eqnarray*}
\lefteqn{
T_n
:=
1 + \frac{\sqrt{2}\,\widetilde{W}_n}{\sqrt{p_n-1}} 
%=
%1 + \frac{\sqrt{2}}{\sqrt{p_n-1}} 
%\, 
%\bigg(\frac{W_n-(p_n-1)}{\sqrt{2(p_n-1)}}\bigg)
=
\frac{W_{n}}{p_n-1}
=
\frac{n^2}{\sum_{i=1}^n V_{ni}^2}
\,
\bar{\Xb}_{n}' (\mathbf{I}_{p_n}-\thetab_{n0}\thetab_{n0}') \bar{\Xb}_n
}
\\[2mm]
& & 
\hspace{33mm} 
%=
%\frac{n^2}{\sum_{i=1}^n V_{ni}^2}
%\,
% \|\Pi_{\thetab_{n0}}\bar{\Xb}_{n}\|^2
=
\frac{n}{f_{n2} L_n}
\,
 \|\Pi_{\thetab_{n0}}\bar{\Xb}_{n}\|^2
=
\frac{n\kappa_n}{(p_n-1)e_{n1} L_n}
\,
 \|\Pi_{\thetab_{n0}}\bar{\Xb}_{n}\|^2
,
\end{eqnarray*}
where we used the identity~$f_{n2}=(p_n-1)e_{n1}/\kappa_n$; see~(\ref{linken1en2tildekappan}). Therefore, (\ref{invlikeproof}) yields
\begin{eqnarray*}
\lefteqn{
\Lambda^{(n){\rm inv}}_{\thetab_n/\thetab_{n0};\kappa_n}
=
\log \frac{d{\rm P}^{(n)\Tb_n}_{1-\frac{1}{2} \nu_n^2 \|\taub_n\|^2,\kappa_n}}{d{\rm P}^{(n)\Tb_n}_{1,\kappa_n}}
}
\\[2mm]
& & 
\hspace{1mm} 
=\,
 n \kappa_n \nu_n (\taub_n\pr \thetab_{n0}) \bar{\Xb}_{n}\pr \thetab_{n0} 
+
\log
H_{\frac{p_n-3}{2}}(n \kappa_n \nu_n  \|\Pi_{\thetab_{n0}}\taub_n\|
\|\Pi_{\thetab_{n0}}\bar{\Xb}_{n}\|)
\\[2mm]
& & 
\hspace{1mm} 
=\,
 n \kappa_n \nu_n (\taub_n\pr \thetab_{n0}) \bar{\Xb}_{n}\pr \thetab_{n0} 
% -\frac{1}{2} n \kappa_n \nu_n^2 \|\taub_n\|^2 \bar{\Xb}_{n}\pr \thetab_{n0} 
%\\[2mm]
%& & 
%\hspace{21mm} 
+
\log
H_{\frac{p_n-3}{2}}(n^{1/2} (p_n-1)^{1/2} \kappa_n^{1/2} \nu_n e_{n1}^{1/2} 
%\|\taub_n\|^2 ( 1- {\textstyle{\frac{1}{4}}} \nu_n^2 \|\taub_n\|^2)
\|\Pi_{\thetab_{n0}}\taub_n\|
 L_n^{1/2} T_n^{1/2})
.
\end{eqnarray*}
Since 
$
\widetilde{W}_n
=
\sqrt{(p_n-1)/2}\times ( T_n - 1 )
$
is asymptotically standard normal (Lemma~\ref{LemJointNormality}), we have that
$
T_n = 1 + o_{\rm P}(1) 
$. Moreover, it directly follows from Lemma~\ref{LemSlutzkyFvML} that~$L_n = 1 + o_{\rm P}(1)$. Consequently, Lemma~\ref{LemGeneralDLFvML} shows that, if~$\nu_n$ satisfies~(\ref{taillenu}), then 
\begin{eqnarray*}
\Lambda^{(n){\rm inv}}_{\thetab_n/\thetab_{n0};\kappa_n}
%&=&
%n \kappa_n \nu_n (\taub_n\pr \thetab_{n0}) \bar{\Xb}_{n}\pr \thetab_{n0} 
%+
%\frac{1}{2(p_n-1)}
%n (p_n-1) \kappa_n \nu_n^2 e_{n1} \|\Pi_{\thetab_{n0}}\taub_n\|^2 T_n^2
%\\[2mm]
%& & 
%-
%\frac{1}{4(p_n-1)^2(p_n+1)}
%n^2 (p_n-1)^2 \kappa_n^2 \nu_n^4 e_{n1}^2 \|\Pi_{\thetab_{n0}}\taub_n\|^4 T_n^4
%\\[2mm]
&=&
n \kappa_n \nu_n (\taub_n\pr \thetab_{n0}) \bar{\Xb}_{n}\pr \thetab_{n0} 
+
\frac{p_n-1}{2(p_n-3)}
n \kappa_n \nu_n^2 e_{n1} \|\Pi_{\thetab_{n0}}\taub_n\|^2 L_n T_n
\\[2mm]
& & 
-
\frac{(p_n-1)^2}{4(p_n-3)^3}
n^2 \kappa_n^2 \nu_n^4 e_{n1}^2 \|\Pi_{\thetab_{n0}}\taub_n\|^4 L_n^{2}T_n^2
+ o_{\rm P}(1)
.
\end{eqnarray*}
Using~(\ref{taillenu}), Lemma~\ref{LemSlutzkyFvML} and the fact that~$T_n=1+o_{\rm P}(1)$ yields
\begin{eqnarray}
\Lambda^{(n){\rm inv}}_{\thetab_n/\thetab_{n0};\kappa_n}
&=&
n \kappa_n \nu_n (\taub_n\pr \thetab_{n0}) \bar{\Xb}_{n}\pr \thetab_{n0} 
+
\frac{1}{2}
n \kappa_n \nu_n^2 e_{n1} \|\Pi_{\thetab_{n0}}\taub_n\|^2 T_n
\label{fromTheo1}
\\[2mm]
& & 
\hspace{3mm} 
-
\frac{n^2 \kappa_n^2 \nu_n^4 e_{n1}^2}{4p_n}
 \|\Pi_{\thetab_{n0}}\taub_n\|^4 
+ o_{\rm P}(1)
.
\nonumber
\end{eqnarray}

Using the definitions of~$Z_n$ and~$T_n$, we obtain
\begin{eqnarray*}
\lefteqn{
\Lambda^{(n){\rm inv}}_{\thetab_n/\thetab_{n0};\kappa_n}
=
n \kappa_n \nu_n (\taub_n\pr \thetab_{n0}) 
\Big(e_{n1}+\frac{\sqrt{\tilde{e}_{n2}}}{n^{1/2}} Z_n \Big)
}
\\[2mm]
& & 
\hspace{-3mm} 
+
\frac{1}{2}
n \kappa_n \nu_n^2 e_{n1} \|\Pi_{\thetab_{n0}}\taub_n\|^2 
\Big( 1+ \frac{\sqrt{2}\,\widetilde{W}_n}{\sqrt{p_n-1}} \Big)
-
\frac{n^2 \kappa_n^2 \nu_n^4 e_{n1}^2}{4p_n}
 \|\Pi_{\thetab_{n0}}\taub_n\|^4 
+ o_{\rm P}(1)
\\[2mm]
& & 
\hspace{10mm} 
=
\sqrt{n} \kappa_n \nu_n \sqrt{\tilde{e}_{n2}} (\taub_n\pr \thetab_{n0})  Z_n 
+
\frac{n \kappa_n \nu_n^2 e_{n1}}{\sqrt{2}(p_n-1)^{1/2}}
 \|\Pi_{\thetab_{n0}}\taub_n\|^2 
 \widetilde{W}_n 
\\[2mm]
& & 
\hspace{-3mm} 
+
n \kappa_n \nu_n (\taub_n\pr \thetab_{n0}) e_{n1}
+
\frac{1}{2}
n \kappa_n \nu_n^2 e_{n1} \|\Pi_{\thetab_{n0}}\taub_n\|^2 
-
\frac{n^2 \kappa_n^2 \nu_n^4 e_{n1}^2}{4p_n}
 \|\Pi_{\thetab_{n0}}\taub_n\|^4 
+ o_{\rm P}(1)
.
\end{eqnarray*}
Using the identities
\begin{equation}
	\label{ide11}
\taub_n\pr \thetab_{n0}
=
- \frac{1}{2} \nu_n \|\taub_n\|^2 
\end{equation}
and
\begin{equation}
	\label{ide12}
\|\Pi_{\thetab_{n0}}\taub_n\|^2 
=
\|\taub_n\|^2 - (\taub_n\pr \thetab_{n0})^2
=
\|\taub_n\|^2 \Big(1- \frac{1}{4} \nu_n^2 \|\taub_n\|^2 \Big)
\end{equation}
provides
\begin{eqnarray*}
\lefteqn{
\Lambda^{(n){\rm inv}}_{\thetab_n/\thetab_{n0};\kappa_n}
=
-\frac{1}{2} \sqrt{n} \kappa_n \nu_n^2 \sqrt{\tilde{e}_{n2}} \|\taub_n\|^2  Z_n 
}
\\[2mm]
& & 
\hspace{-6mm} 
+
\frac{n \kappa_n \nu_n^2 e_{n1}}{\sqrt{2}(p_n-1)^{1/2}}
\|\taub_n\|^2 \Big(1- \frac{1}{4} \nu_n^2 \|\taub_n\|^2 \Big)
  \widetilde{W}_n 
-\frac{1}{2} n \kappa_n \nu_n^2 e_{n1} \|\taub_n\|^2 
\\[2mm]
& & 
\hspace{-6mm} 
+
\frac{1}{2}
n \kappa_n \nu_n^2 e_{n1} 
\|\taub_n\|^2 \Big(1- \frac{1}{4} \nu_n^2 \|\taub_n\|^2 \Big)
-
\frac{n^2 \kappa_n^2 \nu_n^4 e_{n1}^2}{4p_n}
\|\taub_n\|^4  \Big(1- \frac{1}{4} \nu_n^2 \|\taub_n\|^2 \Big)^2
+ o_{\rm P}(1)
.
\end{eqnarray*}
The result then follows from~(\ref{taillenu}) and from the tightness of~$\widetilde{W}_n$ (Lemma~\ref{LemJointNormality}).
\cqfd

%%%%%%%%%%%%%%%%%%%%%%%%%%%%%%%%%%%%%%%%%%%%%%%%%%%%%%

\begin{lemma}
\label{LemEnFvML} 
Let~$(p_n)$ be a sequence of integers that diverges to infinity and~$(\kappa_n)$ be an arbitrary sequence in~$(0,\infty)$. Let~$e_{n1}$ $($resp.,~$\tilde{e}_{n2})$ be the expectation $($resp., the variance$)$ of the distribution with probability density function~(\ref{pdfUFvML}). Then, we have the following: 
(i) 
if~$\kappa_n/p_n\to\infty$, then 
$$
e_{n1}=1+o(1)
,
\quad
\tilde{e}_{n2}=O\Big(\frac{p_n}{\kappa_n^2}\Big)
\quad\textrm{ and }\quad
f_{n2}=\frac{p_n}{\kappa_n}+o\Big(\frac{p_n}{\kappa_n}\Big)
;
$$
(ii) if~$\kappa_n/p_n\to\xi>0$, then, letting~$c_{\xi}:=\frac{1}{2}+\sqrt{\frac{1}{4} + \xi^2}$,  
$$
e_{n1}\to \frac{\xi}{c_\xi}+o(1),
\quad
\tilde{e}_{n2}=O\Big(\frac{1}{p_n}\Big)
\quad\textrm{ and }\quad
f_{n2}=\frac{1}{c_\xi}+o(1)
;
$$
(iii) 
if~$\kappa_n/p_n\to 0$, then 
$$
e_{n1}=\frac{\kappa_n}{p_n}+O\Big(\frac{\kappa_n^3}{p_n^3}\Big)
,
\quad
\tilde{e}_{n2}=\frac{1}{p_n}+o\Big(\frac{1}{p_n}\Big) 
\quad\textrm{ and }\quad
f_{n2}=1+o(1)
.
$$
\end{lemma}
\vspace{1mm}
 
{\sc Proof of Lemma~\ref{LemEnFvML}.} Denoting again as~$\mathcal{I}_\nu(\cdot)$  the order-$\nu$ modified Bessel function of the first kind, we recall (see~(\ref{linken1en2tildekappan})) that 
$$ 
e_{n1}
=
\frac{\mathcal{I}_{\frac{p_n}{2}}(\kappa_n)}{\mathcal{I}_{\frac{p_n}{2}-1}(\kappa_n)}
,
\quad
\tilde{e}_{n2}
=
1-\frac{p_n-1}{\kappa_n}\, e_{n1}-e_{n1}^2
\quad\textrm{ and }\quad
f_{n2}
=
\frac{p_n-1}{\kappa_n}\, e_{n1}
.
$$
%, e.g., Lemma~S.2.1 in \cite{PaiVer17asupp}. 
In each case~(i)--(iii), the claim for~$f_{n2}$ directly follows from the result on~$e_{n1}$, so that it is sufficient to prove the results for~$e_{n1}$ and~$\tilde{e}_{n2}$.  To do so, we will use the bounds
\begin{equation}
\label{Amos}
R^{\rm low}_{\nu}(z)
:=
\frac{z}{\nu+1+\sqrt{(\nu+1)^2+z^2}}
%G_{\nu+1,\nu+1}(z)
\leq
\frac{\mathcal{I}_{\nu+1}(z)}{\mathcal{I}_\nu(z)}
\leq
\frac{z}{\nu+\sqrt{(\nu+2)^2+z^2}}
=:
R^{\rm up}_{\nu}(z)
%G_{\nu,\nu+2}(z)
\end{equation}
and
\begin{equation}
\label{Amostilde}
\tilde{R}^{\rm low}_{\nu}(z)
:=
\frac{z}{\nu+\frac{1}{2}+\sqrt{(\nu+\frac{3}{2})^2+z^2}}
\leq
\frac{\mathcal{I}_{\nu+1}(z)}{\mathcal{I}_\nu(z)}
%\leq
%\frac{z}{\nu+\frac{1}{2}+\sqrt{(\nu+\frac{1}{2})(\nu+\frac{3}{2})+z^2}}
%=:
%\tilde{R}^{\rm up}_{\nu}(z)
;
\end{equation}
see (11) and~(16) in \cite{Amo1974}, respectively.
%where we let~$G_{\alpha,\beta}(t)=z/(\alpha+\sqrt{\beta^2+z^2})$; 
%(the bounds in~(\ref{Amos}) and the lower bound in~(\ref{Amostilde}) were obtained in \cite{Amo1974}, whereas the upper bound in~(\ref{Amostilde}) was derived in \cite{SimSpe1984}). 
%
(i) The lower bound in~(\ref{Amos}) provides
$$
e_{n1}
\geq
\frac{\kappa_n}{\frac{p_n}{2}+\sqrt{\big(\frac{p_n}{2}\big)^2 + \kappa_n^2}}
=
\frac{1}{\frac{p_n}{2\kappa_n}+\sqrt{\big(\frac{p_n}{2\kappa_n}\big)^2 + 1}}
,
$$
which, since~$e_{n1}\leq 1$, establishes the result for~$e_{n1}$. 
%It therefore only remains to prove the result for~$\tilde{e}_{n2}$ in case~(iii).
 Making use of the bound in~(\ref{Amostilde}), we can write
 $$
%1 
%- \frac{p-1}{\kappa} \tilde{R}^{\rm up}_{\frac{p_n}{2}-1}(\kappa_n) 
%- \big(\tilde{R}^{\rm up}_{\frac{p_n}{2}-1}(\kappa_n)\big)^2
%\leq 
\tilde{e}_{n2}
\leq
1 
- \frac{p_n-1}{\kappa_n} \tilde{R}^{\rm low}_{\frac{p_n}{2}-1}(\kappa_n) 
- \big(\tilde{R}^{\rm low}_{\frac{p_n}{2}-1}(\kappa_n)\big)^2
.
$$
Lengthy yet quite straightforward computations allow to rewrite this as
$$
%\frac{
%p_n\big(1-\frac{1}{p_n}\big)
%}{2\kappa_n^2\Big(\frac{p_n-1}{2\kappa_n}+\sqrt{\big(\frac{p_n^2-1}{4\kappa_n^2}\big) + 1}\,\Big)^2}
%\leq 
\tilde{e}_{n2}
\leq
\frac{
p_n
 }{\kappa_n^2 \Big(\frac{p_n-1}{2\kappa_n}+\sqrt{\big(\frac{p_n+1}{2\kappa_n}\big)^2 + 1}\,\Big)^2}
\cdot
$$
It readily follows that~$\kappa_n^2 \tilde{e}_{n2}/p_n$ is~$O(1)$, as was to be showed. 
Let us turn to the proof of~(iii). The bounds in~(\ref{Amos}) readily yield
\begin{equation}
	\label{keymajorE1nUn}
\frac{1}{\frac{1}{2}+\sqrt{\big(\frac{1}{2}\big)^2 + \big(\frac{\kappa_n}{p_n}\big)^2}}
\leq 
\frac{e_{n1}}{\kappa_n/p_n}
\leq 
\frac{1}{\frac{1}{2}-\frac{1}{p_n}+\sqrt{\big(\frac{1}{2}+\frac{1}{p_n}\big)^2 + \big(\frac{\kappa_n}{p_n}\big)^2}}
,
\end{equation}
which provides 
%\newpage
%$$
%\frac{1}{\frac{1}{2}+\sqrt{\big(\frac{1}{2}\big)^2 + \big(\frac{\kappa_n}{p_n}\big)^2}} - 1
%\leq 
%\frac{e_{n1}}{\kappa_n/p_n} - 1
%\leq 
%\frac{1}{\frac{1}{2}-\frac{1}{p_n}+\sqrt{\big(\frac{1}{2}+\frac{1}{p_n}\big)^2 + \big(\frac{\kappa_n}{p_n}\big)^2}} -1
%$$
%%
%$$
%\frac{\frac{1}{2}-\sqrt{\big(\frac{1}{2}\big)^2 + \big(\frac{\kappa_n}{p_n}\big)^2}}{\frac{1}{2}+\sqrt{\big(\frac{1}{2}\big)^2 + \big(\frac{\kappa_n}{p_n}\big)^2}}
%\leq 
%\frac{e_{n1}}{\kappa_n/p_n} - 1
%\leq 
%\frac{\frac{1}{2}+\frac{1}{p_n}-\sqrt{\big(\frac{1}{2}+\frac{1}{p_n}\big)^2 + \big(\frac{\kappa_n}{p_n}\big)^2}}{\frac{1}{2}-\frac{1}{p_n}+\sqrt{\big(\frac{1}{2}+\frac{1}{p_n}\big)^2 + \big(\frac{\kappa_n}{p_n}\big)^2}} 
%$$
%%
$$
\frac{-\big(\frac{\kappa_n}{p_n}\big)^2}{\big( \frac{1}{2}+\sqrt{\big(\frac{1}{2}\big)^2 + \big(\frac{\kappa_n}{p_n}\big)^2}\big)^2}
\leq 
\frac{e_{n1}}{\kappa_n/p_n} - 1
\leq 
\frac{-\big(\frac{\kappa_n}{p_n}\big)^2}{\frac{1}{2}-\frac{1}{p_n}+\sqrt{\big(\frac{1}{2}+\frac{1}{p_n}\big)^2 + \big(\frac{\kappa_n}{p_n}\big)^2}} 
\cdot
$$
This proves the result for~$e_{n1}$. Turning to~$\tilde{e}_{n2}$, the bounds in~(\ref{Amos}) lead to
$$
1 
- \frac{p_n-1}{\kappa_n} R^{\rm up}_{\frac{p_n}{2}-1}(\kappa_n) 
- \big(R^{\rm up}_{\frac{p_n}{2}-1}(\kappa_n)\big)^2
\leq 
\tilde{e}_{n2}
\leq
1 
- \frac{p_n-1}{\kappa_n} R^{\rm low}_{\frac{p_n}{2}-1}(\kappa_n) 
- \big(R^{\rm low}_{\frac{p_n}{2}-1}(\kappa_n)\big)^2
.
$$
As above, heavy but rather straightforward computations allow to rewrite this as
\begin{equation}
	\label{keymajorE2nUn}
	\frac{
\frac{3}{2}    
+
\frac{1}{p_n}   
-
 \sqrt{\big(\frac{1}{2}+\frac{1}{p_n}\big)^2 + \big(\frac{\kappa_n}{p_n}\big)^2}
}
{p_n\Big(\frac{1}{2}-\frac{1}{p_n}+ \sqrt{\big(\frac{1}{2}+\frac{1}{p_n}\big)^2 + \big(\frac{\kappa_n}{p_n}\big)^2}\,\Big)^2} 
\leq
\tilde{e}_{n2}
\leq
 \frac{
1
 }{
 p_n\Big(\frac{1}{2}+\sqrt{\big(\frac{1}{2}\big)^2 + \big(\frac{\kappa_n}{p_n}\big)^2}\,\Big)
 }
 ,
\end{equation}
which establishes that~$p_n\tilde{e}_{n2}=1+o(1)$. 
Finally, the result in~(ii) readily follows from~(\ref{keymajorE1nUn}) and from the upper bound in~(\ref{keymajorE2nUn}). 
\cqfd
\vspace{3mm}
 
%%%%%%%%%%%%%%%%%%%%%%%%%%%%%%%%%%%%%%%%%%%%%%%%%%%%%%%%%%%%%%%%%

{\sc Proof of Theorem~\ref{TheoLAN}.}
Stochastic convergences throughout the proof are as~$n\to\infty$ under~${\rm P}^{(n)}_{\thetab_{n0},\kappa_n}$. 
  %We consider several cases. 
Assume that (i)~$\kappa_n/p_n\to \infty$, 
(ii)~$\kappa_n/p_n\to \xi>0$, or
(iii)~$\kappa_n/p_n\to 0$ with~$\sqrt{n} \kappa_n/p_n\to \infty$, and let~$(\nu_n)$ be the corresponding sequence in the statement of the theorem.  
Using Lemma~\ref{LemEnFvML} and the identity~$\kappa_n f_{n2}=(p_n-1)e_{n1}$, it is then easy to check that~$\nu_n$ satisfies~(\ref{taillenu}), is such that~$\nu_n=o(1)$, and is asymptotically equivalent to $\tilde{\nu}_n=p_n^{3/4}/(\sqrt{n} \kappa_n \sqrt{f_{n2}})$ in the sense that~$\tilde{\nu}_n/\nu_n\to 1$. Theorem~\ref{TheoLAQ} thus applies and yields
	\begin{eqnarray*}
\lefteqn{	
\hspace{-3mm} 
\Lambda^{(n){\rm inv}}_{\thetab_n/\thetab_{n0};\kappa_n}
=
-
\frac{p_n^{3/2}\sqrt{\tilde{e}_{n2}}}{2\sqrt{n} \kappa_n f_{n2}} 
\, \|\taub_n\|^2  Z_n 
+ 
\frac{p_n e_{n1}}{\sqrt{2}\kappa_n f_{n2}} 
\|\taub_n\|^2 
  \widetilde{W}_n
}
\\[2mm]
& &  
\hspace{23mm} 
-
\frac{p_n^3 e_{n1}}{8n \kappa_n^3 f_{n2}^2} 
\|\taub_n\|^4
-
\frac{p_n^2e_{n1}^2}{4\kappa_n^2 f_{n2}^2} 
  \|\taub_n\|^4 
  +o_{\rm P}(1)
.
\end{eqnarray*}
Using again the identity~$\kappa_n f_{n2}=(p_n-1)e_{n1}$, we then obtain
	\begin{eqnarray}
\lefteqn{	
\hspace{7mm} 
\Lambda^{(n){\rm inv}}_{\thetab_n/\thetab_{n0};\kappa_n}
=
-
\frac{p_n^{3/2}\sqrt{\tilde{e}_{n2}}}{2\sqrt{n} (p_n-1) e_{n1}} 
\, \|\taub_n\|^2  Z_n 
+ 
\frac{1}{\sqrt{2}} 
\|\taub_n\|^2 
  \widetilde{W}_n
}
\label{njkls}
\\[2mm]
& &  
\hspace{33mm} 
-
\frac{p_n^3}{8n \kappa_n (p_n-1)^2 e_{n1}} 
\|\taub_n\|^4
-
\frac{1}{4} 
  \|\taub_n\|^4 
  +o_{\rm P}(1)
.
\nonumber
\end{eqnarray}
The result in cases~(i)--(iii) then follows from the fact that Lemma~\ref{LemEnFvML} implies that, in each case, the first and third term of the righthand side of~(\ref{njkls}) are~$o_{\rm P}(1)$. 
%we have that
%$$
%n^{1/2} \kappa_n \nu_n^2 \tilde{e}_{n2}^{1/2}
%=
%n^{1/2} \kappa_n \frac{p_n^{3/2}}{n \kappa_n^2 f_{n2}} \tilde{e}_{n2}^{1/2}
%=
%\frac{p_n}{\sqrt{n}\kappa_n f_{n2}} (p_n\tilde{e}_{n2})^{1/2}
%=
%\Bigg\{
%\begin{array}{ll}
%\frac{p_n}{\sqrt{n}(p_n-1) e_{n1}} \frac{p_n}{\kappa_n} 
%\Big(\frac{\kappa_n^2}{p_n}\tilde{e}_{n2}\Big)^{1/2}
%& \textrm{for (i)}
%\\[2mm]
%\frac{p_n}{\sqrt{n}(p_n-1) e_{n1}} (p_n\tilde{e}_{n2})^{1/2}
%& \textrm{for (ii)--(iii)}
%\end{array}
%$$
%and that
%$$
%n \kappa_n\nu_n^4 e_{n1}
%=
%n \kappa_n \frac{p_n^{3}}{n^2 \kappa_n^4 f_{n2}^2} e_{n1}
%=
%\frac{p_n^{3}}{n \kappa_n^3 f_{n2}^2} e_{n1}
%=
%\frac{p_n}{n \kappa_n}\frac{p_n^{2}}{\kappa_n^2 f_{n2}^2} e_{n1}
%=
%\frac{p_n}{n \kappa_n}\frac{p_n^{2}}{(p_n-1)^2 e_{n1}^2} e_{n1}
%=
%o(1)
%$$
%for~(i)--(iii). Therefore, using the fact that~$\nu_n=o(1)$ in all three cases, LAQ becomes
% \begin{eqnarray*}
% 	\lefteqn{
% \log 
%\Bigg( 
%\frac{d{\rm P}^{(n)\Tb_n}_{\thetab_n,\kappa_n}}{d{\rm P}^{(n)\Tb_n}_{\thetab_{n0},\kappa_n}}
%\Bigg)
%=
%\frac{n \kappa_n\nu_n^2 e_{n1}}{\sqrt{2}p_n^{1/2}}
%\|\taub_n\|^2   \widetilde{W}_n
%-
%\frac{n^2 \kappa_n^2 \nu_n^4 e_{n1}^2}{4p_n}  \|\taub_n\|^4  
%  +o_{\rm P}(1)
%}
%\\[2mm]
%& & 
%\hspace{3mm} 
%=
%\frac{p_n  e_{n1}}{\sqrt{2}\kappa_n f_{n2}}
%\|\taub_n\|^2   \widetilde{W}_n
%-
%\frac{p_n^2 e_{n1}^2}{4 f_{n2}^2\kappa_n^2}  \|\taub_n\|^4  
%  +o_{\rm P}(1)
%=
%\frac{1}{\sqrt{2}}
%\|\taub_n\|^2   \widetilde{W}_n
%-
%\frac{1}{4}  \|\taub_n\|^4  
%+o_{\rm P}(1)
%.
%\end{eqnarray*}

We turn to case~(iv), for which~$\sqrt{n}\kappa_n/p_n=\xi$ (so that, like for all subsequent cases,~$\kappa_n=o(p_n)$). Then, the same argument as above allows to check that $
\nu_n
=
p_n^{3/4}/(\sqrt{n} \kappa_n \sqrt{f_{n2}})$ still satisfies~(\ref{taillenu}) and is such that~$\nu_n=o(1)$, so that, jointly with Lemma~\ref{LemEnFvML}, Theorem~\ref{TheoLAQ} provides
	\begin{eqnarray*}
\Lambda^{(n){\rm inv}}_{\thetab_n/\thetab_{n0};\kappa_n}
&=&
-
\frac{p_n^{3/2} \tilde{e}_{n2}^{1/2}}{2\sqrt{n} \kappa_n f_{n2}}
 \|\taub_n\|^2  Z_n 
+ 
\frac{p_n e_{n1}}{\sqrt{2}\kappa_n f_{n2}}
\|\taub_n\|^2  \widetilde{W}_n
\\[2mm]
& &  
\hspace{23mm} 
-
\frac{p_n^{3} e_{n1}}{8n \kappa_n^3 f_{n2}^2}
 \|\taub_n\|^4
-
\frac{p_n^2 e_{n1}^2}{4\kappa_n^2 f_{n2}^2}  \|\taub_n\|^4  +o_{\rm P}(1)
\\[2mm]
&=&
-
\frac{1}{2\xi} 
\|\taub_n\|^2  Z_n 
+ 
\frac{1}{\sqrt{2}} \|\taub_n\|^2  \widetilde{W}_n
-
\frac{1}{8\xi^2} \|\taub_n\|^4
-
\frac{1}{4}  \|\taub_n\|^4  +o_{\rm P}(1)
,
\end{eqnarray*}
as was to be shown. Consider now case~(v), under which~$\sqrt{n} \kappa_n/p_n\to 0$ with~$\sqrt{n} \kappa_n/\sqrt{p_n}\to \infty$, which still ensures that~$\nu_n=p_n^{1/4}/(n^{1/4} \sqrt{\kappa_n})$ is~$o(1)$ and satisfies~(\ref{taillenu}). Theorem~\ref{TheoLAQ} applies and, by using Lemma~\ref{LemEnFvML} again, yields 
%$$
%\nu_n
%=
%\frac{p_n^{1/4}}{n^{1/4} \sqrt{\kappa_n}}
%=
%\Big( \frac{\sqrt{p_n}}{\sqrt{n}\kappa_n}\Big)^{1/2}
%(=o(1))
%,
%$$
%which satisfies
%$$
%\frac{n^{1/2}\kappa_n}{p_n^{3/4}} \nu_n
%=
%\frac{n^{1/2}\kappa_n}{p_n^{3/4}} 
%\frac{p_n^{1/4}}{n^{1/4}\sqrt{\kappa_n}} 
%=
%\frac{n^{1/4}\kappa_n^{1/2}}{p_n^{1/2}} 
%=
%o(1)
%,
%$$
	\begin{eqnarray*}
\hspace{-1mm} 
\Lambda^{(n){\rm inv}}_{\thetab_n/\thetab_{n0};\kappa_n}
&=&
-
\frac{1}{2}
\sqrt{p_n\tilde{e}_{n2}} \,
\|\taub_n\|^2  Z_n 
+ 
\frac{\sqrt{n} e_{n1}}{\sqrt{2}}
\|\taub_n\|^2  \widetilde{W}_n
%\\[2mm]
%& &
%\hspace{23mm} 
-
\frac{p_n e_{n1}}{8\kappa_n} \|\taub_n\|^4
-
\frac{n e_{n1}^2}{4}
 \|\taub_n\|^4  +o_{\rm P}(1)
\\[2mm]
&=&
-
\frac{1}{2} \|\taub_n\|^2  Z_n 
-
\frac{1}{8} \|\taub_n\|^4
+o_{\rm P}(1)
,
\end{eqnarray*}
which establishes the result in case~(v). If~$\sqrt{n}\kappa_n/\sqrt{p_n}=\xi$ (case~(vi)), then Lemma~\ref{LemEnFvML} implies that~$\nu_n=1$ satisfies~(\ref{taillenu}). Theorem~\ref{TheoLAQ} then provides
	\begin{eqnarray}
\Lambda^{(n){\rm inv}}_{\thetab_n/\thetab_{n0};\kappa_n}
&=&
-
\frac{\sqrt{n} \kappa_n}{2\sqrt{p_n}} \sqrt{p_n\tilde{e}_{n2}}\, \|\taub_n\|^2  Z_n 
+ 
\frac{n \kappa_n e_{n1}}{\sqrt{2}p_n^{1/2}}
\|\taub_n\|^2 \Big(1- \frac{1}{4} \|\taub_n\|^2\Big)
  \widetilde{W}_n
\nonumber
\\[2mm]
& &  
\hspace{1mm} 
-
\frac{1}{8}
n \kappa_n e_{n1} \|\taub_n\|^4
-
\frac{n^2 \kappa_n^2 e_{n1}^2}{4p_n}  \|\taub_n\|^4 \Big(1- \frac{1}{4} \|\taub_n\| ^2\Big)^2 
  +o_{\rm P}(1)
\label{toreferinproofofLAN}
\\[2mm]
&= &  
-
\frac{\xi}{2} \|\taub_n\|^2  Z_n 
- 
\frac{\xi^2}{8} \|\taub_n\|^4 +o_{\rm P}(1)
,
\nonumber
\end{eqnarray} 
where we used Lemma~\ref{LemEnFvML}. Finally, if~$\sqrt{n}\kappa_n/\sqrt{p_n}=o(1)$ (case~(vii)), then~(\ref{taillenu}) again holds for~$\nu_n=1$. Therefore, Theorem~\ref{TheoLAQ} shows that
$\Lambda^{(n){\rm inv}}_{\thetab_n/\thetab_{n0};\kappa_n}$ satisfies the first equality of~(\ref{toreferinproofofLAN}), hence is~$o_{\rm P}(1)$.  
\cqfd
\vspace{3mm}

\noindent {\sc Proof of Theorem~\ref{LANunspec}.} %Unless otherwise mentioned, stochastic convergences in this proof are as~$n\to\infty$ under~${\rm P}^{(n)}_{\thetab_{n0},\kappa_n}$, with~$\kappa_n=p_n\xi/\sqrt{n}$. 
First note that, since~$p_n=o(n^2)$, Lemma~\ref{LemEnFvML}(iii) entails that
\begin{eqnarray}
Z_n
&=&
\frac{\sqrt{n}(\bar{\Xb}_n\pr \thetab_{n0}- e_{n1})}{\sqrt{\tilde{e}_{2n}}}
=
\frac{\sqrt{n}(\bar{\Xb}_n\pr \thetab_{n0}- \xi/\sqrt{n}+O(n^{-3/2}))}{\sqrt{\frac{1}{p_n}+o(\frac{1}{p_n})}}
\nonumber
\\[2mm]
&= & 
%\hspace{6mm} 
%=
%\frac{\sqrt{n}\bar{\Xb}_n\pr \thetab_{n0}- \xi+O(1/n)}{\sqrt{\frac{1}{p_n}+o(\frac{1}{p_n})}}
\frac{\sqrt{p_n}(\sqrt{n}\bar{\Xb}_n\pr \thetab_{n0}- \xi+O(1/n))}{\sqrt{1+o(1)}}
%$$
%$$
%=
%\frac{\sqrt{np_n}\bar{\Xb}_n\pr \thetab_{n0}- \sqrt{p_n} \xi+O(\sqrt{p_n}/n)}{\sqrt{1+o(1)}}
= 
\sqrt{np_n}\bar{\Xb}_n\pr \thetab_{n0}- \sqrt{p_n} \xi
+
o_{\rm P}(1)
\label{repZn}	
\end{eqnarray}
as~$n\to\infty$ under~${\rm P}^{(n)}_{\thetab_{n0},\kappa_n}$. Write then
\begin{eqnarray*}
\Lambda^{(n){\rm inv}}_{\thetab_n, \kappa_{n,s}/\thetab_{n0},\kappa_n} 
&=& 
\Lambda^{(n){\rm inv}}_{\thetab_n, \kappa_{n,s}/\thetab_{n0},\kappa_{n,s}}
+ 
\Lambda^{(n){\rm inv}}_{\thetab_{n0}, \kappa_{n,s}/\thetab_{n0},\kappa_n} 
\\[1mm]
&=& 
\Lambda^{(n){\rm inv}}_{\thetab_n, \kappa_{n,s}/\thetab_{n0},\kappa_{n,s}}
+ 
\log \frac{d{\rm P}\n_{\thetab_{n0},\kappa_{n,s}}}{d{\rm P}\n_{\thetab_{n0},\kappa_{n}}}
\\[2mm]
&\!\!=:\!\!& 
L_{n1}+L_{n2}.
\end{eqnarray*}
Using~(\ref{repZn}), we obtain
\begin{eqnarray}
\label{LLn1}
L_{n2} 
& = &
 n \big( \log(c_{p_n,\kappa_{n,s}})-\log(c_{p_n, \kappa_n}) \big)
 + 
 n(\kappa_{n,s}-\kappa_n) \bar{\Xb}_n\pr \thetab_{n0} 
 \nonumber 
 \\[2.5mm]
& = & 
n \Big[
\log\Big(\frac{c_{p_n,0}}{c_{p_n, \kappa_n}} \Big)
-
\log\Big(\frac{c_{p_n,0}}{c_{p_n, \kappa_{n,s}}} \Big)
\Big]
 + 
 s \sqrt{n p_n}\, \bar{\Xb}_n\pr \thetab_{n0} 
\nonumber 
\\
& = & 
n 
\Big[
\log H_{\frac{p_n}{2}-1} (\kappa_n) 
-
\log H_{\frac{p_n}{2}-1} (\kappa_{n,s})
\Big]
+ 
s \sqrt{p_n} \xi
+
s Z_{n}
+
o_{\rm P}(1)
\nonumber 
\\[2mm]
& =: & 
\tilde{L}_{n2}
+
s Z_{n}
+
o_{\rm P}(1)
\nonumber 
\end{eqnarray}
as~$n\to\infty$ under~${\rm P}^{(n)}_{\thetab_{n0},\kappa_n}$. Since $p_n=o(n^2)$, we can apply Lemma \ref{LemGeneralDLFvML} with~$a_n=\sqrt{n}$ and $T_n \equiv 1$. This yields
%$$b_n/\sqrt{n}=\kappa_n(\xi)=p_n \xi/ \sqrt{n}$$
%meaning that $b_n=p_n$, since $\nu_n=p_n$, the condition is $b_n^6=o(a_n^4p_n^5)$ which is equivalent to $p_n=o(n^2)$
\begin{eqnarray}
\label{LLn2}
\tilde{L}_{n2} 
%&=& 
%\Big(
%\frac{n \kappa_n^2}{4({\textstyle{\frac{p_n}{2}}}-1)}
%-
%\frac{n \kappa_n^4}{32({\textstyle{\frac{p_n}{2}}}-1)^3}
%+
%\xi^2
%\Big)
%-
%\Big(
%\frac{n \kappa_{n,s}^2}{4({\textstyle{\frac{p_n}{2}}}-1)}
%-
%\frac{n \kappa_{n,s}^4}{32({\textstyle{\frac{p_n}{2}}}-1)^3}
%+
%\xi^2
%\Big)
%+
%s\sqrt{p_n} \xi
%+
%o(1)
%\nonumber  
%\\[2mm]
&=& 
\Big(
\frac{n \kappa_n^2}{4({\textstyle{\frac{p_n}{2}}}-1)}
-
\frac{n \kappa_n^4}{32({\textstyle{\frac{p_n}{2}}}-1)^3}
\Big)
-
\Big(
\frac{n \kappa_{n,s}^2}{4({\textstyle{\frac{p_n}{2}}}-1)}
-
\frac{n \kappa_{n,s}^4}{32({\textstyle{\frac{p_n}{2}}}-1)^3}
\Big)
+
s\sqrt{p_n} \xi
+
o(1)
\nonumber  
\\[2mm]
&=& 
-\frac{n(\kappa_{n,s}^2-\kappa_n^2)}{2p_n-4}
+ 
\frac{n(\kappa_{n,s}^4-\kappa_n^4)}{4(p_n-2)^3}
+
s\sqrt{p_n} \xi
+
o(1)
\nonumber  
\\[2mm]
&=& 
-\frac{p_n^2((\xi+ s/\sqrt{p_n})^2-\xi^2)}{2p_n-4}+ \frac{p_n^4((\xi+ s/\sqrt{p_n})^4-\xi^4)}{4n(p_n-2)^3}+s\sqrt{p_n} \xi+o(1) 
\nonumber 
\\
&=& 
-\frac{s^2}{2}+o(1) 
\nonumber 
\end{eqnarray}
as $\ny$. Therefore, 
\begin{equation} 
\label{L2}
L_{n2}
=
 s Z_{n}
- 
\frac{s^2}{2}+o_{\rm P}(1) 
\end{equation}
as~$n\to\infty$ under~${\rm P}^{(n)}_{\thetab_{n0},\kappa_n}$, which 
\vspace{-1mm}
 implies that the sequences of probability measures ${\rm P}^{(n)}_{\thetab_{n0},\kappa_{n,s}}$ and ${\rm P}^{(n)}_{\thetab_{n0},\kappa_n}$ are mutually contiguous (this results from the Le Cam first lemma). 

Now, denote as~$e_{n1,s}$ and~$\tilde{e}_{2n,s}$, respectively, the values of~$e_{n1}$ and~$\tilde{e}_{2n}$ under~${\rm P}^{(n)}_{\thetab_{n0},\kappa_{n,s}}$. Then, proceeding as in~(\ref{repZn}) and using the fact that contiguity  implies that~(\ref{repZn}) also holds under~${\rm P}^{(n)}_{\thetab_{n0},\kappa_{n,s}}$, one obtains
$$
Z_{n,s}
:=
\frac{\sqrt{n}(\bar{\Xb}_n\pr \thetab_{n0}- e_{n1,s})}{\sqrt{\tilde{e}_{2n,s}}}
%=
%\frac{\sqrt{n}(\bar{\Xb}_n\pr \thetab_{n0}- (\xi+s/\sqrt{p_n})/\sqrt{n}+O(n^{-3/2}))}{\sqrt{\frac{1}{p_n}+o(\frac{1}{p_n})}}
%$$
%$$
%=
%\frac{\sqrt{n}\bar{\Xb}_n\pr \thetab_{n0}- \xi -s/\sqrt{p_n}+O(1/n)}{\sqrt{\frac{1}{p_n}+o(\frac{1}{p_n})}}
%=
%\frac{\sqrt{p_n}(\sqrt{n}\bar{\Xb}_n\pr \thetab_{n0}- \xi-s/\sqrt{p_n}+O(1/n))}{\sqrt{1+o(1)}}
%$$
%$$
%=
%\frac{\sqrt{np_n}\bar{\Xb}_n\pr \thetab_{n0}- \sqrt{p_n} \xi-s+O(\sqrt{p_n}/n)}{\sqrt{1+o(1)}}
= 
\sqrt{np_n}\bar{\Xb}_n\pr \thetab_{n0}- \sqrt{p_n} \xi -s
+
o_{\rm P}(1)
=
Z_n
-
s
+
o_{\rm P}(1)
$$
as~$n\to\infty$ under~${\rm P}^{(n)}_{\thetab_{n0},\kappa_{n,s}}$. Consequently, Theorem~\ref{TheoLAN}(iv) implies that
\begin{eqnarray} 
L_{1n} 
&\!\! =\!\! & 
\| \taub_n\|^2 
\bigg(
\frac{\widetilde{W}_n}{\sqrt{2}}
-
\frac{Z_{n,s}}{2\xi}
\bigg)
-\frac{1}{2}
\| \taub_n\|^4 
\bigg(\frac{1}{2}
+
\frac{1}{4\xi^2}\bigg)
+o_{\rm P}(1) 
\nonumber 
\\[3mm]
&\!\! =\!\! & 
\| \taub_n\|^2 
\bigg( 
\frac{\widetilde{W}_n}{\sqrt{2}}
-
\frac{Z_n}{2\xi}
\bigg)
+ \frac{\| \taub_n\|^2 s}{2 \xi}
-\frac{1}{2}
 \| \taub_n\|^4 
\bigg(\frac{1}{2}
+
\frac{1}{4\xi^2}\bigg)
+o_{\rm P}(1) 
\label{L1}
\end{eqnarray}
as~$\ny$ under~${\rm P}^{(n)}_{\thetab_{n0},\kappa_{n,s}}$, hence, from contiguity, also under~${\rm P}^{(n)}_{\thetab_{n0},\kappa_n}$. Combining~\eqref{L2} and~\eqref{L1} establishes the local asymptotic quadraticity result in~(\ref{LAQdoubleiv}). Finally, the asymptotic normality result of~$\Deltab_n$ trivially follows from Lemma~\ref{LemJointNormality}.  
\cqfd
\vspace{3mm}

%%%%%%%%%%%%%%%%%%%%%%%%%%%%%%%%%%%%

The proof of Theorem~\ref{LANunspecregime5} requires both following lemmas.

\begin{lemma}
\label{LemmaLikeTh1}	
Let~$(p_n)$ be a sequence of integers that diverges to infinity. Let~$(\kappa_n)$ and~$(\kappa_{n*})$ be sequences in~$(0,\infty)$ that are~$o(p_n)$ and write~$e_{n1}$ and~$\tilde{e}_{n2}$ $($resp.,~$e_{n1*}$ and~$\tilde{e}_{n2*})$ for the corresponding moments based on~$\kappa_n$ $($resp., on~$\kappa_{n*})$. Let~$(\thetab_{n0})$, $(\nu_n)$ and~$(\taub_n)$ be sequences such that $\thetab_{n0}$ and~$\thetab_n=\thetab_{n0}+\nu_n\taub_n$ belong to~$\mathcal{S}^{p_n-1}$ for any~$n$, with~$(\taub_n)$ bounded and~$(\nu_n)$ such that 
\begin{equation}
\label{taillenulemma}
\nu_n^2 
=
O\Big(\frac{\sqrt{p_n}}{n \kappa_{n*} e_{n1*}}\Big)
.
\end{equation}
Then, with the same~$Z_n$ and~$\widetilde{W}_n$ as in Theorem~\ref{TheoLAQ}, we have that
	\begin{eqnarray*}
\lefteqn{	
\hspace{-8mm} 
\Lambda^{(n){\rm inv}}_{\thetab_n/\thetab_{n0};\kappa_{n*}}
\!\!=
-
\frac{1}{2} \sqrt{n} \kappa_{n*} \nu_n^2 \sqrt{\tilde{e}_{n2}}\, \|\taub_n\|^2  Z_n 
+ 
\frac{n \kappa_{n*}\nu_n^2 e_{n1*}}{\sqrt{2}p_n^{1/2}}
\|\taub_n\|^2 \Big(1- \frac{1}{4} \nu_n^2 \|\taub_n\|^2\Big)
  \widetilde{W}_n
}
\\[2mm]
& &  
\hspace{10mm} 
-
\frac{1}{8}
n \kappa_{n*}\nu_n^4 e_{n1*} \|\taub_n\|^4
-
\frac{n^2 \kappa_{n*}^2 \nu_n^4 e_{n1*}^2}{4p_n}  \|\taub_n\|^4 \Big(1- \frac{1}{4} \nu_n^2 \|\taub_n\| ^2\Big)^2 
\\[2mm]
& &  
\hspace{10mm} 
+
\frac{1}{2} n \kappa_{n*} \nu_n^2 (e_{n1*}-e_{n1}) \|\taub_n\|^2 
  +o_{\rm P}(1)
,
\end{eqnarray*}
as~$n\to\infty$ under~${\rm P}^{(n)}_{\thetab_{n0},\kappa_n}$.
\end{lemma}

{\sc Proof of Lemma~\ref{LemmaLikeTh1}.}
Since~$\kappa_n$ and~$\kappa_{n*}$ are both~$o(n)$, Lemma~\ref{LemEnFvML} ensures that~$f_{n2}/f_{n2*}=1+o(1)$, where~$f_{n2*}$ denotes the quantity~$f_{n2}$ based on~$\kappa_{n*}$. Using this, it can be showed along the exact same lines as in the proof of~(\ref{fromTheo1}) in Theorem~\ref{TheoLAQ} that, as~$n\to\infty$ under~${\rm P}^{(n)}_{\thetab_{n0},\kappa_n}$,  
\begin{eqnarray*}
\Lambda^{(n){\rm inv}}_{\thetab_n/\thetab_{n0};\kappa_{n*}}
&=&
n \kappa_{n*} \nu_n (\taub_n\pr \thetab_{n0}) \bar{\Xb}_{n}\pr \thetab_{n0} 
+
\frac{1}{2}
n \kappa_{n*} \nu_n^2 e_{n1*} \|\Pi_{\thetab_{n0}}\taub_n\|^2 T_n
%\label{fromTheo1}
\\[2mm]
& & 
\hspace{3mm} 
-
\frac{n^2 \kappa_{n*}^2 \nu_n^4 e_{n1*}^2}{4p_n}
 \|\Pi_{\thetab_{n0}}\taub_n\|^4 
+ o_{\rm P}(1)
,
\nonumber
\end{eqnarray*}
where $T_n:=1 +\sqrt{2}\,\widetilde{W}_n/\sqrt{p_n-1}$. If one replaces~$T_n$ by this expression and~$\bar{\Xb}_{n}\pr \thetab_{n0}$ by~$e_{n1}+\sqrt{\tilde{e}_{n2}}Z_n/\sqrt{n}$, then  the result follows by using~(\ref{ide11}), (\ref{ide12}), (\ref{taillenulemma}), and the tightness of~$\widetilde{W}_n$.  
\cqfd
\vspace{3mm} 

%%%%%%%%%%%%

The second lemma reinforces the variance result in Lemma~\ref{LemEnFvML}(iii). 

\begin{lemma}
	\label{Lemme2tilde}
Let~$(p_n)$ be a sequence of integers that diverges to infinity and~$(\kappa_n)$ be a sequence in~$(0,\infty)$ that is~$o(p_n)$. Denote as~$\tilde{e}_{n2}$ the variance of the distribution with probability density function~(\ref{pdfUFvML}). Then, $\sqrt{p_n\tilde{e}_{n2}}-1=O(\kappa_n^2/p_n^2)$ as~$n\to\infty$. 
\end{lemma}
\vspace{1mm}
 
{\sc Proof of Lemma~\ref{Lemme2tilde}.} 
In this proof, $C$ denotes a generic constant that may differ from line to line. Since (\ref{keymajorE2nUn}) yields
$$
\frac{
\sqrt{
\frac{3}{2}    
+
\frac{1}{p_n}   
-
 \sqrt{\big(\frac{1}{2}+\frac{1}{p_n}\big)^2 + \big(\frac{\kappa_n}{p_n}\big)^2}
}
}
{\frac{1}{2}-\frac{1}{p_n}+ \sqrt{\big(\frac{1}{2}+\frac{1}{p_n}\big)^2 + \big(\frac{\kappa_n}{p_n}\big)^2}} 
-
1
\leq
\sqrt{p_n\tilde{e}_{n2}}-1
\leq
 \frac{
1
}{
\sqrt{\frac{1}{2}+\sqrt{\big(\frac{1}{2}\big)^2 + \big(\frac{\kappa_n}{p_n}\big)^2}}
}
-
1
\leq
0
,
$$
we have
\begin{eqnarray*}
	\lefteqn{
|\sqrt{p_n\tilde{e}_{n2}}-1| 
\leq
1
-
\frac{
\sqrt{
\frac{3}{2}    
+
\frac{1}{p_n}   
-
 \sqrt{\big(\frac{1}{2}+\frac{1}{p_n}\big)^2 + \big(\frac{\kappa_n}{p_n}\big)^2}
}
}
{\frac{1}{2}-\frac{1}{p_n}+ \sqrt{\big(\frac{1}{2}+\frac{1}{p_n}\big)^2 + \big(\frac{\kappa_n}{p_n}\big)^2}} 
}
\\[2mm]
& & 
\hspace{3mm} 
\leq
C
\bigg\{
\Big(
{\textstyle{\frac{1}{2}}}-{\textstyle{\frac{1}{p_n}}}+ \sqrt{\big({\textstyle{\frac{1}{2}}}+{\textstyle{\frac{1}{p_n}}}\big)^2 + \big({\textstyle{\frac{\kappa_n}{p_n}}}\big)^2}
\Big)
-
\sqrt{
{\textstyle{\frac{3}{2}}}    
+
{\textstyle{\frac{1}{p_n}}}   
-
 \sqrt{\big({\textstyle{\frac{1}{2}}}+{\textstyle{\frac{1}{p_n}}}\big)^2 + \big({\textstyle{\frac{\kappa_n}{p_n}}}\big)^2}}
 \,
 \bigg\}
 \\[2mm]
& & 
\hspace{3mm} 
\leq
C
\bigg\{
\Big(
{\textstyle{\frac{1}{2}}}-{\textstyle{\frac{1}{p_n}}}+ \sqrt{\big({\textstyle{\frac{1}{2}}}+{\textstyle{\frac{1}{p_n}}}\big)^2 + \big({\textstyle{\frac{\kappa_n}{p_n}}}\big)^2}
\Big)^2
-
\Big(
{\textstyle{\frac{3}{2}}}    
+
{\textstyle{\frac{1}{p_n}}}   
-
 \sqrt{\big({\textstyle{\frac{1}{2}}}+{\textstyle{\frac{1}{p_n}}}\big)^2 + \big({\textstyle{\frac{\kappa_n}{p_n}}}\big)^2}
\Big)
\bigg\}
.
\end{eqnarray*}
Standard computations allow us to rewrite this upper-bound as
\begin{eqnarray*}
	\lefteqn{
|\sqrt{p_n\tilde{e}_{n2}}-1| 
%\leq
%
%$$
%=
%\Big(
%\frac{p_n-2}{2p_n}
%+ 
%\sqrt{\big(\frac{p_n+2}{2p_n}\big)^2 + \big(\frac{\kappa_n}{p_n}\big)^2}
%\Big)^2
%-
%\Big(
%\frac{3p_n+2}{2p_n}
%-
% \sqrt{\big(\frac{p_n+2}{2p_n}\big)^2 + \big(\frac{\kappa_n}{p_n}\big)^2}
%\Big)
%$$
%%
%$$
%=
%\frac{(p_n-2)^2}{4p_n^2}
%+
%\frac{p_n-2}{p_n}
%\sqrt{\big(\frac{p_n+2}{2p_n}\big)^2 + \big(\frac{\kappa_n}{p_n}\big)^2}
%+
%\frac{(p_n+2)^2}{4p_n^2}
%+
%\frac{\kappa_n^2}{p_n^2}
%-
%\frac{3p_n+2}{2p_n}
%+
% \sqrt{\big(\frac{p_n+2}{2p_n}\big)^2 + \big(\frac{\kappa_n}{p_n}\big)^2}
%$$
%%
%$$
%=
%\frac{(p_n-2)^2+(p_n+2)^2+4\kappa_n^2-(3p_n+2)2p_n}{4p_n^2}
%+
%\frac{2(p_n-1)}{p_n}
%\sqrt{\big(\frac{p_n+2}{2p_n}\big)^2 + \big(\frac{\kappa_n}{p_n}\big)^2}
%$$
%%
%$$
%=
%\frac{-4p_n^2-4p_n+8+4\kappa_n^2}{4p_n^2}
%+
%\frac{2(p_n-1)}{p_n}
%\sqrt{\big(\frac{p_n+2}{2p_n}\big)^2 + \big(\frac{\kappa_n}{p_n}\big)^2}
%$$
%%
%$$
%=
%\frac{2(p_n-1)}{p_n}
%\sqrt{\big(\frac{p_n+2}{2p_n}\big)^2 + \big(\frac{\kappa_n}{p_n}\big)^2}
%-
%\frac{4(p_n-1)(p_n+2)-4\kappa_n^2}{4p_n^2}
%$$
%
\leq
C
\bigg\{
{\textstyle{\frac{2(p_n-1)}{p_n}}}
\sqrt{\big({\textstyle{\frac{p_n+2}{2p_n}}}\big)^2 + \big({\textstyle{\frac{\kappa_n}{p_n}}}\big)^2}
-
{\textstyle{\frac{(p_n-1)(p_n+2)-\kappa_n^2}{p_n^2}}}
\bigg\}
}
\\[2mm]
& & 
\hspace{0mm} 
\leq
C
\bigg\{
{\textstyle{\frac{4(p_n-1)^2}{p_n^2}}}
\Big(
\big({\textstyle{\frac{p_n+2}{2p_n}}}\big)^2 + \big({\textstyle{\frac{\kappa_n}{p_n}}}\big)^2
\Big)
-
{\textstyle{\frac{((p_n-1)(p_n+2)-\kappa_n^2)^2}{p_n^4}}} 
\bigg\}
%
%$$
%=
%\frac{4(p_n-1)^2}{p_n^2}
%\Big(
%\frac{(p_n+2)^2+4\kappa_n^2}{4p_n^2} 
%\Big)
%-
%\frac{(p_n-1)^2(p_n+2)^2+\kappa_n^4-2\kappa_n^2(p_n-1)(p_n+2)}{p_n^4}
%$$
%%
%$$
%=
%\frac{4(p_n-1)^2((p_n+2)^2+4\kappa_n^2)-4(p_n-1)^2(p_n+2)^2-4\kappa_n^4+8\kappa_n^2(p_n-1)(p_n+2)}{4p_n^4}
%$$
%%
%$$
%=
%\frac{16(p_n-1)^2\kappa_n^2-4\kappa_n^4+8\kappa_n^2(p_n-1)(p_n+2)}{4p_n^4}
%$$
%%
%$$
%=
%\kappa_n^2\frac{4(p_n-1)^2-\kappa_n^2+2(p_n-1)(p_n+2)}{p_n^4}
%$$
%%
%$$
%=
%\kappa_n^2\frac{ 2(p_n-1) (2(p_n-1)+(p_n+2)) - \kappa_n^2}{p_n^4}
%$$
%%
%$$
%=
%\kappa_n^2\frac{ 6(p_n-1) p_n - \kappa_n^2}{p_n^4}
%$$
%%
%$$
%=
%\frac{6(p_n-1)p_n\kappa_n^2}{p_n^4}
%-
%\frac{\kappa_n^4}{p_n^4}
%$$
%%
%$$
%=
%\bigg(
%6
%-
%\frac{6}{p_n}
%-
%\frac{\kappa_n^2}{p_n^2}
%\bigg)
%\frac{\kappa_n^2}{p_n^2}
%$$
%
%\\[2mm]
%& & 
%\hspace{3mm} 
=
C
\Big(
6
-
{\textstyle{\frac{6}{p_n}}}
-
{\textstyle{\frac{\kappa_n^2}{p_n^2}}}
\Big)
{\textstyle{\frac{\kappa_n^2}{p_n^2}}}
\cdot
\end{eqnarray*}
which, for~$n$ large, is upper-bounded by
$
C\kappa_n^2/p_n^2
$, as was to be proved.
\cqfd
\vspace{3mm}

%%%%%%%%%%%%%%%%%%%%%%%%%%%%%%%%%%%%%%%%%%%%%%%%%%%%%%%%%%%%%%%%%

\noindent {\sc Proof of Theorem~\ref{LANunspecregime5}.}  
Since~$\kappa_n=o(p_n)$, Lemma~\ref{LemEnFvML}(iii) entails that
\begin{eqnarray}
\lefteqn{
\hspace{3mm} 
Z_n
=
\frac{\sqrt{n}(\bar{\Xb}_n\pr \thetab_{n0}- e_{n1})}{\sqrt{\tilde{e}_{2n}}}
=
\frac{\sqrt{n}(\bar{\Xb}_n\pr \thetab_{n0}- \xi r_n/\sqrt{n}+O(r_n^3 n^{-3/2}))}{\sqrt{\frac{1}{p_n}+o(\frac{1}{p_n})}}
}
\nonumber
\\[2mm]
& & 
\hspace{-9mm} 
=
\frac{\sqrt{p_n}(\sqrt{n}\bar{\Xb}_n\pr \thetab_{n0}- \xi r_n+O(r_n^3 n^{-1}))}{\sqrt{1+o(1)}}
%\nonumber
%\\[2mm]
%& & 
= 
\sqrt{np_n}\bar{\Xb}_n\pr \thetab_{n0}- \xi\sqrt{p_n} r_n
+
O_{\rm P}\Big(\frac{\sqrt{p_n}r_n^3}{n}\Big)
\label{repZnregime5}	
\end{eqnarray}
as~$n\to\infty$ under~${\rm P}^{(n)}_{\thetab_{n0},\kappa_n}$. 
Write then
\begin{eqnarray*}
\Lambda^{(n){\rm inv}}_{\thetab_n, \kappa_{n,s,\taub_n}/\thetab_{n0},\kappa_n} 
&=& 
\Lambda^{(n){\rm inv}}_{\thetab_n, \kappa_{n,s,\taub_n}/\thetab_{n0},\kappa_{n,s,\taub_n}}
+ 
\Lambda^{(n){\rm inv}}_{\thetab_{n0}, \kappa_{n,s,\taub_n}/\thetab_{n0},\kappa_n} 
\\[1mm]
&=& 
\Lambda^{(n){\rm inv}}_{\thetab_n, \kappa_{n,s,\taub_n}/\thetab_{n0},\kappa_{n,s,\taub_n}}
+ 
\log \frac{d{\rm P}\n_{\thetab_{n0},\kappa_{n,s,\taub_n}}}{d{\rm P}\n_{\thetab_{n0},\kappa_{n}}}
\\[2mm]
&\!\!=:\!\!& 
L_{n1}+L_{n2}.
\end{eqnarray*}
Letting~$\rho_n:=(1-{\textstyle{\frac{1}{2}}}\nu_n^2\|\taub_n\|^2)^{-1}$ and using~(\ref{repZnregime5}), we obtain
\begin{eqnarray}
%\label{LLn1}
L_{n2} 
& = &
 n \big( \log(c_{p_n,\kappa_{n,s,\taub_n}})-\log(c_{p_n, \kappa_n}) \big)
 + 
 n(\kappa_{n,s,\taub_n}-\kappa_n) \bar{\Xb}_n\pr \thetab_{n0} 
 \nonumber 
 \\[2.5mm]
& = & 
n \Big[
\log\Big(\frac{c_{p_n,0}}{c_{p_n, \kappa_n}} \Big)
-
\log\Big(\frac{c_{p_n,0}}{c_{p_n, \kappa_{n,s,\taub_n}}} \Big)
\Big]
 + 
\rho_n 
(s\sqrt{np_n}+ {\textstyle{\frac{1}{2}}}\xi\sqrt{n} p_n r_n\nu_n^2\|\taub_n\|^2)\, \bar{\Xb}_n\pr \thetab_{n0} 
\nonumber 
\\
& = & 
n 
\Big[
\log H_{\frac{p_n}{2}-1} (\kappa_n) 
-
\log H_{\frac{p_n}{2}-1} (\kappa_{n,s,\taub_n})
\Big]
\nonumber 
\\[1mm]
& & 
\hspace{3mm} 
+ 
\rho_n 
(s\sqrt{np_n}+{\textstyle{\frac{1}{2}}} \xi\sqrt{n} p_n r_n \nu_n^2\|\taub_n\|^2)\, 
\bigg( 
\frac{Z_{n}}{\sqrt{np_n}} + \frac{\xi r_n}{\sqrt{n}} +O_{\rm P}\Big(\frac{r_n^3}{n^{3/2}}\Big)
\bigg)
\nonumber 
\end{eqnarray}
as~$n\to\infty$ under~${\rm P}^{(n)}_{\thetab_{n0},\kappa_n}$. Since~$\sqrt{p_n} r_n^3$ and~$p_n r_n^4 \nu_n^2$ are~$o(n)$, this yields 
%$$
%\frac{p_n r_n^4 \nu_n^2}{n}
%=
%\frac{p_n r_n^4}{n}\,\Big(\frac{p_n^{3/2}}{n\kappa_n^2}\Big)
%=
%\frac{p_n r_n^4}{n}\,\Big(\frac{p_n^{3/2}n}{n\,p_n^2r_n^2}\Big)
%=
%\frac{p_n r_n^4}{n}\,\Big(\frac{1}{p_n^{1/2}r_n^2}\Big)
%=
%\frac{p_n^{1/2} r_n^2}{n}
%=
%o(1)
%$$
%in case~(a) (it is trivially true in cases~(b)--(c)),
\begin{eqnarray}
\label{LLn1}
L_{n2} 
%& = &
%n 
%\Big[
%\log H_{\frac{p_n}{2}-1} (\kappa_n) 
%-
%\log H_{\frac{p_n}{2}-1} (\kappa_{n,s,\taub_n})
%\Big]
%\nonumber 
%\\[1mm]
%& & 
%\hspace{-3mm} 
%+ 
%\rho_n (s + {\textstyle{\frac{1}{2}}} \xi \sqrt{p_n}r_n\nu_n^2\|\taub_n\|^2) Z_n  
%+ \xi s \rho_n r_n\sqrt{p_n}
%+ {\textstyle{\frac{1}{2}}} \xi^2\rho_n p_n r_n^2\nu_n^2\|\taub_n\|^2 
%+ O_{\rm P}\Big(\frac{p_n r_n^4 \nu_n^2}{n}\Big)
%+ o_{\rm P}(1)
% \nonumber 
%\\[2mm]
& = & 
n 
\Big[
\log H_{\frac{p_n}{2}-1} (\kappa_n) 
-
\log H_{\frac{p_n}{2}-1} (\kappa_{n,s,\taub_n})
\Big]
\nonumber 
\\[1mm]
& & 
\hspace{-3mm} 
+ 
\rho_n \big(s + {\textstyle{\frac{1}{2}}} \xi \sqrt{p_n}r_n\nu_n^2\|\taub_n\|^2\big) Z_n 
+ \xi s \rho_n r_n\sqrt{p_n}
+ {\textstyle{\frac{1}{2}}} \xi^2\rho_n p_n r_n^2\nu_n^2\|\taub_n\|^2 
+ o_{\rm P}(1)
 \nonumber 
\\[2mm]
& =: & 
\tilde{L}_{n2}
+
\rho_n (s + {\textstyle{\frac{1}{2}}} \xi \sqrt{p_n}r_n\nu_n^2\|\taub_n\|^2) Z_n
  +
o_{\rm P}(1)
\nonumber 
\end{eqnarray}
as~$n\to\infty$ under~${\rm P}^{(n)}_{\thetab_{n0},\kappa_n}$.
%, where~$\delta=1$ in case~(a), $C^2$ in case~(b), and~$0$ in case~(c). 
%In~(a),
%$$
%\sqrt{p_n} r_n^2\nu_n^2
%=
%\sqrt{p_n} r_n^2 \frac{p_n^{3/2}}{n\kappa_n^2}
%=
%\sqrt{p_n} r_n^2 \frac{p_n^{3/2}n}{n p_n^2 r_n^2 \xi^2}
%=
%\frac{1}{\xi^2}
%$$
Since $p_n=o(n^2 r_n^{-4})$, we can apply Lemma \ref{LemGeneralDLFvML} with~$a_n=\sqrt{n}$ and $T_n \equiv 1$, which provides 
%$$
%\textcolor{red}{ 
%\frac{b_n^6}{a_n^4\nu_n^5}
%=
%\frac{(\sqrt{n}\kappa_n)^6}{n^2p^5}
%=
%\frac{n \kappa_n^6}{p^5}
%=
%\frac{n (p_n r_n
%(\xi+s/(\sqrt{p_n}r_n)+\nu_n^2\|\taub_n\|^2/4)/\sqrt{n})^6}{p^5}
%}
%$$
%$$
%\textcolor{red}{ 
%=
%\frac{p_nr_n^6}{n^2}\,(\xi+s/(\sqrt{p_n}r_n)+\nu_n^2\|\taub_n\|^2/4))^6
%=
%(\xi+s/(\sqrt{p_n}r_n)+\nu_n^2\|\taub_n\|^2/4))^6 o(1)
%=
%o(1)
%}
%$$
\begin{eqnarray}
%\label{LLn2}
\tilde{L}_{n2} 
&=& 
\Big(
\frac{n \kappa_n^2}{4({\textstyle{\frac{p_n}{2}}}-1)}
-
\frac{n \kappa_n^4}{32({\textstyle{\frac{p_n}{2}}}-1)^3}
\Big)
-
\Big(
\frac{n \kappa_{n,s,\taub_n}^2}{4({\textstyle{\frac{p_n}{2}}}-1)}
-
\frac{n \kappa_{n,s,\taub_n}^4}{32({\textstyle{\frac{p_n}{2}}}-1)^3}
\Big)
\nonumber 
\\[2mm]
& & 
\hspace{3mm} 
+ \xi s \rho_n r_n\sqrt{p_n}
+ {\textstyle{\frac{1}{2}}} \xi^2\rho_n p_n r_n^2\nu_n^2\|\taub_n\|^2 
+
o(1)
\nonumber  
\\[2mm]
&=& 
-\frac{n(\kappa_{n,s,\taub_n}^2-\kappa_n^2)}{2p_n-4}
+ 
\frac{n(\kappa_{n,s,\taub_n}^4-\kappa_n^4)}{4(p_n-2)^3}
+ \xi s \rho_n r_n\sqrt{p_n}
+ {\textstyle{\frac{1}{2}}} \xi^2\rho_n p_n r_n^2\nu_n^2\|\taub_n\|^2 
+
o(1)
\nonumber  
\\[2mm]
&=& 
-\frac{p_n^2r_n^2
\{\rho_n^2 (\xi+s/(\sqrt{p_n}r_n))^2-\xi^2\}}{2p_n-4}
+ S_n
%\nonumber 
%\\[2mm]
%& & 
%\hspace{3mm} 
+ \xi s \rho_n r_n\sqrt{p_n}
+ {\textstyle{\frac{1}{2}}} \xi^2\rho_n p_n r_n^2\nu_n^2\|\taub_n\|^2 
+o(1) 
\nonumber 
,
\end{eqnarray}
where, since~$1-\rho_n=-\frac{1}{2}\rho_n\nu_n^2\|\taub_n\|^2$, 
\begin{eqnarray*}
\lefteqn{	
\hspace{2mm} 
S_n
=
\frac{p_n^4r_n^4\{\rho_n^4 (\xi+s/(\sqrt{p_n}r_n))^4-\xi^4\}}{4n(p_n-2)^3}
%=
%\frac{p_n^4r_n^4}{4n(p_n-2)^3}
%\Big( 
%\rho_n^4 
%\Big(\xi+ \frac{s}{\sqrt{p_n}r_n}\Big)^4
%-\xi^4
%\Big)
=
\frac{p_n^4r_n^4}{4n(p_n-2)^3}
\Big( 
(\rho_n^4-1) \xi^4 
+ {\textstyle{ \sum_{\ell=1}^4 {4 \choose \ell} \frac{s^\ell \xi^{4-\ell}}{(\sqrt{p_n}r_n)^\ell} }}
\Big)
}
\\[2mm]
& & 
\hspace{-1mm} 
=
O\Big( \frac{p_n r_n^4}{n}
(\rho_n-1)(\rho_n+1)(\rho_n^2+1) 
\Big)
+ 
O\Big( \frac{\sqrt{p_n} r_n^3}{n} \Big)
=
O\Big( \frac{p_n r_n^4 \nu_n^2}{n} \Big)
+ 
O\Big( \frac{\sqrt{p_n} r_n^3}{n} \Big)
=
o(1)
.
\end{eqnarray*}
Thus, using the identities~$1-\rho_n=-\frac{1}{2}\rho_n\nu_n^2\|\taub_n\|^2$ and~$\rho_n^2-1=\rho_n^2\nu_n^2\|\taub_n\|^2-\frac{1}{4}\rho_n^2\nu_n^4\|\taub_n\|^4$, we have
\begin{eqnarray}
\label{LLn2}
\tilde{L}_{n2} 
&=& 
-\frac{(\rho_n^2-1)\xi^2 p_n^2 r_n^2}{2p_n-4}
-\frac{s^2 \rho_n^2 p_n}{2p_n-4}
-\frac{2\xi s\rho_n^2 p_n^{3/2}r_n}{2p_n-4}
+ \xi s \rho_n r_n\sqrt{p_n}
+ {\textstyle{\frac{1}{2}}} \xi^2\rho_n p_n r_n^2\nu_n^2\|\taub_n\|^2 
+o(1) 
\nonumber 
%
%\\[2mm]
%&=& 
%-\frac{\xi^2\rho_n^2 p_n^2 r_n^2 \nu_n^2\|\taub_n\|^2}{2p_n-4}
%+\frac{\xi^2\rho_n^2 p_n^2 r_n^2 \nu_n^4\|\taub_n\|^4}{4(2p_n-4)}
%-\frac{s^2\rho_n^2 p_n}{2p_n-4}
%-\frac{2\xi s\rho_n^2 p_n^{3/2}r_n}{2p_n-4}
%+ \xi s \rho_n r_n\sqrt{p_n}
%+ {\textstyle{\frac{1}{2}}} \xi^2\rho_n p_n r_n^2\nu_n^2\|\taub_n\|^2 
%+o(1) 
%\nonumber 
%
\\[2mm]
&=& 
-\frac{1}{2} \xi^2\rho_n^2 p_n r_n^2 \nu_n^2\|\taub_n\|^2 
+\frac{1}{8} \xi^2\rho_n^2 p_n r_n^2 \nu_n^4\|\taub_n\|^4 
-\frac{1}{2} s^2 \rho_n^2 
\\[2mm]
& & 
\hspace{20mm} 
- \xi s \rho_n^2 \sqrt{p_n}r_n
+ \xi s \rho_n r_n\sqrt{p_n}
+ {\textstyle{\frac{1}{2}}} \xi^2\rho_n p_n r_n^2\nu_n^2\|\taub_n\|^2 
+o(1) 
\nonumber 
\\[2mm]
&=& 
\frac{1}{2} \xi^2\rho_n(1-\rho_n) p_n r_n^2 \nu_n^2\|\taub_n\|^2 
+\frac{1}{8} \xi^2\rho_n^2 p_n r_n^2 \nu_n^4\|\taub_n\|^4 
-\frac{1}{2} s^2\rho_n^2 
+\xi s\rho_n(1-\rho_n)\sqrt{p_n}r_n 
+o(1) 
\nonumber 
\\[2mm]
&=& 
-\frac{1}{2} s^2 \rho_n^2
-\frac{1}{2}\xi s \rho_n^2 \sqrt{p_n}r_n \nu_n^2 \|\taub_n\|^2
-\frac{1}{8} \xi^2 \rho_n^2 p_n r_n^2 \nu_n^4\|\taub_n\|^4
+o(1) 
\nonumber 
\end{eqnarray}
as $\ny$. Therefore, we proved that, as~$n\to\infty$ under~${\rm P}^{(n)}_{\thetab_{n0},\kappa_n}$,
\begin{eqnarray} 
L_{n2}
&=&
\rho_n (s + {\textstyle{\frac{1}{2}}} \xi \sqrt{p_n}r_n\nu_n^2\|\taub_n\|^2) Z_n
\label{L2}
\\[2mm]
& &
\hspace{5mm} 
-\frac{1}{2} s^2 \rho_n^2
-\frac{1}{2}\xi s \rho_n^2 \sqrt{p_n}r_n \nu_n^2 \|\taub_n\|^2
-\frac{1}{8} \xi^2 \rho_n^2 p_n r_n^2 \nu_n^4\|\taub_n\|^4
+
o_{\rm P}(1) 
.
\nonumber 
%\\[2mm]
%&=&
%\rho_n (s + {\textstyle{\frac{1}{2}}} \xi \sqrt{p_n}r_n\nu_n^2\|\taub_n\|^2) Z_n
%-
%{\textstyle{\frac{1}{2}}}
%\rho_n^2 
%\big(
%s 
%+
%{\textstyle{\frac{1}{2}}} \xi \sqrt{p_n}r_n\nu_n^2\|\taub_n\|^2
%\big)^2
%+
%o_{\rm P}(1) 
%.
%\nonumber 
\end{eqnarray}
%Note that the Le Cam first lemma implies that ${\rm P}^{(n)}_{\thetab_{n0},\kappa_{n,s,\taub_n}}$ and ${\rm P}^{(n)}_{\thetab_{n0},\kappa_n}$ are mutually contiguous sequences of probability measures. 
%Note that
%$$
%p_n r_n^2 \nu_n^4
%=
%p_n r_n^2\,\Big(\frac{p_n^{3}}{n^2\kappa_n^4}\Big)
%=
%p_n r_n^2\,\Big(\frac{p_n^3 n^2}{n^2\,p_n^4r_n^4\xi^4}\Big)
%=
%\frac{1}{r_n^2\xi^4}
%\cdot
%$$

We turn to~$L_{1n}$. Write~$c_{n,s,\taub_n}:=n\nu_n^2 \kappa_{n,s,\taub_n} e_{n1,s,\taub_n}/\sqrt{p_n}$, where~$e_{n1,s,\taub_n}$ and~$\tilde{e}_{2n,s,\taub_n}$ denote the values of~$e_{n1}$ and~$\tilde{e}_{2n}$ under~${\rm P}^{(n)}_{\thetab_{n0},\kappa_{n,s,\taub_n}}\!$. Since~$\sqrt{p_n} r_n^2 \nu_n^2=O(1)$, (\ref{tocallth1inreg5}) below ensures that~$c_{n,s,\taub_n}$ is~$O(1)$. Therefore, Lemma~\ref{LemmaLikeTh1} yields 
\begin{equation}
\label{decompL1}
L_{1n}
=	
L^Z_{1n}+\tilde{L}^Z_{1n}+\bar{L}^Z_{1n}+L^W_{1n}+o_{\rm P}(1)
,
\end{equation}
where we let
$$
L^Z_{1n}
:=
-
\frac{1}{2} \sqrt{n} \kappa_{n,s,\taub_n} \nu_n^2 \sqrt{\tilde{e}_{n2}}\, \|\taub_n\|^2  Z_{n}
,
\quad
\tilde{L}^Z_{1n}
:=
-
\frac{1}{8} \sqrt{p_n}\nu_n^2 c_{n,s,\taub_n} \|\taub_n\|^4,
$$
$$
\bar{L}^Z_{1n}
:=
\frac{1}{2} n \kappa_{n,s,\taub_n} \nu_n^2 (e_{n1,s,\taub_n}-e_{n1}) \|\taub_n\|^2 
,
$$
and
$$
L^W_{1n}
:=
\frac{1}{\sqrt{2}} c_{n,s,\taub_n} \|\taub_n\|^2 \Big(1- \frac{1}{4} \nu_n^2 \|\taub_n\|^2\Big) \widetilde{W}_n
-
\frac{1}{4} c_{n,s,\taub_n}^2  \|\taub_n\|^4 \Big(1- \frac{1}{4} \nu_n^2 \|\taub_n\| ^2\Big)^2 
.
$$
Lemma~\ref{Lemme2tilde} provides 
\begin{eqnarray}
	\lefteqn{
L^Z_{1n}
=
-
\frac{\sqrt{n} \nu_n^2}{2\sqrt{p_n}}
 \|\taub_n\|^2  
\Big( 
\frac{\rho_n p_n r_n}{\sqrt{n}}
\Big(
\xi
+
\frac{s}{\sqrt{p_n}r_n}
\Big)
\Big)
\sqrt{p\tilde{e}_{n2}}\, Z_n
}
\nonumber
\\[2mm]
 & &
 \hspace{3mm} 
=
-
\frac{1}{2} \rho_n \sqrt{p_n} r_n \nu_n^2 \|\taub_n\|^2  
\Big(
\xi
+
\frac{s}{\sqrt{p_n}r_n}
\Big)
\Big(
1+O\Big(\frac{\kappa_n^2}{p_n^2}\Big)
\Big)
Z_n
%
%\nonumber
%\\[2mm]
% & &
% \hspace{3mm} 
%=
%-
%\frac{1}{2} \rho_n \sqrt{p_n} r_n \nu_n^2 \|\taub_n\|^2  
%\Big(
%\xi
%+
%\frac{s}{\sqrt{p_n}r_n}
%+
%O\Big(\frac{\kappa_n^2}{p_n^2}\Big)
%\Big)
%Z_n
%%
%\nonumber
%\\[2mm]
% & &
% \hspace{3mm} 
%=
%- \frac{1}{2} \xi \rho_n \sqrt{p_n} r_n \nu_n^2 \|\taub_n\|^2 Z_n 
%- \frac{1}{2} s \rho_n \nu_n^2 \|\taub_n\|^2 Z_n 
%+ 
%O_{\rm P}\Big(\frac{\kappa_n^2 r_n \nu_n^2}{p_n^{3/2}}\Big)
%
\nonumber
\\[2mm]
 & &
 \hspace{3mm} 
=
- \frac{1}{2} \xi \rho_n \sqrt{p_n} r_n \nu_n^2 \|\taub_n\|^2 Z_n 
- \frac{1}{2} s \rho_n \nu_n^2 \|\taub_n\|^2 Z_n 
+ 
o_{\rm P}(1)
%O_{\rm P}\Big(\frac{\sqrt{p_n} r_n^3 \nu_n^2}{n}\Big)
,
\label{L1Z}
\end{eqnarray}
where we used the fact that~$\sqrt{p_n}r_n^3 \nu_n^2$ is~$o(n)$. 

Now, Lemma~\ref{LemEnFvML}(iii) yields 
\begin{eqnarray}
	\lefteqn{
\hspace{-23mm} 
	e_{n1,s,\taub_n}
=
\frac{\kappa_{n,s,\taub_n}}{p_n}
+
O\Big(\frac{\kappa^3_{n,s,\taub_n}}{p_n^3}\Big)
=
\frac{\rho_n r_n}{\sqrt{n}}
\Big(
\xi
+
\frac{s}{\sqrt{p_n}r_n}
\Big)
+
O\Big( \frac{r_n^3}{n^{3/2}} \Big)
}
\nonumber
\\[2mm]
& & 
\hspace{3mm} 
=
\frac{\xi\rho_n r_n}{\sqrt{n}} 
+
\frac{s\rho_n}{\sqrt{np_n}}
+
O\Big( \frac{r_n^3}{n^{3/2}} \Big)
,
\label{expande1stau}
\end{eqnarray}
so that
\begin{eqnarray}
\lefteqn{
\hspace{-17mm} 
c_{n,s,\taub_n}
=
\frac{n \nu_n^2}{\sqrt{p_n}}
\Big( 
\frac{\rho_n p_n r_n}{\sqrt{n}}
\Big(
\xi
+
\frac{s}{\sqrt{p_n}r_n}
\Big)
\Big)
\Big( 
\frac{\xi\rho_n r_n}{\sqrt{n}} 
+
\frac{s\rho_n}{\sqrt{np_n}}
+
O\Big( \frac{r_n^3}{n^{3/2}} \Big)
\Big)
}
\nonumber
%\\[2mm]
%& & 
%\hspace{-10mm} 
%=
%\rho_n \sqrt{n p_n} r_n \nu_n^2
%\Big(
%\xi
%+
%\frac{s}{\sqrt{p_n}r_n}
%\Big)
%\Big( 
%\frac{\xi\rho_n r_n}{\sqrt{n}} 
%+
%\frac{s\rho_n}{\sqrt{np_n}}
%+
%O\Big( \frac{r_n^3}{n^{3/2}} \Big)
%\Big)
%\nonumber
\\[2mm]
& & 
\hspace{-10mm} 
=
\rho_n^2 \sqrt{p_n} r_n^2 \nu_n^2
\Big(
\xi
+
\frac{s}{\sqrt{p_n}r_n}
\Big)
\Big( 
\xi
+
\frac{s}{\sqrt{p_n}r_n}
+
O\Big( \frac{r_n^2}{n} \Big)
\Big)
\label{tocallth1inreg5}
%\\[2mm]
%& & 
%\hspace{-10mm} 
%=
%\xi^2 \rho_n^2 \sqrt{p_n} r_n^2 \nu_n^2 
%+ 2 \xi s \rho_n^2 r_n \nu_n^2 
%+ s^2 \frac{\rho_n^2 \nu_n^2}{\sqrt{p_n}}
%+
%O\Big( \frac{\sqrt{p_n} r_n^4 \nu_n^2}{n} \Big)
%\nonumber
%\\[2mm]
%& & 
%\hspace{-10mm} 
%=
%\rho_n^2 \delta
%+ 2 \xi s \rho_n^2 r_n \nu_n^2 
%+ s^2 \frac{\rho_n^2 \nu_n^2}{\sqrt{p_n}}
%+
%O\Big( \frac{\sqrt{p_n} r_n^4 \nu_n^2}{n} \Big)
%O\Big( \frac{r_n^2}{n} \Big)
%\nonumber
,
\end{eqnarray}
%where~$\delta=1$ in case~(a), $C$ in case~(b), and~$0$ in case~(c). 
%
which in turn implies that
\begin{eqnarray}
	\tilde{L}^Z_{1n}
%=
%-
%\frac{1}{8} \sqrt{p_n}\nu_n^2 c_{n,s,\taub_n} \|\taub_n\|^4
&=&
-
\frac{1}{8}
\sqrt{p_n}\nu_n^2
\Big(
\xi^2 \rho_n^2 \sqrt{p_n} r_n^2 \nu_n^2
+
2\xi s\rho_n^2 r_n \nu_n^2
+
s^2 \frac{\rho_n^2 \nu_n^2}{\sqrt{p_n}}
+
O\Big( \frac{\sqrt{p_n} r_n^4 \nu_n^2}{n} \Big)
\Big)
\|\taub_n\|^4
\nonumber
\\[2mm]
&=&
- \frac{1}{8} \xi^2 \rho_n^2 p_n r_n^2 \nu_n^4 \|\taub_n\|^4
- \frac{1}{4} \xi s\rho_n^2 \sqrt{p_n} r_n \nu_n^4 \|\taub_n\|^4
- \frac{1}{8} s^2\rho_n^2 \nu_n^4 \|\taub_n\|^4
+ o(1)
.
\quad
\label{Ltilde1Z}
\end{eqnarray}

Using~(\ref{expande1stau}) and applying Lemma~\ref{LemEnFvML}(iii) again, we obtain 
\begin{eqnarray*}
\lefteqn{
e_{n1,s,\taub_n}-e_{n1}
=
\bigg(
\frac{\xi\rho_n r_n}{\sqrt{n}} 
+
\frac{s\rho_n}{\sqrt{np_n}}
+
O\Big( \frac{r_n^3}{n^{3/2}} \Big)
\bigg)
-
\frac{\kappa_{n}}{p_n}
+
O\Big(\frac{\kappa^3_{n}}{p_n^3}\Big)
}
\\[2mm]
%&=&
& & 
\hspace{3mm} 
=
\frac{\xi(\rho_n-1) r_n}{\sqrt{n}} 
+
\frac{s\rho_n}{\sqrt{np_n}}
+
O\Big( \frac{r_n^3}{n^{3/2}} \Big)
%\\[2mm]
%&=&
=
\frac{\xi\rho_n r_n\nu_n^2}{2\sqrt{n}} \|\taub_n\|^2
+
\frac{s\rho_n}{\sqrt{np_n}}
+
O\Big( \frac{r_n^3}{n^{3/2}} \Big)
,
\end{eqnarray*}
so that
\begin{eqnarray}
\lefteqn{
\bar{L}^Z_{1n}
=
\frac{1}{2} n 
\Big( 
\frac{\rho_n p_n r_n}{\sqrt{n}}
\Big(
\xi
+
\frac{s}{\sqrt{p_n}r_n}
\Big)
\Big)
 \nu_n^2 
\bigg(
\frac{\xi\rho_n r_n\nu_n^2}{2\sqrt{n}} \|\taub_n\|^2
+
\frac{s\rho_n}{\sqrt{np_n}}
+
O\Big( \frac{r_n^3}{n^{3/2}} \Big)
\bigg)
\|\taub_n\|^2 
}
\nonumber
%\\[2mm]
%& & 
%\hspace{3mm} 
%=
%\frac{1}{2} \rho_n p_n r_n \sqrt{n}  \nu_n^2 \|\taub_n\|^2
%\Big(
%\xi
%+
%\frac{s}{\sqrt{p_n}r_n}
%\Big)
%\bigg(
%\frac{\xi\rho_n r_n\nu_n^2}{2\sqrt{n}} \|\taub_n\|^2
%+
%\frac{s\rho_n}{\sqrt{np_n}}
%+
%O\Big( \frac{r_n^3}{n^{3/2}} \Big)
%\bigg)
%\nonumber
\\[2mm]
& & 
\hspace{3mm} 
=
\frac{1}{2} \rho_n p_n r_n \sqrt{n}  \nu_n^2 \|\taub_n\|^2
\xi
\bigg(
\frac{\xi\rho_n r_n\nu_n^2}{2\sqrt{n}} \|\taub_n\|^2
+
\frac{s\rho_n}{\sqrt{np_n}}
+
O\Big( \frac{r_n^3}{n^{3/2}} \Big)
\bigg)
\nonumber
\\[2mm]
& & 
\hspace{13mm} 
+
\frac{1}{2} \rho_n p_n r_n \sqrt{n}  \nu_n^2 \|\taub_n\|^2
\frac{s}{\sqrt{p_n}r_n}
\bigg(
\frac{\xi\rho_n r_n\nu_n^2}{2\sqrt{n}} \|\taub_n\|^2
+
\frac{s\rho_n}{\sqrt{np_n}}
+
O\Big( \frac{r_n^3}{n^{3/2}} \Big)
\bigg)
\nonumber
\\[2mm]
& & 
\hspace{3mm} 
=
\frac{1}{4}\xi^2 \rho_n^2 p_n r_n^2 \nu_n^4 \|\taub_n\|^4 
+
\frac{1}{2}\xi s\rho_n^2 \sqrt{p_n} r_n \nu_n^2 \|\taub_n\|^2 
\nonumber
\\[2mm]
& & 
\hspace{13mm} 
+
\frac{1}{4} \xi s \rho_n^2 \sqrt{p_n} r_n \nu_n^4 \|\taub_n\|^4
+
\frac{1}{2} s^2 \rho_n^2 \nu_n^2 \|\taub_n\|^2
+
o(1)
,
\label{Lbar1Z}
\end{eqnarray}
where we used the facts that~$p_n r_n^4 \nu_n^2$ and~$\sqrt{p_n}r_n^3 \nu_n^2$ are~$o(n)$. 

Jointly with~(\ref{L2}), (\ref{L1Z}), (\ref{Ltilde1Z}), and~(\ref{Lbar1Z}), this shows that
%$$
%	\Lambda^{(n){\rm inv}}_{\thetab_n, \kappa_{n,s,\taub_n}/\thetab_{n0},\kappa_n}
%-
%L^W_{1n}
%=
%s \rho_n  (1 - \frac{1}{2} \nu_n^2 \|\taub_n\|^2) Z_n
%- \frac{1}{2} s^2 \rho_n^2
%(1- \nu_n^2 \|\taub_n\|^2+\frac{1}{4} \nu_n^4 \|\taub_n\|^4)
%+ o_{\rm P}(1)
%$$
\begin{eqnarray*}
	\Lambda^{(n){\rm inv}}_{\thetab_n, \kappa_{n,s,\taub_n}/\thetab_{n0},\kappa_n}
-
L^W_{1n}
&=&
s \rho_n  (1 - {\textstyle{\frac{1}{2}}} \nu_n^2 \|\taub_n\|^2) Z_n
- \frac{1}{2} s^2 \rho_n^2
(1- \nu_n^2 \|\taub_n\|^2+{\textstyle{\frac{1}{4}}} \nu_n^4 \|\taub_n\|^4)
+ o_{\rm P}(1)
\\[2mm]
&=&
s Z_n
- \frac{1}{2} s^2
+ o_{\rm P}(1)
.
\end{eqnarray*}
The result thus follows from the definition of
$
L^W_{1n}
$
 and the fact that~(\ref{tocallth1inreg5}) implies that 
 $c_{n,s,\taub_n}=1+o(1)$ in case~(a), 
 $c_{n,s,\taub_n}=\xi^2+o(1)$ in case~(b), and
 $c_{n,s,\taub_n}=o(1)$ in case~(c) (in each case, the asymptotic normality result of~$\Deltab_n$ follows from Lemma~\ref{LemJointNormality}). 
 \cqfd

%%%%%%%%%%%%%%%%%%%%%%%%%%
%%%%%%%%%%%%%%%%%%%%%%%%%%
%%%%%%%%%%%%%%%%%%%%%%%%%%

\section{Technical proofs for Section~3}
\label{SecSupProofBilling}

%%%%%%%%%%%%%%%%%%%%%%%%%%
%%%%%%%%%%%%%%%%%%%%%%%%%%
%%%%%%%%%%%%%%%%%%%%%%%%%%

The proof of Theorem~\ref{TheoBillingH1} requires the following lemma.

\begin{lemma}
\label{LemTrMn}
	Let~${\bf M}_n:=\thetab_n\thetab_n\pr-\thetab_{n0}\thetab_{n0}\pr$, where~$(\thetab_n)$ and~$(\thetab_{n0})$ belong to~$\mathcal{S}^{p_n-1}$. Then, for any real numbers~$a,b,c,d$, we have that
	$
	{\rm tr}\big[
 {\bf M}_n^\ell 
 ( a \thetab_n\thetab_n\pr + b ({\bf I}_{p_n}- \thetab_n\thetab_n\pr) )
 {\bf M}_n^\ell
 ( c\, \thetab_n\thetab_n\pr + d ({\bf I}_{p_n}- \thetab_n\thetab_n\pr) )
  \big] 
$
is equal to~$(ad+bc) (1 - (\thetab_{n0}\pr \thetab_n)^2)+ (a-b)(c-d) (1 - (\thetab_{n0}\pr \thetab_n)^2)^2$ for~$\ell=1$ and to~$(ac+bd) (1 - (\thetab_{n0}\pr \thetab_n)^2)^2$ for~$\ell=2$.  
\end{lemma}

{\sc Proof of Lemma~\ref{LemTrMn}.}
%We have ${\rm tr}[{\bf M}_n]=0$. 
Direct computations yield
$$
{\bf M}_n^2=
\thetab_n\thetab_n\pr+ \thetab_{n0}\thetab_{n0}\pr 
- (\thetab_n\pr\thetab_{n0}) \thetab_{n0}\thetab_n\pr
- (\thetab_n\pr\thetab_{n0}) \thetab_n\thetab_{n0}\pr
\quad
\textrm{ and }
\quad
{\bf M}_n^4
=
(1 - (\thetab_{n0}\pr \thetab_n)^2)
{\bf M}_n^2
.
$$
This provides
$
{\rm tr}[{\bf M}^2_n]
=
2(1 - (\thetab_{n0}\pr \thetab_n)^2)
%=
%2 \nu_n^2 \|\taub_n\|^2 
% (1- {\textstyle{\frac{1}{4} \nu_n^2 \|\taub_n\|^2}})
$
and
$
{\rm tr}[{\bf M}^4_n]
=
2(1 - (\thetab_{n0}\pr \thetab_n)^2)^2
%=
%2 \nu_n^4 \|\taub_n\|^4 
% (1- {\textstyle{\frac{1}{4} \nu_n^2 \|\taub_n\|^2}})^2
,
$
and allows to show that 
$\thetab_n\pr{\bf M}_n^2\thetab_n
%=
%\thetab_n\pr  (\thetab_n\thetab_n\pr-\thetab_{n0}\thetab_{n0}\pr)^2 \thetab_n
%=
%\thetab_n\pr  
%\{
%\thetab_n\thetab_n\pr+ \thetab_{n0}\thetab_{n0}\pr 
%- (\thetab_n\pr\thetab_{n0}) \thetab_{n0}\thetab_n\pr
%- (\thetab_n\pr\thetab_{n0}) \thetab_n\thetab_{n0}\pr
%\} 
%\thetab_n
%=
%1
%+
%(\thetab_n\pr\thetab_{n0})^2
%-
%2
% (\thetab_n\pr\thetab_{n0})^2
=
1
-
(\thetab_n\pr\thetab_{n0})^2
%=
%\thetab_n\pr{\bf M}_n\thetab_n
%=
% \nu_n^2 \|\taub_n\|^2 
% (1-\frac{1}{4} \nu_n^2 \|\taub_n\|^2)
$
and
$
\thetab_n\pr{\bf M}_n^4\thetab_n
%=
%(1 - (\thetab_{n0}\pr \thetab_n)^2) \thetab_n\pr{\bf M}_n^2\thetab_n
=
(1 - (\thetab_{n0}\pr \thetab_n)^2)^2
%=
%(\thetab_n\pr{\bf M}_n\thetab_n)^2
%=
% \nu_n^4 \|\taub_n\|^4 
% (1-\frac{1}{4} \nu_n^2 \|\taub_n\|^2)^2
.
$
Since
$
\thetab_n\pr{\bf M}_n\thetab_n
%=
%\thetab_n\pr  (\thetab_n\thetab_n\pr-\thetab_{n0}\thetab_{n0}\pr) \thetab_n
=
1 - (\thetab_{n0}\pr \thetab_n)^2
%=
% \nu_n^2 \|\taub_n\|^2 -\frac{1}{4} \nu_n^4 \|\taub_n\|^4
%=
% \nu_n^2 \|\taub_n\|^2 
% (1-\frac{1}{4} \nu_n^2 \|\taub_n\|^2)
$, 
this yields
\begin{eqnarray*}
\lefteqn{
{\rm tr}\Big[
 {\bf M}_n 
 ( a \thetab_n\thetab_n\pr + b ({\bf I}_{p_n}- \thetab_n\thetab_n\pr) )
 {\bf M}_n
 ( c\, \thetab_n\thetab_n\pr + d ({\bf I}_{p_n}- \thetab_n\thetab_n\pr) )
  \Big] 
}
  \\[2mm]
  & &
  \hspace{0mm} 
  =
{\rm tr}\Big[
 {\bf M}_n 
 ( b {\bf I}_{p_n} + (a-b) \thetab_n\thetab_n\pr )
 {\bf M}_n
 ( d {\bf I}_{p_n} + (c-d) \thetab_n\thetab_n\pr )
  \Big] 
  \\[2mm]
  & &
  \hspace{0mm} 
  =
bd\, {\rm tr}[{\bf M}_n^2 ]
+ b(c-d) \thetab_n\pr{\bf M}_n^2 \thetab_n
+ (a-b)d \thetab_n\pr{\bf M}_n^2 \thetab_n
+ (a-b)(c-d) (\thetab_n\pr{\bf M}_n \thetab_n)^2
  \\[2mm]
  & &
  \hspace{0mm} 
  =
2bd \thetab_n\pr{\bf M}_n \thetab_n
+ \{ b(c-d) + (a-b)d \} \thetab_n\pr{\bf M}_n \thetab_n
+ (a-b)(c-d) (\thetab_n\pr{\bf M}_n \thetab_n)^2
  \\[2mm]
  & &
  \hspace{0mm} 
  =
(ad+bc) (1 - (\thetab_{n0}\pr \thetab_n)^2)
+ (a-b)(c-d) (1 - (\thetab_{n0}\pr \thetab_n)^2)^2
\end{eqnarray*}
and
\begin{eqnarray*}
\lefteqn{
{\rm tr}\Big[
 {\bf M}_n^2 
 ( a \thetab_n\thetab_n\pr + b ({\bf I}_{p_n}- \thetab_n\thetab_n\pr) )
 {\bf M}_n
 ( c\, \thetab_n\thetab_n\pr + d ({\bf I}_{p_n}- \thetab_n\thetab_n\pr) )
  \Big] 
}
  \\[2mm]
  & &
  \hspace{0mm} 
  =
{\rm tr}\Big[
 {\bf M}_n^2 
 ( b {\bf I}_{p_n} + (a-b) \thetab_n\thetab_n\pr )
 {\bf M}_n^2
 ( d {\bf I}_{p_n} + (c-d) \thetab_n\thetab_n\pr )
  \Big] 
  \\[2mm]
  & &
  \hspace{0mm} 
  =
bd\, {\rm tr}[{\bf M}_n^4 ]
+ b(c-d) \thetab_n\pr{\bf M}_n^4 \thetab_n
+ (a-b)d \thetab_n\pr{\bf M}_n^4 \thetab_n
+ (a-b)(c-d) (\thetab_n\pr{\bf M}_n^2 \thetab_n)^2
  \\[2mm]
  & &
  \hspace{0mm} 
  =
2bd (\thetab_n\pr{\bf M}_n \thetab_n)^2
+ \{ b(c-d) + (a-b)d \} (\thetab_n\pr{\bf M}_n \thetab_n)^2
+ (a-b)(c-d) (\thetab_n\pr{\bf M}_n \thetab_n)^2
  \\[2mm]
  & &
  \hspace{0mm} 
  =
\{
ad+ bc
+ (a-b)(c-d) 
\}
(1 - (\thetab_{n0}\pr \thetab_n)^2)^2
  \\[2mm]
  & &
  \hspace{0mm} 
  =
(ac+bd)
(1 - (\thetab_{n0}\pr \thetab_n)^2)^2
,
\end{eqnarray*}
as was to be showed.
\cqfd
\vspace{3mm}

{\sc Proof of Theorem~\ref{TheoBillingH1}.}
All expectations and variances below are taken under~${\rm P}\n_{\thetab_n,F_n}$, with~$\thetab_n=\thetab_{n0}+\nu_n \taub_n$, and stochastic convergences are under the corresponding sequence of hypotheses. We have
$$
{\rm E}[{\bf X}_{ni}]=e_{n1} \thetab_n
\quad
\textrm{ and }
\quad
{\rm E}[{\bf X}_{ni}{\bf X}_{ni}\pr]
=
e_{n2} \thetab_n\thetab_n\pr +  \frac{f_{n2}}{p_n-1} ({\bf I}_{p_n}- \thetab_n\thetab_n\pr)
;
$$
see the proof of Lemma~B.3 in \cite{PaiVer17a}. Writing~${\bf W}_{ni}:=({\bf I}_{p_n}- \thetab_n\thetab_n\pr) \Xb_{ni}$, this implies that
\begin{equation}
\label{momW}	
{\rm E}[{\bf W}_{ni}]={\bf 0}
\quad
\textrm{ and }
\quad
{\rm E}[{\bf W}_{ni}{\bf W}_{ni}\pr]
%=
%{\rm E}[{\bf W}_{ni}{\bf X}_{ni}\pr]
%=
%{\rm E}[{\bf X}_{ni}{\bf W}_{ni}\pr]
=
  \frac{f_{n2}}{p_n-1} ({\bf I}_{p_n}- \thetab_n\thetab_n\pr)
.
\end{equation}
%\textcolor{red}{under~${\rm P}\n_{\thetab_n,F_n}$.}
Writing~${\bf M}_n=\thetab_n\thetab_n\pr-\thetab_{n0}\thetab_{n0}\pr$ as in Lemma~\ref{LemTrMn} and~$\Yb_{ni}:=({\bf I}_{p_n}- \thetab_{n0}\thetab_{n0}\pr) \Xb_{ni}$, we have ${\bf W}_{ni}=\Yb_{ni}-{\bf M}_n \Xb_{ni}$. This allows to decompose~$W^*_n$ as 
\begin{eqnarray*} 
W^*_n 
&\!\!:=\!\! & \frac{\sqrt{2(p_n-1)}}{nf_{n2}} \sum_{1 \leq i <j \leq n} 
\Yb_{ni}\pr \Yb_{nj} \\
%&=& \frac{\sqrt{2(p_n-1)}}{nf_{n2}} \sum_{1 \leq i <j \leq n} ({\bf W}_{ni}+{\bf M}_n \Xb_{ni})\pr ({\bf W}_{nj}+{\bf M}_n \Xb_{nj}) 
%\\
&\!\! =\!\! & \frac{\sqrt{2(p_n-1)}}{nf_{n2}} \sum_{1 \leq i <j \leq n}( {\bf W}_{ni}\pr {\bf W}_{nj}+\Xb_{ni}\pr{\bf M}_n{\bf W}_{nj}+{\bf W}_{ni}\pr{\bf M}_n\Xb_{nj}+ \Xb_{ni} \pr {\bf M}_n^2 \Xb_{nj})
\\[2mm]
&\!\!=:\!\!& 
W^*_{n0}
+
W^*_{na}
+
W^*_{nb}
+
W^*_{nc}
.
\end{eqnarray*}
From~(\ref{momW}), ${\rm E}[W^*_{na}]={\rm E}[W^*_{nb}]=0$. Now,
\begin{eqnarray*}
{\rm Var}[W^*_{na}]
&=& 
\frac{2(p_n-1)}{n^2f_{n2}^2}  \sum_{1 \leq i <j \leq n} \sum_{1 \leq r <s<n} 
{\rm E}[ \Xb_{ni}\pr{\bf M}_n{\bf W}_{nj} \Xb_{nr}\pr{\bf M}_n{\bf W}_{ns}]
\\[2mm]
%&=& 
%\frac{2(p_n-1)}{n^2f_{n2}^2}  \sum_{1 \leq i <j \leq n} \sum_{1 \leq r <s<n} 
%{\rm E}[{\rm tr}[\Xb_{ni}\pr{\bf M}_n{\bf W}_{nj} {\bf W}_{ns}\pr{\bf M}_n\Xb_{nr}]]
%\\[2mm]
&=& 
\frac{2(p_n-1)}{n^2f_{n2}^2}  \sum_{1 \leq i <j \leq n} \sum_{1 \leq r <s<n} 
{\rm tr}[{\rm E}[{\bf M}_n\Xb_{nr}\Xb_{ni}\pr{\bf M}_n{\bf W}_{nj} {\bf W}_{ns}\pr]]
\\[2mm]
&=& 
\frac{2(p_n-1)}{n^2f_{n2}^2}  \sum_{1 \leq i <j \leq n} \sum_{1 \leq r <s<n} 
c_{n,ijrs}
.
\end{eqnarray*}
Clearly, $c_{n,ijrs}=0$ if~$s\neq j$. Lemma~\ref{LemTrMn} entails that for~$s=j$ and $r \neq i$, we have
\begin{eqnarray*}
 c_{n,ijrs} 
&=& 
 {\rm tr}[{\bf M}_n{\rm E}[\Xb_{nr}\Xb_{ni}\pr] {\bf M}_n{\rm E}[{\bf W}_{nj}{\bf W}_{nj}\pr] ]
 \\[2mm]
&=& 
 {\rm tr}\Big[
 {\bf M}_n 
 (e_{n1}^2 \thetab_n\thetab_n\pr) 
 {\bf M}_n
 \Big(\frac{f_{n2}}{p_n-1} ({\bf I}_{p_n}- \thetab_n\thetab_n\pr)\Big)
  \Big] 
  \\[2mm]
%&=&
%(ad+bc) (1 - (\thetab_{n0}\pr \thetab_n)^2)
%+ (a-b)(c-d) ((1 - (\thetab_{n0}\pr \thetab_n)^2))^2
%\\[2mm]
&=& 
\frac{e_{n1}^2f_{n2}}{p_n-1}
 (1 - (\thetab_{n0}\pr \thetab_n)^2)
-
\frac{e_{n1}^2f_{n2}}{p_n-1}
 (1 - (\thetab_{n0}\pr \thetab_n)^2)^2
\end{eqnarray*}
and that, for~$s=j$ and~$r=i$, we have
\begin{eqnarray*}
 c_{n,ijrs} 
&=& 
 {\rm tr}[{\bf M}_n{\rm E}[\Xb_{ni}\Xb_{ni}\pr] {\bf M}_n{\rm E}[{\bf W}_{nj}{\bf W}_{nj}\pr] ]
 \\[2mm]
&=& 
  {\rm tr}\Big[
 {\bf M}_n 
 \Big( e_{n2} \thetab_n\thetab_n\pr + \frac{f_{n2}}{p_n-1} ({\bf I}_{p_n}- \thetab_n\thetab_n\pr) \Big)
 {\bf M}_n
 \Big(\frac{f_{n2}}{p_n-1} ({\bf I}_{p_n}- \thetab_n\thetab_n\pr)\Big)
  \Big] 
  \\[2mm]
%  & = &
%(ad+bc) (1 - (\thetab_{n0}\pr \thetab_n)^2)
%+ (a-b)(c-d) (1 - (\thetab_{n0}\pr \thetab_n)^2)^2
%  \\[2mm]
%  & = &
%\frac{e_{n2}f_{n2}}{p_n-1}
% (1 - (\thetab_{n0}\pr \thetab_n)^2)
%- 
%\bigg( e_{n2}-\frac{1-e_{n2}}{p_n-1} \bigg)
%\frac{1-e_{n2}}{p_n-1}
% (1 - (\thetab_{n0}\pr \thetab_n)^2)^2
%  \\[2mm]
&=& 
\frac{e_{n2}f_{n2}}{p_n-1}
 (1 - (\thetab_{n0}\pr \thetab_n)^2)
- 
\frac{(p_n e_{n2}-1)f_{n2}}{(p_n-1)^2}
 (1 - (\thetab_{n0}\pr \thetab_n)^2)^2
.
\end{eqnarray*}
We conclude that
\begin{eqnarray*}
\lefteqn{
{\rm Var}[W^*_{na}]
= 
\frac{2(p_n-1)}{n^2f_{n2}^2}  
\Bigg[ 
\frac{n(n-1)(n-2)}{3}
\bigg(
\frac{e_{n1}^2f_{n2}}{p_n-1}
 (1 - (\thetab_{n0}\pr \thetab_n)^2)
-
\frac{e_{n1}^2f_{n2}}{p_n-1}
 (1 - (\thetab_{n0}\pr \thetab_n)^2)^2
\bigg)
}
 \\[2mm]
& & 
\hspace{13mm} 
+
\frac{n(n-1)}{2}
\bigg(
\frac{e_{n2}f_{n2}}{p_n-1}
 (1 - (\thetab_{n0}\pr \thetab_n)^2)
- 
\frac{(p_n e_{n2}-1)f_{n2}}{(p_n-1)^2}
 (1 - (\thetab_{n0}\pr \thetab_n)^2)^2
\bigg)
\Bigg]
\\[2mm]
& & 
\hspace{-2mm} 
=
\frac{n-1}{3n}  
\Bigg[ 
\frac{2(n-2)e_{n1}^2+3e_{n2}}{f_{n2}}
 (1 - (\thetab_{n0}\pr \thetab_n)^2)
-
\bigg(
\frac{2(n-2)e_{n1}^2}{f_{n2}}
+
\frac{3(p_n e_{n2}-1)}{(p_n-1)f_{n2}}
\bigg)
 (1 - (\thetab_{n0}\pr \thetab_n)^2)^2
\Bigg]
.
\end{eqnarray*}
Since 
$
\thetab_n'\thetab_{n0}
=
(\thetab_{n0}+\nu_n\taub_n)'\thetab_{n0}
=
1+\nu_n (\taub_n'\thetab_{n0})
=
{\textstyle{1-\frac{1}{2} \nu_n^2 \|\taub_n\|^2}}
,
$
we have~$1 - (\thetab_{n0}\pr \thetab_n)^2=O(\nu_n^2)$, which, by using the fact that~$\nu_n=O(1)$, yields 
$$
{\rm Var}[W^*_{na}]
 =
\frac{np_n e_{n1}^2+p_n e_{n2}+1}{p_n f_{n2}}
\,
 O(\nu_n^2)
.
$$
The same computations provide~${\rm Var}[W^*_{nb}]={\rm Var}[W^*_{na}]$. 
Turning to~$W^*_{nc}$, 
%recall that
%$$
%W^*_{nc}
%=
%\frac{\sqrt{2(p_n-1)}}{nf_{n2}} \sum_{1 \leq i <j \leq n} \Xb_{ni} \pr {\bf M}_n^2 \Xb_{nj}
%.
%$$
we directly obtain
\begin{eqnarray*}
 {\rm E}[W^*_{nc}] 
&=&
  \frac{\sqrt{2(p_n-1)}}{nf_{n2}} \times
  \frac{n(n-1)}{2} \, e_{n1}^2 \thetab_n\pr{\bf M}_n^2\thetab_n 
\\[2mm] 
&=&
\frac{(n-1)(p_n-1)^{1/2} e_{n1}^2}{\sqrt{2}f_{n2}} 
\,
(1-(\thetab_n\pr\thetab_{n0})^2)
\\[2mm] 
&=&
\frac{np_n^{1/2}e_{n1}^2}{\sqrt{2}f_{n2}} 
\,
\nu_n^2 \|\taub_n\|^2  (1-{\textstyle{\frac{1}{4}}} \nu_n^2 \|\taub_n\|^2)(1+o(1))
%=
% \frac{np_n^{1/2}e_{n1}^2}{\sqrt{2}f_{n2}} 
%\,
%\nu_n^2 (1+o(1)+O(\nu_n^2)) \|\taub_n\|^2 
;
\end{eqnarray*}
see the proof of Lemma~\ref{LemTrMn}. As for the variance, 
\begin{eqnarray*}
{\rm Var}[W^*_{nc}]
&=& 
\frac{2(p_n-1)}{n^2f_{n2}^2}  
\sum_{1 \leq i <j \leq n} \sum_{1 \leq r <s\leq n} 
{\rm Cov}[ \Xb_{ni} \pr {\bf M}_n^2 \Xb_{nj}, \Xb_{nr} \pr {\bf M}_n^2 \Xb_{ns}]
\\[2mm]
&=& 
\frac{2(p_n-1)}{n^2f_{n2}^2}  
\sum_{1 \leq i <j \leq n} \sum_{1 \leq r <s\leq n} 
\Big(
{\rm E}[ \Xb_{ni} \pr {\bf M}_n^2 \Xb_{nj} \Xb_{nr} \pr {\bf M}_n^2 \Xb_{ns}]
-
(e_{n1}^2 \thetab_n\pr{\bf M}_n^2\thetab_n)^2
\Big)
\\[2mm]
&=& 
\frac{2(p_n-1)}{n^2f_{n2}^2}  
\sum_{1 \leq i <j \leq n} \sum_{1 \leq r <s\leq n} 
\Big(
{\rm tr}[{\rm E}[\Xb_{ni}\pr{\bf M}_n^2 {\bf X}_{nj} {\bf X}_{ns}\pr{\bf M}_n^2 \Xb_{nr}]]
-
e_{n1}^4 (1-(\thetab_n\pr\thetab_{n0})^2)^2
\Big)
\\[2mm]
&=& 
\frac{2(p_n-1)}{n^2f_{n2}^2}  
\sum_{1 \leq i <j \leq n} \sum_{1 \leq r <s\leq n} 
\Big(
d_{n,ijrs}
-
e_{n1}^4 (1-(\thetab_n\pr\thetab_{n0})^2)^2
\Big)
.
\end{eqnarray*}
We consider three cases. (1) If $i,j,r,s$ contain two pairs of equal indices (equivalently, if~$r=i$ and~$s=j$), then 
\begin{eqnarray*}
d_{n,ijrs}
%&=&
%{\rm tr}[{\bf M}_n^2 {\rm E}[\Xb_{nr}\Xb_{ni}\pr] {\bf M}_n^2 {\rm E}[{\bf X}_{nj} {\bf X}_{ns}\pr]]
%\\[2mm]
&=&
{\rm tr}[{\bf M}_n^2 {\rm E}[\Xb_{ni}\Xb_{ni}\pr] {\bf M}_n^2 {\rm E}[{\bf X}_{nj} {\bf X}_{nj}\pr]]
\\[2mm]
&=&
{\rm tr}\bigg[
 {\bf M}_n^2 
\Big( 
e_{n2}\thetab_n\thetab_n\pr+ \frac{f_{n2}}{p_n-1}({\bf I}_{p_n}- \thetab_n\thetab_n\pr)
\Big)
 {\bf M}_n^2 
\Big( 
e_{n2}\thetab_n\thetab_n\pr+ \frac{f_{n2}}{p_n-1}({\bf I}_{p_n}- \thetab_n\thetab_n\pr)
\Big)
\bigg]
\\[2mm]
&=&
\Big(
e_{n2}^2
+
\frac{f_{n2}^2}{(p_n-1)^2}
\Big)
(1-(\thetab_n\pr\thetab_{n0})^2)^2
.
\end{eqnarray*}
(2) If $i,j,r,s$ contain exactly one pair of equal indices, then
\begin{eqnarray*}
d_{n,ijrs}
%&=&
%{\rm tr}[{\bf M}_n^2 {\rm E}[\Xb_{nr}\Xb_{ni}\pr] {\bf M}_n^2 {\rm E}[{\bf X}_{nj} {\bf X}_{ns}\pr]]
%\\[2mm]
&=&
{\rm tr}[{\bf M}_n^2 {\rm E}[\Xb_{ni}\Xb_{ni}\pr] {\bf M}_n^2 {\rm E}[{\bf X}_{nj} {\bf X}_{ns}\pr]]
\\[2mm]
&=&
{\rm tr}\bigg[
 {\bf M}_n^2 
\Big( 
e_{n2}\thetab_n\thetab_n\pr+ \frac{f_{n2}}{p_n-1}({\bf I}_{p_n}- \thetab_n\thetab_n\pr)
\Big)
 {\bf M}_n^2 
\big( 
e_{n1}^2\thetab_n\thetab_n\pr
\big)
\bigg]
\\[2mm]
&=&
e_{n1}^2 e_{n2}
(1-(\thetab_n\pr\thetab_{n0})^2)^2
.
\end{eqnarray*}
(3) If the indices~$i,j,r,s$ are pairwise different, then 
\begin{eqnarray*}
d_{n,ijrs}
&=&
{\rm tr}[{\bf M}_n^2 {\rm E}[\Xb_{nr}\Xb_{ni}\pr] {\bf M}_n^2 {\rm E}[{\bf X}_{nj} {\bf X}_{ns}\pr]]
\\[2mm]
&=&
{\rm tr}\Big[
 {\bf M}_n^2 
\big( 
e_{n1}^2\thetab_n\thetab_n\pr
\big)
 {\bf M}_n^2 
\big( 
e_{n1}^2\thetab_n\thetab_n\pr
\big)
\Big]
\\[2mm]
&=&
e_{n1}^4
(1-(\thetab_n\pr\thetab_{n0})^2)^2
.
\end{eqnarray*}
Therefore,
\begin{eqnarray*}
\lefteqn{
{\rm Var}[W^*_{nc}]
= 
\frac{2(p_n-1)}{n^2f_{n2}^2}  
\bigg[
\frac{n(n-1)}{2}
\Big(
e_{n2}^2
+
\frac{f_{n2}^2}{(p_n-1)^2}
\Big)
+
n(n-1)(n-2)
e_{n1}^2 e_{n2}
}
\\[2mm]
& & 
\hspace{23mm} 
+
\frac{n(n-1)(n-2)(n-3)}{4}
\,
e_{n1}^4
-
\frac{n^2(n-1)^2}{4}
\,
e_{n1}^4
\bigg]
(1-(\thetab_n\pr\thetab_{n0})^2)^2
\\[2mm]
%%
%%
%&=& 
%\frac{(p_n-1)}{nf_{n2}^2}  
%\bigg[
%(n-1)
%\Big(
%e_{n2}^2
%+
%\frac{f_{n2}^2}{(p_n-1)^2}
%\Big)
%+
%2(n-1)(n-2)
%e_{n1}^2 e_{n2}
%\\[2mm]
%& & 
%\hspace{23mm} 
%-
%(n-1)(2n-3)
%e_{n1}^4
%\bigg]
%(1-(\thetab_n\pr\thetab_{n0})^2)^2
%\\[2mm]
%%
%%
%&=& 
%\bigg[
%\frac{(n-1)(p_n-1)}{nf_{n2}^2}  
%\Big(
%e_{n2}^2
%+
%\frac{f_{n2}^2}{(p_n-1)^2}
%\Big)
%+
%\frac{2(n-1)(n-2)
%(p_n-1)e_{n1}^2 e_{n2}
%}{nf_{n2}^2}  
%\\[2mm]
%& & 
%\hspace{23mm} 
%-
%\frac{(n-1)(2n-3)(p_n-1)e_{n1}^4}{nf_{n2}^2}  
%\bigg]
%(1-(\thetab_n\pr\thetab_{n0})^2)^2
%\\[2mm]
%%
%%
%&=& 
%\bigg[
%\frac{n-1}{n(p_n-1)}  
%+
%\frac{(n-1)(p_n-1)e_{n2}^2}{nf_{n2}^2}  
%+
%\frac{2(n-1)(n-2)
%(p_n-1)e_{n1}^2 e_{n2}
%}{nf_{n2}^2}  
%\\[2mm]
%& & 
%\hspace{23mm} 
%-
%\frac{(n-1)(2n-3)(p_n-1)e_{n1}^4}{nf_{n2}^2}  
%\bigg]
%(1-(\thetab_n\pr\thetab_{n0})^2)^2
%\\[2mm]
%%
%%
%&=& 
%\frac{n-1}{n}
%\bigg[
%\frac{(p_n-1)e_{n2}^2}{f_{n2}^2}  
%+
%\frac{2(n-2)
%(p_n-1)e_{n1}^2 e_{n2}
%}{f_{n2}^2}  
%-
%\frac{(2n-3)(p_n-1)e_{n1}^4}{f_{n2}^2}  
%\bigg]
%(1-(\thetab_n\pr\thetab_{n0})^2)^2 
%+
%o(1)
%\\[2mm]
%
%
%& &
%\hspace{1mm} 
%= 
%\frac{(n-1)(p_n-1)\tilde{e}_{n2}(e_{n2}+(2n-3)e_{n1}^2)}{nf_{n2}^2}
%(1-(\thetab_n\pr\thetab_{n0})^2)^2
%+
%o(1)
%\\[2mm]
%
%
& &
\hspace{1mm} 
= 
\frac{(n-1)(p_n-1)\tilde{e}_{n2}(\tilde{e}_{n2}+2(n-1)e_{n1}^2)}{nf_{n2}^2}
(1-(\thetab_n\pr\thetab_{n0})^2)^2
+
o(1)
.
\end{eqnarray*}
This finally yields
$$
{\rm Var}[W^*_{nc}]
 =
\frac{p_n \tilde{e}_{n2}^2+np_n e_{n1}^2  \tilde{e}_{n2}}{f_{n2}^2}
\,
 O(\nu_n^4)
.
$$

%& &
%\hspace{1mm} 
%= 
%\frac{n-1}{n}
%\bigg[
%\frac{2n
%(p_n-1)e_{n1}^2 \tilde{e}_{n2}
%}{f_{n2}^2}  
%+
%\frac{(p_n-1)(e_{n2}^2-4e_{n1}^2e_{n2}+6e_{n1}^4)}{f_{n2}^2}  
%\bigg]
%(1-(\thetab_n\pr\thetab_{n0})^2)^2
%+
%o(1)

%%%%%%%%%%%%%%%%%%%%%%

Summarizing,
$
W^*_{n}
=
W^*_{n0}+W^*_{na}+W^*_{nb}+W^*_{nc}
,
$
where~$W^*_{n0}$ is asymptotically standard normal (see Theorem~3.1 from \cite{Leyetal2015}), 
$$
{\rm E}[W^*_{na}]={\rm E}[W^*_{nb}]=0
, 
\quad
 {\rm E}[W^*_{nc}] 
= 
\frac{np_n^{1/2}e_{n1}^2}{\sqrt{2}f_{n2}} 
\,
\nu_n^2 \|\taub_n\|^2  (1-{\textstyle{\frac{1}{4}}} \nu_n^2 \|\taub_n\|^2)(1+o(1))
,
$$
$$
{\rm Var}[W^*_{na}]
=
{\rm Var}[W^*_{nb}]
=
\frac{np_n e_{n1}^2+p_n e_{n2}+1}{p_n f_{n2}}
\,
 O(\nu_n^2)
$$
and
$$
%\
%\textrm{ and }
%\
{\rm Var}[W^*_{nc}]
=
\frac{p_n \tilde{e}_{n2}^2+np_n e_{n1}^2  \tilde{e}_{n2}}{f_{n2}^2}
\,
 O(\nu_n^4)
.
$$  

We can now consider the several cases of the theorem. In cases~(i)--(iii), the sequence~$(\nu_n)$ involved, namely~$\nu_n=\sqrt{f_{n2}}/(\sqrt{n}p_n^{1/4}e_{n1})$, is~$o(1)$, so that~${\rm E}[W^*_{nc}]=t^2/\sqrt{2}+o(1)$. In all three cases, one checks that ${\rm Var}[W^*_{n\ell}]=o(1)$ for~$\ell=a,b,c$ (note that in cases~(ii)--(iii), the fact that $e_{n2}\leq e_{n1}$ implies that both~$e_{n2}$ and~$\tilde{e}_{n2}$ are~$o(1)$), which establishes that 
$
	W^*_n 
	\stackrel{\mathcal{D}}{\longrightarrow}
	\mathcal{N}\big(\frac{t^2}{\sqrt{2}},1\big)
$
	under~${\rm P}\n_{\thetab_{n0}+\nu_n \taub_n,F_n}$. 
In case~(iv), we have, with~$\nu_n=1$, 
${\rm E}[W^*_{nc}]=
	\frac{\xi^2t^2}{\sqrt{2}}
	\big(1-\frac{t^2}{4} \big)
+o(1)$. 
Since~$\sqrt{p_n} e_{n2}=o(1)$ 
\vspace{-.6mm}
 by assumption, one can check that~${\rm Var}[W^*_{n\ell}]=o(1)$ for~$\ell=a,b,c$, which yields
$
	W^*_n 
	\stackrel{\mathcal{D}}{\longrightarrow}
	\mathcal{N}\big(
	\frac{\xi^2t^2}{\sqrt{2}}
	\big(1-\frac{t^2}{4} \big)
	,1
	\big)
$
under~${\rm P}\n_{\thetab_{n0}+\nu_n \taub_n,F_n}$, as was to be showed. 
%with~$
%\nu_n
%=
%f_{n2}^{1/2}/(n^{1/2}p_n^{1/2}e_{n1})
%(=o(1))
%$, we have ${\rm Var}[W^*_{n\ell}]=o(1)$ for~$\ell=a,b,c$ and
%${\rm E}[W^*_{nc}]=o(1)$, which implies that~$W^*_n$ is asymptotically standard normal. In the same regime, if
%
Finally, in case~(v),
\vspace{.3mm}
  still with~\mbox{$\nu_n=1$}, we have~${\rm E}[W^*_{nc}]=o(1)$. One can again check that~${\rm Var}[W^*_{n\ell}]=o(1)$ for~$\ell=a,b,c$, which yields that~$W^*_n$ is asymptotically standard normal. This establishes the result. 
\cqfd
\vspace{3mm}

We turn to the proof of Theorem~\ref{TheoGeneralSlutzky}, that will make use of the following lemma.

\begin{lemma}
\label{LemMomentsForGeneralSlutzky}
Under~${\rm P}\n_{\thetab_n,F_n}$,  
$$
{\rm E}[(\Xb_{n1}'\thetab_{n0})^2]
=
e_{n2} (\thetab_{n0}\pr\thetab_n)^2 
+
\frac{f_{n2}}{p_n-1} (1-(\thetab_{n0}\pr\thetab_n)^2 )
$$
and
$$
{\rm E}[(\Xb_{n1}'\thetab_{n0})^4]
=
e_{n4}
(\thetab_n'\thetab_{n0})^4
+
\frac{6(e_{n2}-e_{n4})}{p_n-1}
(\thetab_n'\thetab_{n0})^2
 (1-(\thetab_{n0}\pr\thetab_n)^2 )
+ 
\frac{3f_{n4}}{p_n^2-1} (1-(\thetab_{n0}\pr\thetab_n)^2 )^2 
.
$$
\end{lemma}

{\sc Proof of Lemma~\ref{LemMomentsForGeneralSlutzky}.}
All computations in this proof are performed under~${\rm P}\n_{\thetab_n,F_n}$, which leads us to consider the tangent-decomposition~$\Xb_{n1}=U_{n1}\thetab_n+V_{n1}\Sb_{n1}$ of~$\Xb_{n1}$ with respect to~$\thetab_n$. Since~$\Xb_{n1}$ is rotationally symmetric with respect to~$\thetab_n$, $\Sb_{n1}$ is equal in distribution to~$\Gamb_{\thetab_n}\Ub_n$, where~$\Ub_n$ is uniformly distributed over the unit sphere~$\mathcal{S}^{p_n-2}$ in~$\R^{p_n-1}$ and where~$\Gamb_{\thetab_n}$ is an arbitrary $p_n\times (p_n-1)$ matrix whose columns form an orthonormal basis of the orthogonal complement of~$\thetab_n$ in~$\R^{p_n}$ (so that~$\Gamb_{\thetab_n}'\Gamb_{\thetab_n}=\mathbf{I}_{p_n-1}$ and $\Gamb_{\thetab_n}\Gamb_{\thetab_n}'=\mathbf{I}_{p_n}-\thetab_n\thetab_n'$).    In particular, 
$$
{\rm E}[\Sb_{n1}]=\mathbf{0}
\quad\textrm{and}
\quad
{\rm E}[\Sb_{n1}\Sb_{n1}']
=
\frac{1}{p_n-1} 
\,
({\bf I}_{p_n}- \thetab_n\thetab_n\pr)
.
$$
This readily yields 
\begin{eqnarray*}
{\rm E}[(\Xb_{n1}'\thetab_{n0})^2]
&=&
{\rm E}[(U_{n1}\thetab_n'\thetab_{n0}+V_{n1} \Sb_{n1}'\thetab_{n0})^2]
\nonumber
\\[2mm]
&=&
{\rm E}[U_{n1}^2] (\thetab_n'\thetab_{n0})^2
+ 2 {\rm E}[U_{n1} V_{n1}] {\rm E}[\Sb_{n1}'\thetab_{n0}](\thetab_n'\thetab_{n0})
+ {\rm E}[V_{n1}^2] \thetab_{n0}'{\rm E}[\Sb_{n1}\Sb_{n1}'] \thetab_{n0}
\nonumber
\\[2mm]
&=&
e_{n2} (\thetab_{n0}\pr\thetab_n)^2 
+
\frac{f_{n2}}{p_n-1} (1-(\thetab_{n0}\pr\thetab_n)^2 )
.
\end{eqnarray*}
Using the identity~$U_{n1}^2 V_{n1}^2=U_{n1}^2-U_{n1}^4$, we obtain similarly
\begin{eqnarray}
\lefteqn{
\hspace{-3mm} 
{\rm E}[(\Xb_{n1}'\thetab_{n0})^4]
=
{\rm E}[(U_{n1}\thetab_n'\thetab_{n0}+V_{n1} \Sb_{n1}'\thetab_{n0})^4]
}
\nonumber
\\[2mm]
& & 
\hspace{3mm} 
=
e_{n4}
(\thetab_n'\thetab_{n0})^4
+
6 {\rm E}[U_{n1}^2 V_{n1}^2 (\Sb'_{n1}\thetab_{n0})^2]
(\thetab_n'\thetab_{n0})^2
+
f_{n4} {\rm E}[ (\Sb'_{n1}\thetab_{n0})^4]
\nonumber
\\[2mm]
& & 
\hspace{3mm} 
=
e_{n4}
(\thetab_n'\thetab_{n0})^4
+
6 (e_{n2}-e_{n4}) 
\,
( \thetab_{n0}'{\rm E}[\Sb_{n1}\Sb_{n1}'] \thetab_{n0})
(\thetab_n'\thetab_{n0})^2
+
f_{n4} {\rm E}[ (\Sb'_{n1}\thetab_{n0})^4]
\nonumber
\\[2mm]
& & 
\hspace{3mm} 
=
e_{n4}
(\thetab_n'\thetab_{n0})^4
+
\frac{6(e_{n2}-e_{n4})}{p_n-1}
(\thetab_n'\thetab_{n0})^2
 (1-(\thetab_{n0}\pr\thetab_n)^2 )
+ 
f_{n4} {\rm E}[ (\Sb'_{n1}\thetab_{n0})^4]
%\frac{3f_{n4}}{p_n^2-1}
% (1-(\thetab_{n0}\pr\thetab_n)^2 )^2 
.
\label{almostlem}
\end{eqnarray}
%It therefore only remains to show that
%\begin{equation}
%	\label{toshowforlemmaslutzky}
%{\rm E}[ (\Sb'_{n1}\thetab_{n0})^4]
%=
%\frac{3}{p_n^2-1} (1-(\thetab_{n0}\pr\thetab_n)^2 )^2 
%.
%\end{equation}
Standard formulas for the Kronecker product yield
\begin{eqnarray*}
\lefteqn{
	{\rm E}[(\Sb'_{n1}\thetab_{n0})^4]
=
{\rm E}[(\thetab_{n0}\pr \Sb_{n1}\Sb_{n1}\pr \thetab_{n0})^2
%\thetab_{n0}\pr \Sb_{n1}\Sb_{n1}\pr \thetab_{n0}
]
=
(\thetab_{n0}\otimes \thetab_{n0})\pr 
{\rm E}[{\rm vec}(\Sb_{n1}\Sb_{n1}\pr){\rm vec}'(\Sb_{n1}\Sb_{n1}\pr)]
(\thetab_{n0}\otimes \thetab_{n0})
}
\\[3mm]
& & 
\hspace{1mm} 
=
(\thetab_{n0}\otimes \thetab_{n0})\pr 
(\Gamb_{\thetab_n}\otimes \Gamb_{\thetab_n})
{\rm E}[
{\rm vec}(\Ub_{n}\Ub_{n}\pr)
{\rm vec}'(\Ub_{n}\Ub_{n}\pr)
]
(\Gamb_{\thetab_n}\pr\otimes \Gamb_{\thetab_n}\pr)
(\thetab_{n0}\otimes \thetab_{n0})
\\[2mm]
& & 
\hspace{1mm} 
=
\frac{1}{p_n^2-1}
(\thetab_{n0}\otimes \thetab_{n0})\pr 
(\Gamb_{\thetab_n}\otimes \Gamb_{\thetab_n})
\big(
\mathbf{I}_{(p_n-1)^2}
+
\mathbf{K}_{p_n-1}
+
\mathbf{J}_{p_n-1}
\big)
(\Gamb_{\thetab_n}\pr\otimes \Gamb_{\thetab_n}\pr)
(\thetab_{n0}\otimes \thetab_{n0})
,
\end{eqnarray*}
where~$\mathbf{K}_{\ell}$ is the $\ell\times\ell$ commutation matrix and where we let~$\mathbf{J}_{\ell}=({\rm vec}\,\mathbf{I}_\ell)({\rm vec}\,\mathbf{I}_\ell)'$; see \cite{Tyl1987}, page~244. Using the fact that~$\mathbf{K}_\ell(\mathbf{A}\otimes \mathbf{B})=(\mathbf{A}\otimes \mathbf{B})\mathbf{K}_{\ell'}$ for~$\ell\times \ell'$ matrices~$\Ab$ and~$\Bb$, along with the identity~$\mathbf{K}_1=1$, we obtain
\begin{eqnarray*}
	{\rm E}[(\Sb'_{n1}\thetab_{n0})^4]
&=&
\frac{2}{p_n^2-1}
(\thetab_{n0}\otimes \thetab_{n0})\pr 
(\Gamb_{\thetab_n}\otimes \Gamb_{\thetab_n})
(\Gamb_{\thetab_n}\pr\otimes \Gamb_{\thetab_n}\pr)
(\thetab_{n0}\otimes \thetab_{n0})
\\[2mm]
& & 
\hspace{5mm} 
+
\,
\frac{1}{p_n^2-1}
(\thetab_{n0}\otimes \thetab_{n0})\pr 
{\rm vec}(\Gamb_{\thetab_n}\Gamb_{\thetab_n}\pr) 
{\rm vec}'(\Gamb_{\thetab_n}\Gamb_{\thetab_n}\pr) 
(\thetab_{n0}\otimes \thetab_{n0})
\\[2mm]
&=&
\frac{3}{p_n^2-1}
 (1-(\thetab_{n0}\pr\thetab_n)^2 )^2 
.
\end{eqnarray*}
Plugging this in~(\ref{almostlem}) provides the result.
\cqfd
\vspace{3mm}
 
{\sc Proof of Theorem~\ref{TheoGeneralSlutzky}.}
Fix a sequence of hypotheses~${\rm P}\n_{\thetab_{n0}+\nu_n \taub_n,F_n}$ associated with a given regime~(i) to~(v) in Theorem~\ref{TheoBillingH1}. Throughout the proof, stochastic convergences, expectations and variances refer to this sequence of hypotheses. In view of the decomposition
$
\widetilde{W}_n
-
W^*_n
=
L_n^{-1}(1-L_n) W^*_n
$
from~(\ref{differencedeW}), it is sufficient to show that~$L_n$ converges to one in quadratic mean (note indeed that Theorem~\ref{TheoBillingH1} indeed implies that~$W^*_n$ is~$O_{\rm P}(1)$). In order to do so, write\begin{eqnarray*}
{\rm E}
\big[
(L_{n}-1)^2
\big]
&=&
\frac{1}{f_{n2}^2}
\,
{\rm E}
\Bigg[
\Bigg(
 f_{n2}
-
\bigg[ \frac{1}{n} \sum_{i=1}^{n} V_{ni}^2\bigg]
\Bigg)^2
\Bigg]
=
\frac{1}{f_{n2}^2}
\,
{\rm E}
\Bigg[
\Bigg(
\bigg[ \frac{1}{n} \sum_{i=1}^{n} (\Xb_{ni}'\thetab_{n0})^2\bigg]
- e_{n2}
\Bigg)^2
\Bigg]
\\[2mm]
&=&
\frac{1}{f_{n2}^2}
\,
{\rm E}
\Bigg[
\Bigg(
\bigg[ \frac{1}{n} \sum_{i=1}^{n} (\Xb_{ni}'\thetab_{n0})^2\bigg]
- {\rm E}[(\Xb_{n1}'\thetab_{n0})^2]
+ {\rm E}[(\Xb_{n1}'\thetab_{n0})^2]
- e_{n2}
\Bigg)^2
\Bigg]
\\[2mm]
&\!\!\!\leq\!\!\!&
\frac{2}{f_{n2}^2}
\,
{\rm Var}
\Bigg[
\frac{1}{n} \sum_{i=1}^{n} (\Xb_{ni}'\thetab_{n0})^2
\Bigg]
+
\frac{2}{f_{n2}^2}
\big( {\rm E}[(\Xb_{n1}'\thetab_{n0})^2] - e_{n2}\big)^2
\,
\\[2mm]
&\!\!\!\leq\!\!\!&
\frac{2}{nf_{n2}^2}
\,
\Big(
{\rm E}
\big[
(\Xb_{n1}'\thetab_{n0})^4
\big]
-
\big(
{\rm E}
\big[
(\Xb_{n1}'\thetab_{n0})^2
\big]
\big)^2
\Big)
+
\frac{2}{f_{n2}^2}
\big( {\rm E}[(\Xb_{n1}'\thetab_{n0})^2] - e_{n2}\big)^2
\,
\\[3mm]
&\!\!\!=:\!\!\!&
2T_{na}
+
2T_{nb}
,
\end{eqnarray*}
say. Since~$f_{n2}=1-e_{n2}$, Lemma~\ref{LemMomentsForGeneralSlutzky} provides
$$
{\rm E}[(\Xb_{n1}'\thetab_{n0})^2] - e_{n2}
%=
%e_{n2} ((\thetab_{n0}\pr\thetab_n)^2-1) 
%+
%\frac{1-e_{n2}}{p_n-1} (1-(\thetab_{n0}\pr\thetab_n)^2 )
%$$
%$$
=
\bigg(
\frac{f_{n2}}{p_n-1}-e_{n2}\bigg)
 (1-(\thetab_{n0}\pr\thetab_n)^2 )
=
\frac{1-p_ne_{n2}}{p_n-1}
 (1-(\thetab_{n0}\pr\thetab_n)^2 )
 ,
$$
which yields
$$
T_{nb}
=
\frac{(1-p_ne_{n2})^2}{(p_n-1)^2f_{n2}^2}
 (1-(\thetab_{n0}\pr\thetab_n)^2 )^2
=
\frac{1+p_n^2e_{n2}^2}{p_n^2 f_{n2}^2}
\,
O(\nu_n^4)
.
$$
In each of the regimes considered in Theorem~\ref{TheoBillingH1}, we thus obtain that~$T_{nb}=o(1)$, irrespective of the fact that~$\sqrt{p_n}e_{n2}=o(1)$ or not. Turning to~$T_{na}$, Lemma~\ref{LemMomentsForGeneralSlutzky} yields
\begin{eqnarray*}
nf_{n2}^2 T_{na}
&=&
e_{n4}
(\thetab_n'\thetab_{n0})^4
+
\frac{6(e_{n2}-e_{n4})}{p_n-1}
(\thetab_n'\thetab_{n0})^2
 (1-(\thetab_{n0}\pr\thetab_n)^2 )
+
\frac{3f_{n4}}{p_n^2-1}
 (1-(\thetab_{n0}\pr\thetab_n)^2 )^2
\\[2mm]
& & 
\hspace{-6mm} 
-
\bigg(
e_{n2}^2 (\thetab_{n0}\pr\thetab_n)^4 
+
\frac{2e_{n2}f_{n2}}{p_n-1} (\thetab_{n0}\pr\thetab_n)^2 
 (1-(\thetab_{n0}\pr\thetab_n)^2 )
+
\frac{f_{n2}^2}{(p_n-1)^2} (1-(\thetab_{n0}\pr\thetab_n)^2 )^2
\bigg)
\\[2mm]
&=&
(e_{n4}-e_{n2}^2)
(\thetab_n'\thetab_{n0})^4
+
\bigg(
\frac{6(e_{n2}-e_{n4})}{p_n-1}
-
\frac{2e_{n2} f_{n2}}{p_n-1}
\bigg)
 (\thetab_{n0}\pr\thetab_n)^2
 (1-(\thetab_{n0}\pr\thetab_n)^2 )
\\[2mm]
& & 
\hspace{3mm} 
+
\bigg(
\frac{3f_{n4}}{p_n^2-1}
-
\frac{f_{n2}^2}{(p_n-1)^2} 
\bigg)
(1-(\thetab_{n0}\pr\thetab_n)^2 )^2
.
\end{eqnarray*}
Using the facts that 
$
e_{n4}-e^2_{n2}
=
{\rm Var}[U_{n1}^2]
%=
%{\rm Var}[1-U_{n1}^2]
=
{\rm Var}[V_{n1}^2]
\leq
{\rm E}[V_{n1}^4]
=
f_{n4}
$
and that
$
e_{n2}-e_{n4}
=
{\rm E}[U_{n1}^2(1-U_{n1}^2)]
\leq
{\rm E}[1-U_{n1}^2]
%=
%1-e_{n2}
=
f_{n2}
,
$
we then obtain
\begin{eqnarray*}
T_{na}
&\!\!\leq \!\!&
\frac{f_{n4}}{nf_{n2}^2}
\,
(\thetab_{n0}\pr\thetab_n)^2
+
\frac{6-2e_{n2}}{n(p_n-1)f_{n2}}
 (\thetab_{n0}\pr\thetab_n)^2
 (1-(\thetab_{n0}\pr\thetab_n)^2 )
\\[2mm]
& & 
\hspace{23mm} 
+
\bigg(
\frac{3f_{n4}}{n(p_n^2-1)f_{n2}^2}
-
\frac{1}{n(p_n-1)^2} 
\bigg)
(1-(\thetab_{n0}\pr\thetab_n)^2 )^2
\\[2mm]
&\!\!=\!\! &
o(1)
+
\frac{1}{np_nf_{n2}} \, O(\nu_n^2)
+
o(\nu_n^4)
=
o(1) 
+
\frac{1}{np_nf_{n2}}
\,
O(\nu_n^2)
.
\end{eqnarray*}
Trivially, we then have~$T_{na}=o(1)$ in each of the regime considered in Theorem~\ref{TheoBillingH1}, still irrespective of the fact that~$\sqrt{p_n}e_{n2}=o(1)$ or not. This establishes the result. 
\cqfd
\vspace{3mm}
 
%%%%%%%%%%%%%%%%%%%%%%%%%%%%%%%%%%%%%%%%%%%%%%%%%%%%%%%

%\begin{acknowledgements}
%\end{acknowledgements}

%%%%%%%%%Bibliography

\bibliographystyle{spmpsci}      % mathematics and physical sciences
\bibliography{Paper2.bib}   % name your BibTeX data base

\end{document}